\input amstex 
\documentstyle{amsppt}
\magnification 1200
\NoBlackBoxes

\topmatter
\title 
Modified proof of 
a local analogue of the Grothendieck
Conjecture 
\endtitle 
\rightheadtext{Analogue of Grothendieck Conjecture}

\author Victor Abrashkin 
\endauthor

\email victor.abrashkin@durham.ac.uk 
\endemail 
\address Maths Dept., Durham University, Sci. Laboratories, 
South Rd., Durham, DH1 3LE, U.K. 
\endaddress 
\abstract 
A local analogue of the Grothendieck Conjecture is an equivalence
of the category of complete discrete valuation fields $K$ 
with finite residue fields of characteristic $p\ne 0$ 
and the
category of absolute Galois groups of fields $K$ together with their 
ramification filtrations. The case
of characteristic 0 fields $K$ was considered by Mochizuki several years
ago. 
Then the author 
proved it by different method if $p>2$ (but $\operatorname{char}K=0$ 
or $p$). This paper represents a modified approach: it 
covers the case $p=2$, contains 
considerable technical simplifications and replaces the Galois group 
of $K$ by its maximal pro-$p$-quotient. Special attention is 
paid to the procedure of recovering  
field isomorphisms coming from isomorphisms of Galois groups, which
are 
compatible with corresponding ramification filtrations. 

\ \ 
\newline 
{\smc R\'esum\'e.} 
Un analogue local de la conjecture de Grothendieck est une \' equivalence
entre la cat\'egorie des corps $K$ complets pour une valuation discr\` ete
\` a corps r\'esiduels finis
de caract\' eristique $p\ne 0$, et la cat\'egorie des groupes
galoisiens
absolus de corps $K$ munis de la filtration de ramification.
Le cas des corps de caract\' eristique 
0 a \' et\' e consid\' er\' e par Mochizuki il y a quelques ann\' ees.
Par la suite, le pr\' esent auteur
a demontr\' e l'\'equivalence par une m\' ethode diff\'erente si $p>2$
(mais $\operatorname{char}K=0$ or $p$). Dans l'article pr\' esent\' e
ici, une modification de l'approche pr\' ec\' edente
est envisag\' ee: elle couvre le cas $p=2$, contient des 
simplifications consid\' erables et remplace le group galoisien absolu
de $K$ par son
pro-$p$-quotient maximal. Une attention particuli\` ere 
est accorde\' e au 
proc\' ed\' e de reconstruction d'isomorphisme de corps obtenu 
a partir 
d'isomorphisme de groupes du Galois qui sont compatibles avec les 
filtrations de ramification correspondantes. 
\endabstract 

\endtopmatter

\def\adm{\operatorname{adm}}
\def\tr{\operatorname{tr}}
\def\Tr{\operatorname{Tr}}

\def\Ker{\operatorname{Ker}}
\def\ab{\operatorname{ab}}
\def\Gal{\operatorname{Gal}}
\def\sep{\operatorname{sep}}
\def\Hom{\operatorname{Hom}}
\def\id{\operatorname{id}}

\def\max{\operatorname{max}}

\def\adm{\operatorname{adm}}

\def\ur{\operatorname{ur}}

\def\Aut{\operatorname{Aut}}

\def\Iso{\operatorname{Iso}}

\def\char{\operatorname{char}}
\def\m{\operatorname{m}}

\def\char{\operatorname{char}}

\def\Fr{\operatorname{Fr}}
\def\char{\operatorname{char}}

\def\Tr{\operatorname{Tr}}
\def\Iso{\operatorname{Iso}}
\def\an{\operatorname{an}}

\document 

\subhead 0. Introduction 
\endsubhead 
\medskip

Throughout all this paper $p$ is a prime number. 
If $E$ is a complete discrete valuation field 
then we shall assume that 
its residue field has characteristic $p$, $E$ is considered as a
subfield of its fixed separable closure $E_{\sep }$, 
$\Gamma _E=\Gal (E_{\sep }/E)$. 
$E(p)$ will denote the maximal $p$-extension 
of $E$ in $E_{\sep }$ and $\Gamma _E(p)=\Gal (E(p)/E)$. 

Assume that $E,E'$ 
are complete discrete valuation fields with 
finite residue fields and there is a continuous 
field isomorphism $\mu :E\longrightarrow E'$. Then $\mu $ 
can be extended to a field isomorphism 
$\bar\mu :E(p)\longrightarrow E'(p)$. 
The correspondence $\tau\mapsto \bar\mu ^{-1}\tau\bar\mu $  
(cf. the agreement about compositions of morphisms 
in the end of this Introduction)  
defines a continuous group isomorphism 
$\bar\mu ^*:\Gamma _E(p)\longrightarrow\Gamma _{E'}(p)$ such that 
for any $v\geqslant 0$, 
$\bar\mu ^*(\Gamma _E(p)^{(v)})=\Gamma _{E'}(p)^{(v)}$. Here 
$\Gamma _E(p)^{(v)}$ is the ramification subgroup of 
$\Gamma _E(p)$ in the upper numbering.

The principal result of this paper is the following theorem. 

\proclaim{Theorem A} Suppose $E,E'$ are complete discrete valuation
fields with finite residue fields and there is 
a continuous group isomorphism 
$g:\Gamma _E(p)\longrightarrow\Gamma _{E'}(p)$ such that 
for any $v\geqslant 0$, $g(\Gamma _E(p)^{(v)})=\Gamma
_{E'}(p)^{(v)}$. Then there is 
a continuous field isomorphism $\bar\mu :E(p)\longrightarrow
E'(p)$ such that $\bar\mu (E)=E'$ and 
$g=\bar\mu ^*$. 
\endproclaim

This theorem implies easily a corresponding statement, where 
the maximal $p$-extensions $E(p)$ and $E'(p)$ and 
their Galois groups $\Gamma _E(p)$ and $\Gamma _{E'}(p)$ are 
replaced, respectively, by the separable closures $E_{\sep }$ and $E_{\sep }'$ and 
the Galois groups $\Gamma _E$ and $\Gamma _{E'}$. 
Such a statement is known as  a   
local analogue of the Grothendieck Conjecture. 
Mochizuki [Mo] proved it for local fields 
of characteristic 0. His method is based on an elegant
application of Hodge-Tate theory. Under the restriction 
$p>2$ the case of local fields of arbitrary characteristic 
was proved by another method by the author [Ab3]. This proof is based 
on an explicit description of the ramification subgroups 
$\Gamma _K(p)^{(v)}$ modulo the subgroup 
$C_3(\Gamma _K(p))$ of commutators 
of order $\geqslant 3$ in $\Gamma _K(p)$, where 
$K=k((t))$, and $k$ is a finite field of characteristic 
$p>2$. The restriction $p\ne 2$ appears because the proof  
uses the equivalence of the
category of $p$-groups and of Lie $\Bbb Z_p$-algebras of nilpotent class
2, 
which holds only under the assumption $p>2$. 

The statement of Theorem A is free from 
the restriction $p\ne 2$. Its proof follows mainly the 
strategy from [Ab3] but there are 
several essential changes. 

Firstly, instead of working with the ramification subgroups 
$\Gamma _K(p)^{(v)}$, $v\geqslant 0$, we fix the simplest possible embedding of 
$\Gamma _K(p)$ into its Magnus's algebra $\Cal A$ and 
study the induced fitration by the ideals $\Cal A^{(v)}$, $v\geqslant
0$, of $\Cal A$. As a result, we obtain an explicit description of the
ideals $\Cal A^{(v)}\operatorname{mod}\Cal J^3$, where 
$\Cal J$ is the augmentation ideal in $\Cal A$. This corresponds to the
description of the groups $\Gamma
_K(p)^{(v)}\operatorname{mod}C_3(\Gamma _K(p))$ in [Ab1] 
but it is easier to obtain and it works for all prime numbers $p$ 
including $p=2$. 

Secondly, any continuous group automorphism of $\Gamma _K(p)$
which is compatible with the ramification filtration induces 
a continuous algebra automorphism $f$ of $\Cal A$ such that 
for any $v\geqslant 0$, $f(\Cal A^{(v)})=\Cal A^{(v)}$. 
Similarly to [Ab3], the conditions 
$f(\Cal A^{(v)})\operatorname{mod}\Cal J^3=\Cal
A^{(v)}\operatorname{mod}\Cal J^3$ imply non-trivial properties of the restriction
of the original automorphism of $\Gamma _K(p)$ 
to the inertia subgroup $I_K(p)^{\ab }$ of the Galois group 
of the maximal abelian extension of $K$. 
These properties are studied in detail in this paper.  
This allows us to give a more detailed 
and effective version of 
the final stage of the proof of the local analogue of the Grothendieck
Conjecture even in the case $p\ne 2$. In particular, 
this clarifies why it holds with the absolute Galois groups 
replaced by the Galois groups of maximal $p$-extensions.

The methods of this paper can be helpful for understanding the  
relations between fields and their Galois groups in
the context of the global Grothendieck Conjecture. 
For example, suppose $F$ is an algebraic number field, $\bar F$ is its 
algebraic closure, $\Gamma _F=\Gal (\bar F/F)$, 
$\wp $ is a prime divisor in $F$, 
$\bar\wp $ is its extension to $\Bar F$ and $F_{\wp }$, 
$\bar F_{\bar\wp }$ are the corresponding completions of $F$ and 
$\bar F$, respectively. Then $\Gamma _{F,\bar\wp }=\Gal (\bar F_{\bar\wp
}/F_{\wp })\subset\Gamma _F$ is the decomposition group of $\bar\wp
$. 
Suppose $F$ is Galois over $\Bbb Q$ and 
$g_{\wp }:\Gamma _{F,\bar\wp }\longrightarrow\Gamma _{F,\bar\wp }$ 
is a continuous group automorphism which is compatible with the 
ramification filtration on $\Gamma _{F,\bar\wp }$. 
By the local analogue of the Grothendieck Conjecture, $g_{\wp }$ 
is induced by a field automorphism $\bar\mu _{\wp }:
\bar F_{\bar\wp }\longrightarrow\bar F_{\bar\wp }$ such that 
$\bar\mu :=\bar\mu _{\bar\wp }|_{\bar F}$ maps $\bar F$ to $\bar F$ 
(because $\bar\mu (\Bbb Q)=\Bbb Q$), 
and, therefore, $F$ to $F$ (because $F$ is Galois over $\Bbb Q$). 
So, $\bar\mu $
induces a group automorphism $g$ of $\Gamma _F$, which extends   
the automorphism $g_{\wp }$ of $\Gamma _{F,\bar\wp }$, and we obtain 
the following criterion: 
\medskip 

{\it $g_{\wp }\in\Aut\Gamma _{F,\bar\wp }$ can be extended 
to $g\in\Aut\Gamma _F$ if and only if $g_{\wp }$ is compatible with
the ramification filtration on $\Gamma _{F,\bar\wp }$.}
\medskip 

It would be interesting to understand how  ``global'' 
information about the 
embedding of $\Gamma _{F,\wp }$ into 
$\Gamma _F$ is reflected in ``local'' properties of the ramification 
filtration of $\Gamma_{F,\bar\wp }$. 
\medskip

Everywhere in the paper we use the following agreement about
compositions of morphisms: if $f:A\longrightarrow B$ and 
$g:B\longrightarrow C$ are morphisms then their composition will be
denoted by $fg$, in other words, if $a\in A$ then $(fg)(a)=g(f(a))$. 
One of the reasons is that when
operating with morphisms (rather that their values in $a\in A$) 
the notation $fg$ reflects much better the reality that 
$f$ is the first morphism and $g$ is the second. 
\medskip 

The author is very grateful to Ruth Jenni for very careful
checking 
of the final version of this paper and pointing out various
inexactitudes and misprints. 
\medskip 

\subhead 1. An analogue of the Magnus algebra for $\Gamma (p)$ 
\endsubhead 

In this section $K=k((t_K))$ is the local field of formal Laurent 
series 
with residue field $k=\Bbb F_{q_0}$, where 
$q_0=p^{N_0}$, $N_0\in\Bbb N$, and $t_K$ is a fixed uniformiser of $K$ 
(in most cases $t_K$ will be denoted just by $t$). We fix a choice of a
separable closure $K_{\sep }$ of $K$, denote by $K(p)$ the maximal 
$p$-extension of $K$ in $K_{\sep }$ and set 
$\Gamma =\Gal (K_{\sep }/K)$, $\Gamma (p)=\Gal (K(p)/K)$. 

\subsubhead{1.1 Liftings}
\endsubsubhead 
Notice first, that the uniformiser $t_K$ 
of $K$ can be taken as a $p$-basis for any 
finite extension $L$ of $K$ in $K_{\sep }$. For $M\in\Bbb N$, 
set 
$$O_M(L)=W_{M}(\sigma ^{M-1}L)[t_{K,M}]\subset W_{M}(L),$$ 
where $W_{M}$
is the functor of Witt vectors of length $M$, $\sigma $ 
is the $p$-th power map and $t_{K,M}=[t_K]=(t_K,0,\dots ,0)\in W_{M}(L)$ is 
the Teichm\" uller representative of $t_K$. Very often we shall 
use the simpler notation $t$ for $t_{K,M}$ (as well as for $t_K$). 
$O_M(L)$ is a lifting of
$L$ modulo $p^{M}$ or, in other words, it is a flat 
$W_{M}(\Bbb F_p)$-module such that $O_M(L)\operatorname{mod}p=L$. 
This is a special case of the construction of liftings in  
[B-M].  

Let $O_M(K_{\sep})$ be the inductive limit of all $O_M(L)$, where
$L\subset K_{\sep}$, $[L:K]<\infty $. Then we have a natural action of
$\Gamma $ on $O_M(K_{\sep})$ and $O_M(K_{\sep})^{\Gamma }=O_M(K)=
W_M(k)((t))$. We shall use again the notation 
$\sigma $ for the natural extension of $\sigma $ to 
$O_M(K_{\sep })$. Clearly,  $O_M(K_{\sep})|_{\sigma =\id }=W_{M}(\Bbb
F_p)$. 
Introduce the absolute liftings 
$O(K)=\mathbin{\underset{M}\to\varprojlim}O_M(K)$ and 
$O(K_{\sep })=\mathbin{\underset{M}\to\varprojlim}O_M(K_{\sep })$. 
Again we have $O(K_{\sep })^{\Gamma }=O(K)$ and 
$O(K_{\sep })|_{\sigma =\id }=W(\Bbb F_p)$. We can also consider 
the liftings $O_M(K(p))$ and $O(K(p))$ with the natural action of $\Gamma
(p)$ and similar properies.

Notice that   
for any $j\in O(K(p))$ there is an 
$i\in O(K(p))$ such that 
$\sigma (i)-i=j$.
\medskip 

\subsubhead{1.2. The algebra $\Cal A$}
\endsubsubhead 

Set $\Bbb Z(p)=\{a\in\Bbb N\ |\ (a,p)=1\}$ and $\Bbb Z^0(p)=\Bbb
Z(p)\cup\{0\}$. 
Let $\Cal A_k$ be the profinite associative $W(k)$-algebra 
with the set of pro-free generators 
$\{D_{an}\ |\ a\in\Bbb Z(p), n\in\Bbb
Z\operatorname{mod}N_0\}\cup\{D_0\}$. 

This means that $\Cal A_k=\mathbin{\underset{C,M}\to\varprojlim}\Cal
A_{CMk}$, where $C,M\in\Bbb N$, 
$$\Cal A_{CMk}=W_{M}(k)[[\{D_{an}\ |\ a\leqslant C,n\in\Bbb
Z\operatorname{mod}N_0\}]]$$ 
and the connecting morphisms 
$\Cal A_{C_1M_1k}\longrightarrow\Cal A_{C_2M_2k}$ are defined for
$C_1\geqslant C_2$, $M_1\geqslant M_2$ and  induced 
by the
correspondences $D_{an}\mapsto 0$ if $C_2<a\leqslant C_1$ and  
$D_{an}\mapsto D_{an}$ if $a\leqslant C_2$, and by the morphism 
$W_{M_1}(k)\longrightarrow W_{M_2}(k)$
of reduction modulo $p^{M_2}$. 
 
Denote again by $\sigma $ the extension
of the automorphism $\sigma $ of $W(k)$ to $\Cal A_k$ via 
the correspondences $\sigma :D_{an}\mapsto D_{a,n+1}$, where $a\in\Bbb
Z(p)$ , $n\in\Bbb Z\operatorname{mod}N_0$, and the correspondence 
$D_0\mapsto D_0$. 
Then $\Cal A:=\Cal A_k|_{\sigma =\id }$ is a pro-free $\Bbb
Z_p$-algebra: 
if $\beta _1,\dots ,\beta _{N_0}$ is a $\Bbb Z_p$-basis 
of $W(k)$ and, for $a\in\Bbb Z(p)$ and $1\leqslant r\leqslant N_0$,  
$$D_a^{(r)}:=\sum\Sb n\in\Bbb Z\operatorname{mod}N_0\endSb 
\sigma ^n(\beta _r)D_{an},$$
 then 
$\{D_a^{(r)}\ |\ a\in\Bbb Z(p), 1\leqslant r\leqslant
N_0\}\cup\{D_0\}$ 
is a set of pro-free generators of $\Cal A$. Notice also that if 
$\alpha _1,\dots ,\alpha _{N_0}\in W(k)$ is a dual basis for 
$\beta _1,\dots ,\beta _{N_0}$ (i.e. $\Tr (\alpha _i\beta _j)=\delta
_{ij}$, where 
$1\leqslant i,j\leqslant N_0$ and $\Tr $ is the trace of the field
extension 
$W(k)\otimes\Bbb Q_p$ over $\Bbb Q_p$) then 
for any $a\in\Bbb Z(p)$ and $n\in\Bbb Z\operatorname{mod}N_0$, it
holds 
$$D_{an}=\sum\Sb 1\leqslant r\leqslant N_0\endSb \sigma ^n(\alpha _r)
D_a^{(r)}.$$

Denote by $\Cal J$, resp. $\Cal J_{CM}$, the augmentation ideal in $\Cal
A$, 
resp. $\Cal A_{CM}$. 
Set $\Cal A_K:=\Cal A\hat\otimes O(K)$, $\Cal A_{CMK}=
\Cal A_{CM}\hat\otimes O(K)$, $\Cal A_{K(p)}=
\Cal A\hat\otimes O(K(p))$. We shall also use similar notation in other cases of 
extensions of scalars, e.g. 
$\Cal J_k=\Cal J\hat\otimes W(k)$, $\Cal J_K=\Cal J\hat\otimes O(K)$, 
$\Cal J_{K(p)}=\Cal J\hat\otimes O(K(p))$. 
\medskip

\subsubhead {\rm 1.3.} The embeddings $\psi _f$ 
\endsubsubhead 

Take $\alpha _0\in W(k)$ such that $\operatorname{Tr}(\alpha _0)=1$, 
where again $\operatorname{Tr}$ is the trace of the field extension 
$W(k)\otimes \Bbb Q_p\supset \Bbb Q_p$. 
For all $n\in\Bbb Z\operatorname{mod}N_0$, set 
$D_{0n}=\sigma^n(\alpha _0)D_0$ and 
introduce the element 
$$e=1+\sum\Sb a\in\Bbb Z^0(p)\endSb t^{-a}D_{a0}\in 1+\Cal J_K.$$
We shall use the same notation $e$ for the  
projections of $e$ to any of 
\linebreak 
$\Cal A_{CMK}\operatorname{mod}\Cal J_{CMK}^n$, 
where $C,M,n\in\Bbb N$.

\proclaim{Proposition 1.1} There is an $f\in 1+\Cal J_{K(p)}$ such that 
$\sigma (f)=fe$.
\endproclaim 

\demo{Proof} For $C,M,n\in\Bbb N$, set 
$$S_{CMn}=\left\{ f\in 1+\Cal J_{CMK(p)}\operatorname{mod}\Cal J_{CMK(p)}^n 
\ |\ \sigma f=fe\operatorname{mod}\Cal J_{CMK(p)}^n\right\}.$$

We use induction on $n\in\Bbb N$ to prove that for all $C,M,n\in\Bbb N$, 
$S_{CMn}\ne\emptyset $. 

Clearly, $S_{CM1}=\{1\}\ne\emptyset $. 

Suppose that $S_{CMn}\ne\emptyset $ and $f\in S_{CMn}$. Then 
$\sigma (f)=fe\operatorname{mod}\Cal J_{CMK(p)}^n$. Let 
$$\pi :1+\Cal J_{CMK(p)}\operatorname{mod}\Cal J_{CMK(p)}^{n+1}
\longrightarrow 1+\Cal J_{CMK(p)}\operatorname{mod}\Cal J_{CMK(p)}^n$$
be the natural projection. If $f_1\in 1+\Cal J_{CMK(p)}
\operatorname{mod}\Cal J_{CMK(p)}^{n+1}$ is such that 
$\pi (f_1)=f$ then $\sigma (f_1)=f_1e+j
\operatorname{mod}\Cal J_{CMK(p)}^{n+1}$, where 
$j\in \Cal J^n_{CMK(p)}$. 
There is an $i\in \Cal J^n_{CMK(p)}$ such that 
$\sigma (i)-i=j$, cf. n.1.1.  Therefore,  
$$\sigma (f_1-i)=f_1e+j-(i+j)=(f_1-i)e\operatorname{mod}\Cal J_{CMK(p)}^{n+1},$$
using that $ie=i\operatorname{mod}\Cal J_{CMK(p)}^{n+1}$,  and   
$S_{CM,n+1}\ne\emptyset $ because it contains $f_1-i$. 

Notice that each $S_{CMn}$ is a finite set and each $f\in S_{CMn}$ has
a finite field of definition. This follows from the fact that 
for any $C,M,n\in\Bbb N$, the $\Bbb Z_p$-module 
$\Cal A_{CM}\operatorname{mod}\Cal J_{CM}^n$ has finitely 
many free generators 
and, therefore, the equation $\sigma f=fe$ is equivalent to finitely 
many usual polynomial equations. 
Also notice that 
$\{S_{CMn}\ |\ C,M,n\in\Bbb N\}$ has a natural structure of projective system. 
Therefore, 
$\mathbin{\underset{C,M,n}\to\varprojlim}S_{CMn}\ne\emptyset $, and
any element  
$f$ of this projective limit 
satisfies $f\in 1+\Cal J_{K(p)}$ and $\sigma (f)=fe$. 

The proposition is proved.
\enddemo 

For any $f\in 1+\Cal J_{K(p)}$ such that $\sigma (f)=fe$ 
and $\tau\in\Gamma (p)$, set 
$\psi _f(\tau )=(\tau f)f^{-1}$. Clearly, 
$\sigma (\psi _f(\tau ))=\tau (\sigma f)(\sigma f)^{-1}=(\tau 
f)ee^{-1}f=\psi _f(\tau )$.  
Therefore, $\psi _f(\tau )\in (1+\Cal J_{K(p)})|_{\sigma =\id
}=1+\Cal J$. 

\proclaim{Proposition 1.2} 
{\rm a}) $\psi _f$ is a closed group embedding of 
$\Gamma (p)$ into $(1+\Cal J)^{\times }$. 
\newline 
{\rm b}) $\psi _f$ induces an isomorphism $\psi _f^{\ab }$ 
of the topological groups 
$\Gamma (p)^{\ab }$ and 
\linebreak 
$(1+\Cal J)^{\times }\operatorname{mod}\Cal J^2$. 
\newline 
{\rm c}) If $f_1\in 1+\Cal J_{K(p)}$ is such that $\sigma (f_1)=f_1e$ then 
there is an element $c\in 1+\Cal J$ such that for any $\tau\in\Gamma (p)$, 
$\psi _{f_1}(\tau )=c\psi _f(\tau )c^{-1}$.
\newline 
{\rm d}) $\psi _f$ induces an embedding of the group of all continuous automorphisms 
$\Aut \Gamma (p)$ into the group $\Aut\Cal A$ of 
continuous automorphisms of the $\Bbb Z_p$-algebra $\Cal A$. 
\endproclaim 

\demo{Proof} 
a) Clearly, $\psi _f$ 
can be treated as a pro-$p$-version of 
the embedding of the group $\Gamma (p)$ into its 
Magnus algebra. Therefore, by [Se, Ch.1, n.6], 
$\psi _f$ induces, for all $n\in\Bbb N$, the 
closed embeddings of the
quotients 
$C_n(\Gamma (p))/C_{n+1}(\Gamma (p))$ of 
commutator subgroups in
$\Gamma (p)$ into 
$1+\Cal J^n\operatorname{mod}\Cal J^{n+1}$. 
This implies that $\psi _f$ induces, for all $n\geqslant 1$, the 
closed 
group embeddings of $\Gamma (p)/C_n(\Gamma (p))$ into
$1+\Cal J\operatorname{mod}\Cal J^n$, 
and therefore, $\psi _f$ is a 
closed group monomorphism. 
\medskip 
 
 b) Consider the profinite $\Bbb Z_p$-basis 
$\{D_a^{(r)}\ |\ a\in\Bbb Z(p), 1\leqslant r\leqslant N_0\ \}\cup\{D_0\}$ for 
$\Cal J\operatorname{mod}\Cal J^2$ from n.1.2. 
For $1\leqslant r\leqslant N_0$, as earlier, consider $\alpha _r\in W(k)$, 
which form the dual basis  
of the basis $\{\beta _r\ |\ 1\leqslant r\leqslant N_0\}$  
chosen in n.1.2 to define the generators $D_a^{(r)}$. 
Then 
$$e=1+\sum\Sb 1\leqslant r\leqslant N_0 \\ a\in\Bbb Z(p) \endSb 
\alpha _rt^{-a}D_a^{(r)}+\alpha _0D_0$$ 
and 
$$f=1+\sum\Sb 1\leqslant r\leqslant N_0 \\ a\in\Bbb Z(p) \endSb 
f_a^{(r)}D_a^{(r)}+f_0D_0\operatorname{mod}\Cal J_{K(p)}^2,$$
where for $1\leqslant r\leqslant N_0$ and $a\in\Bbb Z(p)$, 
$f_a^{(r)}$ and $f_0$ belong to 
$O(K(p))\subset W(K(p))$ and satisfy the 
equations  
$\sigma f_a^{(r)}-f_a^{(r)}=\alpha _rt^{-a}$ and 
$\sigma f_0-f_0=\alpha _0$. 

Then for any $\tau\in\Gamma (p)$, 
$$\psi _f(\tau )=1+\sum\Sb a,r\endSb (\tau f_a^{(r)}-f_a^{(r)})D_a^{(r)}+
(\tau f_0-f_0)D_0\operatorname{mod}\Cal J_{K(p)}^2$$
and the identification 
$\psi _f:\Gamma (p)^{\ab }\simeq (1+\Cal J)^{\times
}\operatorname{mod}\Cal J^2$ 
is equivalent to the identifications of Witt-Artin-Schreier theory 
$$\oplus \Sb a\in\Bbb Z(p)\endSb W(k)t^{-a}\oplus W(\Bbb F_p)\alpha _0=
O(K)/(\sigma -\id )O(K)=\Hom _{\operatorname{cts}}(\Gamma (p), W(\Bbb F_p)).$$

c) Clearly, $\sigma (f_1f^{-1})=\sigma (f_1)\sigma (f)^{-1}=f_1ee^{-1}f^{-1}=f_1f^{-1}$. 
Therefore, 
$$f_1f^{-1}=c\in (1+\Cal J_{K(p)})\cap\Cal A=1+\Cal J$$ 
and for any $\tau\in\Gamma (p)$, 
$$\psi _{f_1}(\tau )=\tau (f_1)f_1^{-1}=\tau (cf)(cf)^{-1}=c(\tau
f)f^{-1}c^{-1}=c\psi _f(\tau )c^{-1}.$$

d) This also follows from the above mentioned 
interpretation of $\Cal A$ as a profinite analogue of the Magnus algebra for $\Gamma (p)$.
\enddemo 

\subsubhead {\rm 1.4.} The identification $\psi _f^{\ab }$ 
\endsubsubhead 

As it was already mentioned in the proof of 
proposition 1.2 the identification $\psi _f^{\ab }$ 
comes from the isomorphism of 
Witt-Artin-Schreier theory 
$$\Gamma (p)^{\ab }=\Hom (O(K)/(\sigma -\id )O(K), W(\Bbb F_p))$$
and does not depend on the choice of $t=t_K$ and $f\in 1+\Cal J_{K(p)}$. 
Suppose $\tau _0\in\Gamma (p)^{\ab }$ is such that 
$\psi _f^{\ab }(\tau _0)=1+D_0$ and 
for $a\in\Bbb Z(p)$ and $1\leqslant r\leqslant N_0$, the elements
$\tau _a^{(r)}\in\Gamma (p)^{\ab }$ are such that 
$\psi _f^{\ab }(\tau _a^{(r)})=1+D_a^{(r)}\operatorname{mod}\Cal J^2$.  
Then the element 
$$e=1+\alpha _0D_0+\sum\Sb a,r\endSb \alpha _rt^{-a}D_a^{(r)}$$ 
corresponds to the diagonal element 
$\alpha _0\otimes \tau _0+\sum\Sb a,r \endSb \alpha _r t^{-a}\otimes\tau _a^{(r)}$ from 
$O(K)\otimes\Gamma (p)^{\ab }=$
$$O(K)\otimes\Hom (O(K)/(\sigma -\id )O(K),\Bbb Z_p)
=\Hom (O(K)/(\sigma -\id )O(K), O(K)),$$ 
which comes from the following natural embedding  
$$O(K)/(\sigma -\id )O(K)=\oplus\Sb a\in\Bbb Z(p) \endSb
W(k)t^{-a}\oplus W(\Bbb F_p)\alpha _0\subset O(K).$$

 The above elements $\tau _0$, resp. $\tau
_a^{(r)}$, correspond 
to $t$, resp. 
$E(\beta _r, t^{a})^{1/a}$, by the reciprocity map of local class
field theory. (Here $\beta _1,\dots ,\beta _{N_0}\in W(k)$ were chosen
in n.1.2 and for $\beta\in W(k)$, 
$$E(\beta ,X)=\exp (\beta X+(\sigma \beta )X^p/p+\dots +(\sigma
^n\beta )X^{p^n}/p^n+\dots )\in W(k)[[X]]$$
is the generalisation of the Artin-Hasse exponential introduced by
Shafarevich [Sh].) This fact follows from the Witt explicit reciprocity
law, cf. [Fo]. Then the elements $D_{an}$, where 
$a\in\Bbb Z(p)$ and $n\in\Bbb Z\operatorname{mod}N_0$, 
correspond to 
$$\sum\Sb 1\leqslant r\leqslant N_0 \endSb \sigma ^n(\alpha _r)\otimes
E(\beta _r, t^a)^{1/a}\in W(k)\otimes _{\Bbb Z_p}\Cal G _a,$$
where the (multiplicative) group  
$\Cal G _a:=\{E(\gamma ,t^{a})\ |\ \gamma\in W(k)\}$ is  
identified with the $\Bbb Z_p$-module of Witt vectors  
$W(k)$ via the map $E(\gamma
,t^a)^{1/a}\mapsto\gamma $. Consider the identification 
$$W(k)\otimes _{\Bbb Z_p}W(k)=\oplus \Sb m\in\Bbb
Z\operatorname{mod}N_0\endSb W(k)_m$$
given by the correspondence $\alpha\otimes\beta\mapsto  \{\sigma
^{-m}(\alpha )\beta \}_{m\in\Bbb Z\operatorname{mod}N_0}$. Under this 
identification  the element $D_{an}$ corresponds to 
the vector $\delta _n\in \oplus _mW(k)_m$, which has $n$-th coordinate 
1 and 
all remaining coordinates 0. 
This interpretation of the generators $D_{an}$ 
will be applied below in the following
situation. Suppose $[k':k]<\infty $, $k'\simeq\Bbb F_{q'_0}$ with
$q'_0=p^{N'_0}$. Clearly, $N'_0\equiv 0\operatorname{mod}N_0$. 
For $a\in\Bbb Z(p)$ and $n\in\Bbb Z\operatorname{mod}N'_0$ denote by 
$D'_{an}$ an analogue of $D_{an}$ constructed for $K'=k'((t_{K'}))$ 
with $t_{K'}=t$. 
Let $\Gamma '=\Gal (K_{\sep }/K')$ and let 
$\Gamma '(p)$ be the Galois group of the 
maximal $p$-extension $K'(p)$ of $K'$ in $K_{\sep }$.  
With the above notation we have the following property:

\proclaim{Proposition 1.3} For any $a\in\Bbb Z(p)$ and 
$n\in\Bbb Z\operatorname{mod}N'_0$, 
$D'_{an}$ is mapped to $D_{a,n\operatorname{mod}N_0}$ 
under the map $\Gamma '(p)^{\ab }
\longrightarrow \Gamma (p)^{\ab }$, which 
is induced by the natural embedding 
$\Gamma '\subset\Gamma  $. 
\endproclaim 
\medskip

\subhead 2. Action of analytic automorphisms on $I^{\ab }(p)$ 
\endsubhead 
\medskip 

As earlier, $K=k((t))$, $k\simeq\Bbb F_{q_0}$ with 
$q_0=p^{N_0}$ and $\Gamma (p)=\Gal (K(p)/K)$. Let 
$I(p)$ be the inertia subgroup of $\Gamma (p)$ and let 
$I(p)^{\ab }$ be its image in the maximal abelian quotient $\Gamma
(p)^{\ab }$ of $\Gamma (p)$. 
\medskip 

2.1. Consider the group $\Aut K$ of continuous field automorphisms 
of $K$. Let $\Fr (t)\in\Aut K$ be such that 
$\Fr (t)|_k=\sigma $ and $\Fr
(t):t\mapsto t$. 
Then any  element of $\Aut K$ is the composition of a power 
$\Fr (t)^{n}$, where  $n\in\Bbb Z\operatorname{mod}N_0$,  
and  a field automorphism from $\Aut ^0(K):=
\{\eta\in\Aut K\ |\ \eta |_k=\id \}$. 
Notice that any $\eta\in\Aut ^0K$ is uniquely determined by 
the image 
$\eta (t)$ of $t$, 
which is again a uniformizer in $K$.

Let 
$\Aut _KK(p)$ be the group of continuous automorphisms $\bar\eta $
of $K(p)$ such that $\bar\eta |_K\in\Aut K$. 
Then $\Aut _KK(p)$ acts on $\Gamma (p)$: if 
$\bar\eta\in\Aut _KK(p)$ and $\tau\in\Gamma (p)$ then the action of
$\bar\eta $ is given by the correspondence   
$\tau\mapsto\bar\eta ^*(\tau )=\bar\eta ^{-1}\tau\bar\eta $, i.e. 
$\bar\eta ^*(\tau ): K(p)\mathbin{\overset{\bar\eta
^{-1}}\to\longrightarrow} 
K(p)\mathbin{\overset{\tau }\to\longrightarrow} K(p)
\mathbin{\overset{\bar\eta }\to\longrightarrow} K(p)$, cf. the 
introduction for the 
agreement about compositions of maps. 
The action induced by $\bar\eta ^*\in\Aut _KK(p)$ on $\Gamma (p)^{\ab }$ 
depends only on $\eta :=\bar\eta |_K$ and will be denoted simply by
$\eta ^*$. 

\medskip 

2.2.   
Let $\Cal M=I (p)^{\ab }\otimes\Bbb F_p$. 
If $U_K$ is the group of principal units 
in $K$ then we shall use the identification 
$\Cal M=U_K/U_K^{p}$, which is given by the reciprocity 
map of local class field theory. 
Notice that, with respect of this identification,  
for any 
$\eta\in\Aut K$, the action $\eta ^*$ comes 
from the natural
action of $\eta $ on $K$. 
We shall  denote the $k$-linear extension of the action 
of $\eta $ to $\Cal M_k:=\Cal M\otimes _{\Bbb F_p}k$ 
 by the same symbol 
$\eta ^*$. 

Use the map $m\mapsto (\psi _f^{\ab }(m)-1)\operatorname{mod}p$ 
to identify $\Cal M_k$ with a submodule of 
$\Cal J_k\operatorname{mod}(p,\Cal J_k^2)$. 
For $a\in\Bbb Z(p)$ and $n\in\Bbb Z\operatorname{mod}N_0$, 
consider the images of the elements  
$D_{an}$, where 
$a\in\Bbb Z(p)$ and 
$n\in\Bbb Z\operatorname{mod}N_0$ (cf. n.1),  
in $\Cal J_k\operatorname{mod}(p,\Cal J_k^2)$.   
Denote these images by same symbols. Then 
they give a set of free topological generators of 
the $k$-module $\Cal M_k$.  The action of 
$\eta\in\Aut K$ on $\Cal M_k$ in terms of these generators is as follows.

\proclaim{Proposition 2.1} {\rm 1)} $\Fr (t)^*(D_{an})=D_{a,n-1}$;
\newline 
{\rm 2)} if $\eta\in\Aut ^0K$, then 
$$\sum\Sb a\in\Bbb Z(p)  \endSb t^{-a}\eta ^*(D_{a0})\equiv 
\sum\Sb a\in\Bbb Z(p) \endSb \eta ^{-1}(t)^{-a}D_{a0}
\operatorname{mod}(k+(\sigma -\id )K)\otimes\Cal M.$$
\endproclaim 

\demo{Proof} 1) Consider the generators 
$\alpha _rD_a^{(r)}$ of $\Cal A$ from n.1.2,  
 where $a\in\Bbb Z(p), 1\leqslant r\leqslant
N_0$. 
Note that the residue of the corresponding 
element $e-1$ modulo 
\linebreak 
$(\sigma -\id )K\otimes (\Cal J\operatorname{mod}\Cal J^2)$ 
does not depend on the choice of $t$ or of the elements 
$\alpha _1,\alpha _2,\dots ,\alpha _{N_0}$, because this is the diagonal
element of Artin-Schreier duality. Therefore, if 
$\Fr (t)^*(D_a^{(r)})=D^{\prime (r)}_a$ and $\Fr (t)^*(D_0)=D_0'$ 
then 
$$e-1\ \ \equiv \ \ \sigma (\alpha _0)\otimes D'_0+
\sum\Sb a,r \endSb \sigma (\alpha _r)t^{-a}
\otimes D_a^{\prime (r)}$$
$$\equiv \alpha _0\otimes D_0+
\sum\Sb a,r \endSb \alpha _rt^{-a}\otimes D_a^{(r)}
\operatorname{mod}(\sigma -\id )K\otimes 
(\Cal J\operatorname{mod}\Cal J^2).\tag {2.1}$$
So, for any $a\in\Bbb Z(p)$, we see that  in 
$k\otimes _{\Bbb F_p}\Cal M=\Cal M_k$  
$$D_{a0}=\sum\Sb r \endSb \alpha _r\otimes D_a^{(r)}=\sum\Sb 
r \endSb \sigma (\alpha _r)\otimes D^{\prime (r)}_a.$$
Denoting the $k$-linear extension of 
$\Fr (t)^*$ by the same symbol, as usual, we have 
$$\Fr (t)^*(D_{a0})=\sum\Sb r\endSb 
\alpha _r\otimes\Fr (t)^*(D_a^{(r)})
=\sum\Sb r \endSb \alpha _r\otimes D_a^{\prime
(r)}=\sigma ^{-1}D_{a0}=D_{a,-1}.$$
Therefore, for any $a\in\Bbb Z(p)$ and 
$n\in\Bbb Z\operatorname{mod}N_0$, 
$\Fr (t)^*(D_{an})=D_{a,n-1}$. Notice also that 
congruence (2.1) implies that 
$\Fr (t)^*D_0=D_0$. 
\medskip 

2) Using that $\eta $ is a $k$-linear automorphism of $K$ and 
proceeding similarly to the above part 1) we obtain that 
$$\sum\Sb a\in\Bbb Z(p)^0 \endSb \eta (t)^{-a}\eta ^*(D_{a0})\equiv
\sum \Sb a\in\Bbb Z(p)^0 \endSb t^{-a}D_{a0}
\operatorname{mod}(\sigma -\id )K\otimes \Cal M.$$
Now apply $(\eta ^{-1}\otimes\id )$ to both sides of this
congruence 
and notice that we can omit the terms with index $a=0$ 
when working 
modulo $(k+(\sigma -\id )K)\otimes\Cal M$, because they 
belong to $\Cal M_k$. 
The lemma is proved. 
\enddemo 
\medskip

2.3. If $f$ is a continuous automorphism of the $\Bbb F_p$-module 
$\Cal M$, we agree to use the same notation $f$ for its 
$k$-linear extension to an automorphism of $\Cal M_k$. 
For any $a\in\Bbb Z(p)$, set

$$f(D_{a0})=\sum\Sb b\in\Bbb Z(p), 
n\in\Bbb Z\operatorname{mod}N_0\endSb \alpha _{abn}(f)D_{bn}.$$
Then all coefficients $\alpha _{abn}(f)$ are in $k$. 
 Sometimes we shall 
use the notation $\alpha _{abn}(f)$ if $a$ or $b$ are divisible by
$p$, then it is assumed that $\alpha _{abn}(f)=0$. Notice that 
for any $m\in\Bbb Z\operatorname{mod}N_0$, 

$$f(D_{am})=\sum\Sb b\in\Bbb Z(p), n\in\Bbb Z\operatorname{mod}N_0
\endSb \sigma ^m(\alpha _{abn}(f))D_{b,n+m}.$$

\definition{Definition} For any $v\in\Bbb N$, let 
$\Cal M^{(v)}$ be the minimal closed 
$\Bbb F_p$-submodule in $\Cal M$ such that 
$\Cal M_k^{(v)}:=\Cal M^{(v)}\otimes k$  
is topologically generated over $k$ by all 
$D_{an}$, where $a\in\Bbb Z(p)$, $a\geqslant v$ and 
$n\in\Bbb Z\operatorname{mod}N_0$. (Notice that 
$\Cal M=\Cal M^{(1)}$.) 
\enddefinition

\definition{Definition}  
$\Aut _{\adm }\Cal M$ is the subset in the group $\Aut \Cal M$, 
consisting of all 
continuous 
$\Bbb F_p$-linear automorphisms $f$ satifying 
$\alpha _{a,b,m\operatorname{mod}N_0}(f)=0$ if 
$bp^m<a$, 
for any $a,b\in\Bbb Z(p)$ and $-N_0<m\leqslant 0$.  
\enddefinition 

It is easy to see that:  

\roster  
\item $\Aut _{\adm}\Cal M$ is a subgroup 
of $\Aut\Cal M$;  
\medskip 

\item if 
$f\in\Aut _{\adm }\Cal M$ then for any $a\in\Bbb N$, 
$f(\Cal M^{(a)})\subset \Cal M^{(a)}$, 
i.e. $f$ is compatible with the image of the ramification 
filtration in $\Cal M$;
\medskip 

\item  if $f\in\Aut _{\adm }\Cal M$ then 
for any $a\in\Bbb Z(p)$, $\alpha _{aa0}\in
k^*$ and 
$\alpha _{aan}(f)=0$ if $n\ne 0$. 
\endroster 
\medskip

\definition{Definition} For   
$f\in\Aut \Cal M$,  
let $f_{\an }\in \operatorname{End}\Cal M$ be
such that 
for all $a\in\Bbb Z(p)$, 
$$f_{\an }(D_{a0})=\sum\Sb b\in\Bbb Z(p)\endSb 
\alpha _{ab0}(f)D_{b0}.$$
\enddefinition 

\proclaim{Proposition 2.2} If $f,g\in\Aut_{\adm}\Cal M$ then 
for any $a,b\in\Bbb Z(p)$ such that 
$a\leqslant b<ap^{N_0}$, 
$$\alpha _{ab0}(fg)=\sum\Sb c \endSb \alpha _{ac0}(f)\alpha
_{cb0}(g).$$
\endproclaim

\proclaim{Corollary 2.3} If $v< p^{N_0}$ then 
the correspondence $f\mapsto f_{\an }$ 
is a group homomorphism from $\Aut _{\adm }\Cal M$ to 
$\Aut _{\adm }\Cal M\operatorname{mod}\Cal M^{(v)}$. 
\endproclaim 

\demo{Proof of proposition} We have 
$$\alpha _{ab0}(fg)=\sum\Sb m+n\equiv 0\operatorname{mod}N_0 \\
0\geqslant n,m>-N_0\endSb 
\alpha _{a,c,n\operatorname{mod}N_0}(f)\sigma ^n
(\alpha _{c,b,m\operatorname{mod}N_0}(g))D_{b,(m+n)\operatorname{mod}N_0}.$$
Then $\alpha _{a,c,n\operatorname{mod}N_0}(f)\ne 0$ implies that 
$cp^n\geqslant a$ and 
$\alpha _{c,b,m\operatorname{mod}N_0}(g)\ne 0$ implies 
that $bp^m\geqslant c$. So, if the
corresponding coefficient for $D_{b,(m+n)\operatorname{mod}N_0}$ is not zero then 
$bp^{m+n}\geqslant a$, i.e. $m+n>-N_0$ and, therefore, $m=n=0$. 
\enddemo 

The following proves that $\Aut ^0K\subset\Aut _{\adm}\Cal M$. 

\proclaim{Proposition 2.4} If $\eta\in\Aut ^0K$ then 
$\eta ^*\in\Aut_{\adm}\Cal M$. 
\endproclaim 

\demo{Proof} 
For $a\in\Bbb Z(p)$, set 
$$\eta ^{-1}(t)^{-a}\equiv \sum 
\Sb b\in\Bbb Z(p), s\geqslant 0\endSb \gamma _{abs}t^{-bp^s}
\operatorname{mod}k[[t]].$$
Clearly,  $\gamma _{abs}=0$ if $bp^s> a$. It follows from part 2) of 
proposition 2.1
that 
$$\eta ^*(D_{b0})=\sum\Sb a\in\Bbb Z(p), s\geqslant 0\endSb 
\sigma ^{-s}(\gamma _{abs})D_{a,-s\operatorname{mod}N_0}.$$
Therefore, for $0\leqslant m<N_0$, 
$$\alpha _{b,a,-m\operatorname{mod}N_0}(\eta ^*)=
\sum\Sb s\equiv m\operatorname{mod}N_0 \\
s\geqslant 0 \endSb 
\sigma ^{-s}(\gamma _{abs})$$
and $a/p^m<b$ implies for $s\equiv m\operatorname{mod}N_0$, 
$s\geqslant 0$, that 
$a/p^s<b$. So, $bp^s>a$, $\gamma _{abs}=0$ and 
$\alpha _{b,a,-m\operatorname{mod}N_0}(\eta ^*)=0$. 

The proposition is proved.
\enddemo  
 
2.4. In this subsection we prove three technical propositions. 
Notice that in proposition 2.5 we treat the case of fields of 
characteristic $p\ne 2$ and  in proposition 2.6 
the characteristic of $K$ is 2. 
Propositions 2.5-2.7 will be used   
later in section 5. If $a,b\in\Bbb N$ then $\delta _{ab}$ is the
Kronecker symbol. 

\proclaim{Proposition 2.5} Suppose $p\ne 2$, 
$w_0\in\Bbb N$, $w_0+1\leqslant
p^{N_0}$ and 
$f\in\Aut _{\adm }\Cal M$ is such that 
$\alpha _{1a0}(f)=\delta _{1a}$ if $1\leqslant a<w_0$ and 
$\alpha _{2a0}(f)=0$ if $a\equiv 1\operatorname{mod}p$ and 
$a\leqslant
w_0$. Then there is an $\eta\in\Aut ^0K$ such that $\eta (t)\equiv
t\operatorname{mod}t^{w_0}$, 
$\alpha _{1a0}(f\eta ^*)=\delta _{1a}$ if $1\leqslant a<w_0+1$, and 
$\alpha _{2a0}(f\eta ^*)=0$ if $a\equiv 1\operatorname{mod}p$  
and $a\leqslant w_0+1$.
\endproclaim 

\demo{Proof} Take $\eta\in\Aut ^0K$ such that 
$\eta ^{-1}(t)=t(1+\gamma t^{w_0-1})$ with $\gamma\in k$.  
Then for any $a\in\Bbb Z(p)$, 
$\eta ^{-1}(t^{-a})=t^{-a}(1-a\gamma t^{w_0-1})
\operatorname{mod}t^{-a+w_0}$, and part 2) of  
proposition 2.1  implies that 
$\alpha _{aa0}(\eta ^*)=1$, $\alpha _{ab0}(\eta ^*)=0$ 
if $a<b<a+w_0-1$, $\alpha _{a,a+w_0-1,0}(\eta ^*)=
-(a+w_0-1)\gamma $. 

Therefore, by proposition 2.2 $\alpha _{1a0}(f\eta ^*)=\delta _{1a}$ if 
$1\leqslant a<w_0$ and 
$\alpha _{2a0}(f\eta ^*)=0$ if $a\equiv 1\operatorname{mod}p$, 
$a\leqslant w_0$. 

Suppose $w_0\not\equiv 0\operatorname{mod}p$. Then 
by proposition 2.2 
$$\alpha _{1w_00}(f\eta ^*)=-w_0\gamma +
\alpha _{1w_00}(f)=0$$ 
if $\gamma =w_0^{-1}\alpha _{1w_00}(f)$. 
This proves the proposition in the case 
$w_0\not\equiv 0\operatorname{mod}p$,  because 
$w_0+1\not\equiv 1\operatorname{mod}p$ and 
no conditions are required for $\alpha _{2,w_0+1,0}(f\eta ^*)$. 

Suppose $w_0\equiv 0\operatorname{mod}p$. Then there are no 
conditions for 
$\alpha _{1w_00}(f\eta ^*)$ and by proposition 2.2 
$$\alpha _{2,w_0+1,0}(f\eta ^*)=
\alpha _{220}(f)\alpha _{2,w_0+1,0}(\eta ^*)+
\alpha _{2,w_0+1,0}(f)\alpha _{w_0+1,w_0+1,0}(\eta ^*)$$
$$\qquad\qquad\qquad\qquad =
-\alpha _{220}(f)\gamma +\alpha _{2,w_0+1,0}(f)=0$$ 
if $\gamma =\alpha _{2,w_0+1,0}(f)\alpha _{220}(f)^{-1}$. 
(Using that 
$f\in\Aut _{\adm }\Cal M$ hence 
 $\alpha _{220}(f)\in k^*$.)

The proposition is proved. 
\enddemo 

\proclaim{Proposition 2.6} Let $M\in\Bbb N$, $p=2$, $w_0=4M$ and
$w_0+1<2^{N_0}$. Suppose $f\in\Aut _{\adm }\Cal M$ is such that 
$\alpha _{1a0}(f)=\delta _{1a}$ if $1\leqslant a\leqslant w_0-3$ and 
$\alpha _{3a0}(f)=\delta _{3a}$ if $3\leqslant a\leqslant w_0-1$. 
Then there is an $\eta\in\Aut ^0K$ such that 
$\alpha _{1a0}(f\eta ^*)=\delta _{1a}$ and 
$\alpha _{3a0}(f\eta ^*)=\delta _{3a}$ if $a\leqslant w_0+1$. 
\endproclaim 

\demo{Proof} {\it 1st step.} 

Take $\eta _1\in\Aut ^0K$ such that $\eta _1^{-1}(t)=t(1+\gamma
_1t^{4M-2})$ with $\gamma _1\in k$. Then for $a\in\Bbb Z(2)$, 
$\eta _1^{-1}(t^{-a})\equiv t^{-a}(1+\gamma
_1t^{4M-2})\operatorname{mod}t^{-a+4M}$ and 
by part 2) of proposition 2.1,   
$\alpha _{aa0}(\eta _1^*)=1$, 
$\alpha _{ab0}(\eta _1^*)=0$ if $a<b<a+4M-2$, and 
$\alpha _{a,a+4M-2,0}(\eta _1^*)=\gamma _1$. 

So by proposition 2.2,  
$\alpha _{1a0}(f\eta _1^*)=\alpha _{1a0}(f)$ if 
$a\leqslant 4M-3=w_0-3$, 
$\alpha _{3a0}(f\eta _1^*)=\alpha _{3a0}(f)$ if 
$a\leqslant 4M-1=w_0-1$, 
$\alpha _{1,w_0-1,0}(f\eta _1^*)=\alpha _{1,w_0-1,0}(f)+
\alpha _{1,w_0-1,0}(\eta _1^*)=0$ if 
$\gamma _1=\alpha _{1,w_0-1,0}(f)$. 
\medskip 

{\it 2nd step.}

By the above first step we can now assume that 
$\alpha _{1,w_0-1,0}(f)=0$. 

Take $\eta _2\in\Aut ^0K$ such that $\eta _2^{-1}(t)=t(1+\gamma
_2t^{2M-1})$. Then for $a\in\Bbb Z(2)$, 
$\eta _2^{-1}(t^{-a})\equiv t^{-a}(1+\gamma
_2t^{2M-1}+\delta (a)\gamma
_2^2t^{4M-2})\operatorname{mod}t^{-a+4M}$, where 
$\delta (a)=a(a+1)/2$. 

So by part 2) of proposition 2.1,   
$\alpha _{aa0}(\eta _2^*)=1$, 
$\alpha _{ab0}(\eta _2^*)=0$ if $a<b<a+4M-2$ 
(notice that $-a+2M-1\equiv 0\operatorname{mod}2$), and 
$\alpha _{a,a+4M-2,0}(\eta _2^*)=\delta (a+4M-2)\gamma _2^2$ 
(notice that $\delta (a+4M-2)=0$ if $a\equiv 1\operatorname{mod}4$ and
$\delta (a+4M-2)=1$ if $a\equiv 3\operatorname{mod}4$). 

Again by proposition 2.2, $\alpha _{1a0}(f\eta _2^*)=\alpha _{1a0}(f)$ if 
$a\leqslant 4M-1=w_0-1$ (use that 
$\alpha _{1,w_0-1,0}(f)=\alpha _{1,w_0-1,0}(\eta _2^*)=0$), 
$\alpha _{3a0}(f\eta _2^*)=\alpha _{3a0}(f)$ if 
$a\leqslant 4M-1=w_0-1$, 
$\alpha _{3,w_0+1,0}(f\eta _2^*)=\alpha _{3,w_0+1,0}(f)+
\alpha _{3,w_0+1,0}(\eta _2^*)=0$ if 
$\gamma _2\in k$ is such that $\gamma _2^2=\alpha _{3,w_0+1,0}(f)$. 
\medskip 

{\it 3rd step.}

Now we can assume that $\alpha _{1,w_0-1,0}(f)=\alpha _{3,w_0+1,0}(f)=0$.

Take $\eta _3\in\Aut ^0K$ such that $\eta _3^{-1}(t)=t(1+\gamma
_3t^{4M})$. Then for $a\in\Bbb Z(2)$, 
$\eta _3^{-1}(t^{-a})\equiv t^{-a}(1+\gamma
_3t^{4M})\operatorname{mod}t^{-a+4M+2}$,  
$\alpha _{aa0}(\eta _3^*)=1$, 
$\alpha _{ab0}(\eta _3^*)=0$ if $a<b<a+4M$, and 
$\alpha _{a,a+4M,0}(\eta _3^*)=\gamma _3$. 

This implies that $\alpha _{1a0}(f\eta _3^*)=\alpha _{1a0}(f)$ if 
$a\leqslant 4M-1=w_0-1$, 
$\alpha _{1,w_0+1,0}(f\eta _3^*)=\alpha _{1, w_0+1,0}(f)+
\alpha _{1,w_0+1,0}(\eta _3^*)=0$ if $\gamma _3=\alpha
_{1,w_0+1,0}(f)$ and   
$\alpha _{3a0}(f\eta _3^*)=\alpha _{3a0}(f)$ if 
$a\leqslant w_0+1$. 

The proposition is proved. 
\enddemo 
\medskip

\proclaim{Proposition 2.7} Suppose $a\in\Bbb Z(p)$,  
$w_0\leqslant ap^{N_0}$, where $w_0\in p\Bbb N$, 
$w_0>a+1$  if $p\ne 2$ and 
$w_0\in 4\Bbb N$, $w_0>a+2$  if $p=2$. 
Suppose $\eta ,\eta _1\in\Aut ^0K$ are such that 
for any $b,c\in\Bbb Z(p)$ satisfying the restrictions  
$a\leqslant c\leqslant b<w_0\leqslant ap^{N_0}$, we have the equality 
$$\alpha _{cb0}(\eta ^*)=\alpha _{cb0}(\eta _1^*).$$
Then $\eta (t)\equiv \eta _1(t)\operatorname{mod}t^{v_0}$, where 
$v_0=w_0-a+1$ if $p\ne 2$ and $v_0=(w_0-a+1)/2$ 
if $p=2$.     
\endproclaim 

\remark{Remark} With notation from n.2.3 this proposition 
implies  that if 
\linebreak 
${\eta _1}^*_{\an }\equiv \eta _{\an }^*\operatorname{mod}\Cal
M^{(w_0)}$ then $\eta (t)\equiv \eta _1(t)\operatorname{mod}t^{v_0}$. 
\endremark 

\demo{Proof} Use proposition 2.2 to reduce the proof to the case $\eta
_1(t)=t$. 

Suppose, first, that  $\eta ^{-1}(t)=\alpha
t\operatorname{mod}t^2$. Then 
$$\alpha _{cc0}(\eta ^*)=\alpha ^{-c}=1. \tag {2.2} $$. 
If $a+1\in\Bbb Z(p)$ then $p\ne 2$ and we can 
use formula (2.2) for $c=a,a+1$ to prove that 
$\alpha =1$. 
Suppose $a+1\notin\Bbb Z(p)$. If $p=2$ use (2.2) 
for $c=a,a+2<w_0$, and if 
$p\ne 2$ use (2.2) for $c=a+2,a+3<w_0$ to prove again that $\alpha =1$. 

Assume now that $p\ne 2$. 

Suppose $\eta ^{-1}(t)\equiv t+\alpha t^{v-1}\operatorname{mod}t^v$
with $v\geqslant 3$ and $\alpha\in k^*$. If $a+v-2\in\Bbb Z(p)$ then 
by part 2) of proposition 2.1  
$\alpha _{a,a+v-2,0}(\eta ^*)\ne 0$. This implies that  
$a+v-2\geqslant w_0+1$, i.e. $v\geqslant w_0-a+1$,  
as required. If 
$a+v-2\equiv 0\operatorname{mod}p$ then 
by part 2) of proposition 2.1  
$\alpha _{a+1,a+v-1,0}(\eta ^*)\ne 0$.This 
implies that  
$a+v-1\geqslant w_0+1$ and 
$v\geqslant w_0-a+2>w_0-a+1$. The case $p\ne 2$ is considered. 

Assume now that $p=2$. 

Suppose that $M\in\Bbb N$ is such that 
$$\eta ^{-1}(t)=t\left (1+\sum\Sb r\geqslant 2M-1 \endSb \gamma
_rt^r\right )\equiv t\operatorname{mod}t^{2M}$$
with either $\gamma _{2M-1}\ne 0$ or $\gamma _{2M}\ne 0$.

Therefore,  if $r\equiv 0\operatorname{mod}2$, 
$r\geqslant 2M-1$ and 
$a+r<ap^{N_0}$ then by part 2) of proposition 2.1  
$\alpha _{a,a+r,0}(\eta ^*)=\gamma _r$. This implies that 
either $2M\geqslant w_0$ (and the proposition 
is proved) or $2M\leqslant w_0-2$, 
$\gamma _{2M}=0$ and $\gamma _{2M-1}\ne 0$. 

Suppose $a+4M<w_0$. Then with the notation from the second step 
in the proof of proposition 2.6, we have 
 
$$\alpha _{a,a+4M-2,0}(\eta ^*)=\gamma _{4M-2}
+\gamma _{2M-1}^2\delta (a+4M-2)$$
 
$$\alpha _{a+2,a+4M,0}(\eta ^*)=\gamma _{4M-2}
+\gamma _{2M-1}^2\delta (a+4M).$$

The sum of the right hand sides of the above two 
equalities is $\gamma
_{2M-1}^2\ne 0$,  
because $\delta (a+4M-2)+\delta (a+4M)=1$. 
 Therefore, at least one of their left hand sides is
not zero. This means that the assumption about 
$a+4M<w_0$ was wrong. Therefore, $4M>w_0-a$ and 
$2M\geqslant (w_0-a+1)/2$. 

The proposition is proved. 
\enddemo 
\medskip

\subhead 3. Compatible systems of group morphisms 
\endsubhead

For any $s\in\Bbb Z_{\geqslant 0}$, let $K_s$ be the unramified extension 
of $K$ in $K(p)$ of degree $p^s$. Then 
$K_s=k_s((t))$, where $t=t_K$ is a fixed uniformiser, 
$k\subset k_s$, $[k_s:k]=p^s$, 
$k_s\simeq\Bbb F_{q_s}$, $q_s=p^{N_s}$ with $N_s=N_0p^s$. 

Let $K_{\ur }$ be the  
union of all $K_s$, $s\geqslant 0$. 
This is the maximal unramified extension of $K$ in $K(p)$ and 
its residue field coincides with the residue field $k(p)$ of $K(p)$. 
Let $I_{K_{\ur }}(p)^{\ab }$, resp. $I_{K_s}(p)^{\ab }$,  for $s\in\Bbb
Z_{\geqslant 0}$, be the images of the inertia
subgroups of $\Gal (K(p)/K_{\ur })$, resp. $\Gal (K(p)/K_s)$, in the
corresponding maximal abelian quotients. 
Then 
$I_{K_{\ur }}(p)^{\ab }=
\mathbin{\underset{s}\to\varprojlim}I_{K_s}(p)^{\ab }$. 

3.1. For $s\geqslant 0$, introduce the $\Bbb F_p$-modules 
$\Cal M_{Ks}=I_{K_s}(p)^{\ab }\otimes \Bbb
F_p$   
and $\Cal M_{K{\ur }}=I_{K_{\ur }}(p)^{\ab }\otimes \Bbb F_p$ with the
corresponding 
$k(p)$-modules $\bar {\Cal M}_{Ks}=\Cal M_{Ks}\hat\otimes _{\Bbb F_p} k(p)$ and 
$\bar\Cal M_{K{\ur }}=\Cal M_{K{\ur }}\hat\otimes _{\Bbb F_p} k(p)$. Then 
for all $s\geqslant 0$, we have  
natural connecting morphisms 
$j_s:\Cal M_{K,{s+1}}\longrightarrow\Cal M_{Ks}$ and 
$\bar\jmath _s:\bar \Cal M_{K,{s+1}}\longrightarrow\bar\Cal M_{Ks}$ 
(both 
are induced by the natural group embeddings 
$\Gamma_{K_{s+1}}\longrightarrow\Gamma _{K_s}$). Therefore, 
we have  projective systems 
$\{\Cal M_{Ks},j_s\}$ and $\{\bar\Cal M_{Ks}, \bar\jmath _s\}$ and 
 natural identifications 
$\Cal M_{K{\ur }}=\mathbin{\underset{s}\to\varprojlim}\Cal M_{Ks}$ and 
$\bar\Cal M_{K_{\ur }}=\mathbin{\underset{s}\to\varprojlim}\bar\Cal M_{K_s}$. 

Let $\Cal M_{K\infty }$ be the $k(p)$-submodule in $\bar\Cal M_{K{\ur
}}$ which is topologically generated by all 
$D^{\infty }_{an}:=\mathbin{\underset{s}\to\varprojlim} 
D^{(s)}_{a,n\operatorname{mod}N_s}$, where 
$a\in\Bbb Z(p)$ and $n\in\Bbb Z$. Here for $a\in\Bbb Z(p)$ and 
$n\in\Bbb Z\operatorname{mod}N_s$, 
$D^{(s)}_{an}$ are generators for $\bar\Cal M_{Ks}$, which
are 
analogues of the generators $D_{an}$ 
introduced in section 2 for the $k$-module $\Cal M_k$. 
Notice that the generators $D^{(s)}_{an}$  depend on the choice of the 
uniformising element $t$ in $K$. 

\proclaim{Proposition 3.1} The $k(p)$-submodule $\Cal M_{K\infty }$ 
of $\bar\Cal M_{K{\ur }}$ does not depend on the choice of $t$. 
\endproclaim 

\demo{Proof} Let $t_1$ be another uniformiser in $K$. Introduce 
$\eta\in\Aut ^0(K_{\ur })$ such that 
$\eta (t)=t_1$. The proposition will be proved if we show that 
$\eta ^*(\Cal M_{K\infty })=\Cal M_{K\infty }$. 

For $s\geqslant 0$, let $\eta _s=\eta |_{K_s}\in\Aut ^0K_s$. 
Then for $a\in\Bbb Z(p)$ and $n\in\Bbb Z\operatorname{mod}N_s$, 
$$\eta _s^*(D^{(s)}_{an})=\sum\Sb b\in\Bbb Z(p), 
m\in\Bbb Z\operatorname{mod}N_s\endSb \sigma ^n\alpha
_{abm}(\eta ^*_s)D_{b,m+n}^{(s)},$$
where the coefficients  $\alpha _{abm}(\eta ^*_s)\in k_s$ satisfy the 
following compatibility
conditions (using that $j_s(D^{(s)}_{an})=
D^{(s-1)}_{a,n\operatorname{mod}N_{s-1}}$): 
\medskip 

{\it if $a,b\in\Bbb Z(p)$ and 
$m\in\Bbb Z\operatorname{mod}N_{s-1}$ then  }
$$\sum\Sb n\operatorname{mod}N_{s-1}=m\endSb\alpha _{abn}(\eta ^*_s) 
=\alpha _{abm}(\eta ^*_{s-1}).$$ 
\medskip 

By proposition 2.4, if 
$0\leqslant m<N_s$ and $b/p^m<a$ then 
$\alpha _{a,b,-m\operatorname{mod}N_s}(\eta ^*_s)=0$. Therefore, 
if $s$ is such that $b/p^{N_s}<a$ then  
$\alpha ^{\infty }_{a,b,-m}(\eta ^*):=
\alpha _{a,b,-m\operatorname{mod}N_s}(\eta ^*_s)$ 
does not depend on $s$ 
and 
for any $a\in\Bbb Z(p)$ and $n\in\Bbb Z_{\geqslant 0}$, 
$$\eta ^*(D^{\infty }_{an})=\sum\Sb b\in\Bbb Z(p), m\geqslant 0 \endSb 
\sigma ^n\alpha ^{\infty }_{a,b,-m}(\eta ^*)D^{\infty }_{b,n-m}\in\Cal
M_{K\infty }.$$
 The proposition is proved. 
\enddemo 
\medskip 

3.2. Consider the identification of class field theory 
$I_{K_s}(p)^{\ab }=U_{K_s}$, where 
$U_{K_s}$ is the group of principal units of $K_s$. Define the 
continuous morphism of topological $k(p)$-modules 
$$\pi _{Ks}:\bar\Cal M_{Ks}=I_{K_s}(p)^{\ab }\hat\otimes k(p)
\longrightarrow\hat\Omega ^1_{O_{K_{\ur }}},$$
by 
$\pi _{Ks}(u\otimes\alpha )=\alpha \operatorname{d}(u)/u$ 
for $u\in U_{K_s}$ and $\alpha\in k(p)$.
Here $\hat\Omega ^1_{O_{K_{\ur }}}$ is the completion of 
the module of differentials of the valuation ring $O_{K_{\ur }}$ 
with respect to the $t$-adic topology.  
 Notice that for any 
$a\in\Bbb Z(p)$ and $0\leqslant n<N_s$, 
$$D^{(s)}_{a,n\operatorname{mod}N_s}=\sum\Sb 0\leqslant i<N_s\endSb 
u_i\otimes (\sigma ^n\alpha _i\operatorname{mod}p).$$
Here $\{\alpha _i\ |\ 1\leqslant i\leqslant N_s\}$ is a 
$\Bbb Z_p$-basis of $W(k_s)$. If 
$\{\beta _i\ |\ 1\leqslant i\leqslant N_s\}$ is its 
dual basis then for $1\leqslant i\leqslant N_s$,  
$u_i=E(\beta _i,t^a)^{1/a}$, cf. n.1.4.  Therefore, 

$$\pi _{Ks}(D^{(s)}_{a,n\operatorname{mod}N_s})=\left (\sum\Sb
i\geqslant 0\endSb t^{ap^{n+iN_s}}\right )\frac{\operatorname{d}(t)}{t}.$$
It is easy to see that $\pi _{K{\ur }}:=\varprojlim\pi _{Ks}$ is 
a continuous map from $\bar\Cal M_{K{\ur }}$ to $\hat\Omega ^1_{O_{K_{\ur
}}}$. 

Notice that if $\bar n=\varprojlim
(n_s\operatorname{mod}N_s)\in
\mathbin{\underset{s}\to\varprojlim}\Bbb Z/N_s\Bbb Z$, where all 
$n_s\in[0,N_s)$ and if  
$D^{\infty }_{a\bar n}=
\mathbin{\underset{s}\to\varprojlim} 
D^{(s)}_{a,n_s\operatorname{mod}N_s}$,  for $a\in\Bbb Z(p)$, then 
$\pi _{K{\ur }}(D^{\infty }_{a\bar n})=0$ if 
$\bar n\notin\Bbb Z_{\geqslant 0}\subset\varprojlim\Bbb Z/N_s\Bbb Z$,  
and $\pi _{K{\ur }}(D^{\infty }_{an})=t^{ap^n-1}\operatorname{d}(t)$ if 
$\bar n=n\in\Bbb Z_{\geqslant 0}$. 

Let $\pi _{K\infty }:=\pi _{K{\ur }}|_{\Cal M_{K\infty }}$. Then one
can easily prove the following proposition. 

\proclaim{Proposition 3.2} {\rm 1)} $\pi _{K\infty }:
\Cal M_{K\infty }\longrightarrow\hat\Omega ^1_{O_{K_{\ur }}}$ is a
continuous epimorphism of $k(p)$-modules; 
\medskip 
{\rm 2)} $\ker\pi _{K\infty }$ is the 
$k(p)$-submodule in $\Cal M_{K\infty }$ topologically generated by 
all $D^{\infty }_{an}$ with $n<0$.
\endproclaim 
\medskip 

\subsubhead {\rm 3.3.} Admissible systems of group morphisms 
\endsubsubhead 

Suppose $K'=k((t'))\subset K(p)$ has the same residue field as
$K$. 
Using $K'$ instead of $K$ 
 we can introduce analogues $\Cal M_{K's}$, $\bar\Cal M_{K's}$, 
$\Cal M_{K'\infty }$, etc. of $\Cal M_{Ks}$, $\bar\Cal M_{Ks}$, $\Cal
M_{K\infty }$, etc. 

\definition{Definition} $f_{KK'}=\{f_{KK's}\}_{s\geqslant 0}$ is a family 
of continuous morphisms of $\Bbb F_p$-modules 
$f_{KK's}:\Cal M_{Ks}\longrightarrow\Cal M_{K's}$ which are always assumed to be
compatible, i.e. for all $s\geqslant 0$, 
$f_{KK',s+1}j'_s=j_sf_{KK's}$. Here $j_s:\Cal
M_{K,s+1}\longrightarrow\Cal M_{Ks}$ and 
$j'_s:\Cal M_{K',s+1}\longrightarrow\Cal M_{K's}$ are 
connecting morphisms.  
\enddefinition 

We shall denote the $k(p)$-linear extension of $f_{KK's}$ by the same
symbol $f_{KK's}$. Set 
$$f_{KK'\ur }:=\mathbin{\underset{s}\to\varprojlim} f_{KK's}:\bar\Cal M_{K{\ur }}
\longrightarrow \bar\Cal M_{K'{\ur }}.$$

\definition{Definition} With the above notation 
$f_{KK'}$ is called admissible if:
\medskip 
{\bf A1}.  There is a continuous $k(p)$-linear isomorphism 
$f_{KK'\infty }:\hat\Omega ^1_{O_{K_{\ur }}}\longrightarrow 
\hat\Omega ^1_{O_{K'_{\ur }}}$ such that 
$f_{KK'\ur }\pi _{K'{\ur }}=\pi _{K{\ur }}f_{KK'\infty }$; 
\medskip 
{\bf A2}.  $f_{KK'\infty }$ commutes with the Cartier operators $C$ and $C'$  
on $\hat\Omega ^1_{O_{K_{\ur }}}$ and, resp., $\hat\Omega ^1_{O_{K'_{\ur }}}$;
\medskip 
{\bf A3}.  For all $m\in\Bbb N$, 
$f_{KK'\infty }\left (t^m\hat\Omega ^1_{O_{K_{\ur }}}\right )\subset 
t^{\prime m}\hat\Omega ^1_{O_{K'_{\ur }}}.$
\enddefinition 

\remark{Remark} Recall that the Cartier operator 
$C:\hat\Omega ^1_{O_{K_{\ur }}}\longrightarrow\hat\Omega ^1_{O_{K_{\ur }}}$ 
is uniquely determined by the following properties:
\newline  
a) $C(\operatorname{d}(\hat O_{K_{\ur }}))=0$;
\newline  
b) if $f\in t\hat O_{K_{\ur }}$ then 
$C(f^p\operatorname{d}(t)/t)=f\operatorname{d}(t)/t$. 
\newline  
It can be shown that the definition of $C$ does not depend on the
choice of the uniformiser $t$, 
$C$ is $\sigma ^{-1}$-linear and $\Ker C=\operatorname{d}(\hat O_{K_{\ur }})$.
\endremark 

The following properties of admissible systems
$f_{KK'}=\{f_{KK's}\}_{s\geqslant 0}$
follow directly from the above definition:

\roster 

\item  the map $f_{KK'\infty }$ is uniquely determined by
$f_{KK'\ur }$;
\medskip 

\item  if $K''=k((t''))\subset K(p)$ and
$g_{K'K''}=\{g_{K'K''s}\}_{s\geqslant 0}$ is admissible then 
so is the composition 
$(fg)_{KK''}:=\{f_{KK's}g_{K'K''s}\}_{s\geqslant 0}$ and 
it holds $(fg)_{KK''\infty }=f_{KK'\infty}
g_{K'K''\infty }$;
\medskip 

\item \ $f_{KK'\infty }(\operatorname{d}\hat O_{K_{\ur }})\subset 
\operatorname{d}\hat O_{K'_{\ur }}$; 
\medskip 

\item \  for all $a,b\in\Bbb Z(p)$ and $m\in\Bbb Z_{\geqslant 0}$,
there are unique $\alpha ^{\infty }_{a,b,-m}(f_{KK'})\in k(p)$ such that if
$n\geqslant 0$ then 
$$f_{KK'\infty }\left (t^{ap^n}\frac{\operatorname{d}(t)}{t}\right )=
\sum\Sb 0\leqslant m\leqslant n\endSb 
\sigma ^n\alpha ^{\infty }_{a,b,-m}(f_{KK'})
t^{\prime bp^{n-m}}\frac{\operatorname{d}(t')}{t'};\tag{$3.1$}$$
\medskip 

\item \ the above coefficients $\alpha ^{\infty }_{a,b,-m}(f_{KK'})$ 
satisfy the following property: 
if $b/p^m<a$ then $\alpha ^{\infty }_{a,b,-m}(f_{KK'})=0$.
\endroster 
\medskip 

\definition{Definition} With the above notation 
an admissible compatible system $f_{KK'}$ 
will be called special admissible if 
$f_{KK'\ur }(\Cal M_{K\infty })\subset\Cal M_{K'\infty }$. 
\enddefinition 

Notice that the composition of special admissible systems is again 
special admissible. 

\subsubhead{\rm 3.4.} Characterisation of special admissible systems 
\endsubsubhead 

Let $f_{KK'}=\{f_{KK's}\}_{s\geqslant 0}$ be a compatible system. Then for any 
$s\geqslant 0$, the $k(p)$-linear morphism 
$f_{KK's}:\bar\Cal M_{Ks}\longrightarrow\bar\Cal M_{K's}$ 
is defined over $\Bbb F_p$, i.e. it 
comes from a $\Bbb F_p$-linear 
morphism $f_{KK's}:\Cal M_{Ks}\longrightarrow\Cal M_{K's}$. Therefore, 
in terms of the standard generators $D^{(s)}_{an}$ and $D^{\prime
(s)}_{an}$ (which correspond to the uniformisers 
$t=t_K$ and, resp., $t'=t_{K'}$), we have for any $s\geqslant 0$ and $a\in\Bbb
Z(p)$ 
that 
$$f_{KK's}(D^{(s)}_{a0})=\sum\Sb b\in\Bbb Z(p), m
\in\Bbb Z\operatorname{mod}N_s\endSb 
\alpha _{abm}(f_{KK's})D^{\prime
(s)}_{bm},$$
where all $\alpha _{abm}(f_{KK's})\in k_s\subset k(p)$. 
Notice that for all 
$n\in\Bbb Z\operatorname{mod}N_s$, it  
holds  
$$f_{KK's}(D^{(s)}_{an})=\sum\Sb b\in\Bbb Z(p), m\in\Bbb
Z\operatorname{mod}N_s
\endSb \sigma ^n\alpha
_{abm}(f_{KK's})D^{\prime (s)}_{b,m+n}.$$

\proclaim{Proposition 3.3} Suppose 
$f_{KK'}=\{f_{KK's}\}_{s\geqslant 0}$ is a compatible
system. Then 
it is special admissible if and only if for any $s\geqslant 0$, 
there are $v_s\in\Bbb N$ such that $v_s\to \infty $ if $s\to\infty $,
and 
if $a,b<v_s$, $m\geqslant 0$ and $b/p^m<a$ then 
$\alpha _{a,b,-m\operatorname{mod}N_s}(f_{KK's})=0$. 
\endproclaim 

\demo{Proof} Suppose $f_{KK'}$ is special admissible. 
Then $f_{KK'\ur }(\Cal M_{K\infty })\subset\Cal M_{K'\infty }$ and 
for all $a\in\Bbb Z(p)$ and $n\in\Bbb Z$, 
$$f_{KK'\ur }(D^{\infty }_{an})=\sum\Sb b\in\Bbb Z(p), m\in\Bbb Z\endSb 
\beta _{anbm}D^{\prime\infty }_{b,n+m}.$$
Here all coefficients $\beta _{anbm}\in k(p)$ and because $f_{KK'\ur }$ commutes with
$\sigma $, 
there are $\gamma _{abm}\in k(p)$ such that 
$\beta _{anbm}=\sigma ^n(\gamma _{abm})$. Therefore, 
if $a,b\in\Bbb Z(p)$, $m\in\Bbb Z$ and $\gamma _{abm}\ne 0$ then $m\leqslant 0$
and 
$\alpha _{abm}^{\infty }(f_{KK'})=\gamma _{abm}$. 

If $s\geqslant 0$, $a\in\Bbb Z(p)$,  
$$f_{KK's}(D^{(s)}_{a0})=
\sum\Sb b\in\Bbb Z(p), m\in\Bbb Z\operatorname{mod}N_s\endSb \alpha
_{a,b,-m}(f_{KK's})D^{\prime (s)}_{b,-m}$$ 
and $b/p^{N_s}<a$ then for any $m\geqslant 0$,  
$\alpha _{a,b,-m\operatorname{mod}N_s}(f_{KK's})=\alpha ^{\infty
}_{a,b,-m}(f_{KK'})$.  
This implies that $\alpha _{a,b,-m\operatorname{mod}N_s}(f_{KK's})=0$ if 
$a,b<p^{N_s}$ and 
$b/p^m<a$. Therefore, we can take $v_s=p^{N_s}$. 
This proves the \lq\lq only if \rq\rq\ part of the proposition. 

Suppose now that $v_s\to \infty $ if $s\to\infty $ and for $a,b\in\Bbb Z(p)$,
$m\geqslant 0$, 
$$\alpha _{a,b,-m\operatorname{mod}N_s}(f_{KK's})=0$$ 
if $a,b<v_s$ and $b/p^m<a$. 
If in addition $p^{N_s}>b$ then 
$\alpha _{a,b, -m\operatorname{mod}N_s}(f_{KK's})$ does not depend on $s$
and can be denoted by $\alpha ^{\infty }_{a,b,-m}$. 
Clearly, $\alpha ^{\infty }_{a,b,-m}=0$ if $b/p^m<a$. 
Let $a\in\Bbb Z(p)$ and 
$$d=f_{KK'\ur }(D^{\infty }_{a0})-\sum\Sb b\in\Bbb Z(p), m\geqslant 0\endSb 
\alpha ^{\infty
}_{a,b,-m}D^{\prime\infty }_{b,-m}.$$ 
Let $s\geqslant 0$ and let 
$d_s\in\bar\Cal M_{Ks}$ be the image of $d$ 
under the natural projection 
$\bar\Cal M_{K{\ur }}\longrightarrow\bar\Cal M_{Ks}$.  
If $s_1\geqslant s$ then the corresponding projection 
$d_{s_1}\in\bar\Cal M_{K_{s_1}}$ 
is a linear combination of $D^{(s_1)}_{bm}$ with
$b>p^{N_{s_1}}$. Therefore, $d_s$ also does not contain the terms 
$D^{(s)}_{bm}$ for which $b>p^{N_{s_1}}$. Because  
$\lim _{s_1\to\infty }N_{s_1}=\infty $, this implies that 
 $d_s=0$ for all $s\geqslant 0$ and, therefore, 
$d=0$. So, $f_{KK'\ur }(\Cal M_{K\infty })\subset \Cal M_{K'\infty }$.

Set $\alpha ^{\infty }_{a,b,-m}(f_{KK'}):=\alpha^{\infty } _{a,b,-m}$ and 
define 
$f_{KK'\infty }:\hat\Omega ^1_{O_{K_{\ur }}}\longrightarrow\hat\Omega
^1_{O_{K'_{\ur }}}$ by formula (3.1). 
It is easy to see that $f_{KK'\infty }$ satisfies the 
requirements {\bf A1-A3} from the definition of admissible system in n.3.3. 
This proves the \lq\lq if\rq\rq\ part of our proposition.  

\enddemo 

\remark{Remark} Any special admissible $f_{KK'}$ can be defined as 
a $k(p)$-linear isomorphism 
$f_{KK'\ur }:\Cal M_{K\infty }\longrightarrow\Cal
M_{K'\infty }$ such that 

\roster 

\item \ $f_{KK'\ur }$ commutes with $\sigma $;
\medskip 
\item  \ if $a\in\Bbb Z(p)$ then  
$$f_{KK'\ur }(D^{\infty }_{a0})=
\sum\Sb b\in\Bbb Z(p), m\geqslant 0\endSb \alpha _{a,b,-m}
D^{\prime\infty }_{b,-m}$$ 
where 
$\alpha _{a,b,-m}=0$ if $b/p^m<a$.
\endroster 
\endremark 

\subsubhead{\rm 3.5.} Analytic compatible systems 
\endsubsubhead 

Suppose  $K,K'\subset K(p)$. Then the corresponding residue fields 
$k$ and $k'$ are subfields of the residue field 
$k(p)\subset \bar\Bbb
F_{q_0}$. Therefore, if $K\simeq K'$ then $k=k'$ and we can introduce 
the set $\Iso ^0(K,K')$ of 
field isomorphisms $\eta :K\longrightarrow K'$ such that 
$\eta |_k=\id $. Notice that any  $\eta \in\Iso ^0(K,K')$ induces a $k(p)$-linear map 
$\Omega ^1(\eta ):\hat\Omega ^1_{O_{K_{\ur }}}\longrightarrow\hat\Omega ^1_{O_{K_{\ur
}'}}$. 

For all $s\geqslant 0$, any $\eta\in\Iso ^0(K,K')$ can be naturally 
extended  to $\eta _s\in\Iso ^0(K_s,K'_s)$. Then 
$\eta ^*_{KK'}=\{\eta ^*_s\}_{s\geqslant 0}$ is a compatible 
system and $\eta _{KK'\infty }=\Omega ^1(\eta )$.  
Propositions 2.4 and 3.3 imply that 
$\eta ^*_{KK'}$ is a special admissible system. 

Consider the opposite situation. Choose a uniformiser 
$t_K$ in $K$ and introduce 
$\Fr (t_K)\in\Aut (K_{\ur })$ such that 
$\Fr (t_K):t_K\mapsto t_K$ and $\Fr (t_K)|_{k(p)}=\sigma $. Then for 
all $s\geqslant 0$, $\Fr (t_K)$ induces an 
automorphism of $K_s$ which will be denoted  
by $\Fr (t_K)_s$. Then 
$\Fr (t_K)^*=\{\Fr (t_K)_s\}_{s\geqslant 0}$ 
is a compatible system, but this system is not admissible: 
the corresponding map $\Fr (t_K)_{\infty }$ coincides 
with the Cartier operator and, therefore, is not $k(p)$-linear.  

More generally, consider a compatible system 
$\theta _{KK'}=\{\theta _{KK's}\}_{s\geqslant 0}$ 
where for all $s\geqslant 0$, 
$\theta _{KK's}=\theta _s^*$ and $\theta _s\in\Iso (K_s,K'_s)$. 
Then after choosing a uniformising element 
$t_{K'}$ in $K'$ we have  
$\theta _s=\eta _s\Fr (t_{K'})^{n_s}$, for all $s\geqslant 0$, 
where $\eta _s\in \Iso ^0
(K_s,K'_s)$ and $n_{s+1}\equiv n_s\operatorname{mod}N_s$. If 
$\bar n=\mathbin{\underset{s}\to\varprojlim}n_s\in
\mathbin{\underset{s}\to\varprojlim}\Bbb Z/N_s\Bbb Z$ then 
$\theta _{KK'}$ is the composite of the special admissible system 
$\{\eta _s^*\}_{s\geqslant 0}$ and the system 
$\Fr (t_{K'})^{\bar n*}$ which is special admissible 
if and only if $\bar n=0$. Therefore, 
$\theta _{KK'}$ is special admissible if and only if it comes from 
a compatible system of field isomorphisms 
$\eta _s\in\Iso ^0(K_s,K'_s)$.  
\medskip 

\subsubhead {\rm 3.6.} Locally analytic systems
\endsubsubhead 

\definition{Definition} If $f_{KK'}$ is an admissible system,  
then $f_{KK'\an }:=f_{KK'\infty }|_{\operatorname{d}(\hat O_{K_{\ur }})}$. 
\enddefinition 

\remark{Remark} Notice the following similarity to the definition 
of $f_{\an }$ for $f\in\Aut\Cal M$ from n.2.3. If 
$f_{KK}=\{f_{KKs}\}_{s\geqslant 0}$ is any admissible system 
then 
$g_{KK}:=\{f_{KKs\an }\}_{s\geqslant 0}$ is also 
admissible and $f_{KK\an }=g_{KK\an }$.
\endremark 
\medskip 

\definition{Definition} An admissible system 
$f_{KK'}=\{f_{KK's}\}_{s\geqslant 0}$ will be called locally analytic if
for 
any $s\geqslant
0$, 
there are $v_s\in\Bbb N$ and $\eta _s\in\Iso ^0(K,K')$ such that 
$v_s\to +\infty $ as $s\to\infty $ and $f_{KK'\an }\equiv
\operatorname{d}(\eta _s)\hat\otimes _kk(p)\operatorname{mod}t^{\prime v_s}$. 
\enddefinition 

\proclaim{Proposition 3.4} 
Suppose that 
$f_{KK'}=\{f_{KK's}\}_{s\geqslant 0}$ is special admissible and locally
analytic. 
Then there is an $\eta\in\Iso ^0(K,K')$ such that 
$f_{KK'\an }=\operatorname{d}(\eta )\hat\otimes _kk(p)$. 
\endproclaim 

\demo{Proof} If $s\geqslant 0$ and $a,b\in\Bbb Z(p)$ are such that 
$v_s/p^{N_0}<a,b<v_s$, then 

$$\alpha ^{\infty }_{ab0}(f_{KK'})=\alpha _{ab0}(\eta _s^*)=\alpha
_{ab0}(f_{KK's})=\alpha _{ab0}(f_{KK'0})\in k.$$ 
Therefore, by Proposition 2.7, all conjugates of $\eta _s$ over $K$ 
are congruent 
modulo $t^{\prime v_s(1-p^{-N_0})/\delta _p}$, and  
$\eta _s(t)\in k[[t']]\operatorname{mod}t^{\prime
v_s(1-p^{-N_s})/\delta _p}$, where 
$\delta _p$ is 1 if $p\ne 2$ and $\delta _p=2$ if $p=2$. 
This implies that $\alpha _{ab0}(f_{KK's})\in k$ if 
$a,b<v_s(1-p^{-N_s})/\delta _p$.

If $b<p^{N_s}$ then $\alpha _{ab0}(f_{KK's})=
\alpha ^{\infty }_{ab0}(f_{KK'})$. 
So, $\alpha ^{\infty }_{ab0}(f_{KK'})\in k$ if 
$b<c_s:=$
\linebreak 
$\min \left\{p^{N_s}, v_s(1-p^{-N_s})/\delta _p\right\}$. 
But $c_s\to\infty $ if $s\to\infty $ and, therefore, 
$\alpha ^{\infty }_{ab0}(f_{KK'})\in k$ for all $a,b\in\Bbb Z(p)$. 

As we have already noticed, if $b<\min\{p^{N_s}, v_s\}$ then 
$$\alpha _{ab0}(f_{KK's})=\alpha _{ab0}(\eta
_s^*)=\alpha ^{\infty }_{ab0}(f_{KK'}).$$
Therefore, by Proposition 2.7  there exists  
$\mathbin{\underset{s}\to\varprojlim} \eta _s:=\eta\in \Iso ^0(K,K')$ and 
$f_{KK'\an }=\operatorname{d}(\eta )\hat\otimes _kk(p)$. 

The proposition is proved. 
\enddemo 

\subsubhead{\rm 3.7.} Comparability of admissible  systems  
\endsubsubhead  

With the above notation 
suppose $L,L'$ are finite field extensions of 
$K$, resp. $K'$, in $K(p)$. Let $g_{LL'}=\{g_{LL's}\}_{s\geqslant 0}$ 
be a compatible family of
continuous field 
isomorphisms  
$g_{LL's}:L_s\longrightarrow L'_s$. 
Then the natural embeddings 
$\Gamma _L(p)\subset \Gamma _K(p)$ and 
$\Gamma _{L'}(p)\subset\Gamma _{K'}(p)$ 
induce embeddings 
$\Gamma _{L_s}(p)\subset\Gamma _{K_s}(p)$ and 
$\Gamma _{L'_s}(p)\subset\Gamma
_{K'_s}(p)$,  for any $s\geqslant 0$.  

\definition{Definition} With the above assumptions the systems 
$g_{LL'}$ and $f_{KK'}$ will be called comparable if, for 
all $s\geqslant 0$, there is the following commutative diagram  

$$
\CD 
\Cal M_{Ls}  @>{g_{LL's}}>>    \Cal M_{L's} \\ 
@V{j_s}VV                   @VV{j_s'}V \\ 
\Cal M_{Ks}  @>{f_{KK's}}>>    \Cal M_{K's} 
\endCD 
$$
where the vertical arrows $j_s$ and $j'_s$ are induced by the 
embeddings $\Gamma _{L_s}(p)\subset\Gamma _{K_s}(p)$ and, resp.,  
$\Gamma _{L'_s}(p)\subset\Gamma _{K'_s}(p)$. 
\enddefinition 

If $g_{LL'}$ and $f_{KK'}$ are comparable then we have 
the following commutative diagram  

$$
\CD 
\bar\Cal M_{L{\ur }}  @>{g_{LL'\ur }}>>    \bar\Cal M_{L'{\ur }} \\ 
@V{j_{\ur }}VV                   @VV{j_{\ur }'}V \\ 
\bar\Cal M_{K{\ur }}  @>{f_{KK'\ur }}>>    \bar\Cal M_{K'{\ur }} 
\endCD 
\tag{3.2}$$
where $j_{\ur }:=\mathbin{\underset{s}\to\varprojlim} j_s\hat\otimes _{k_s}k(p)$ 
and $j'_{\ur }:=\mathbin{\underset{s}\to\varprojlim} j'_s\hat\otimes _{k_s}k(p)$. 
Notice that $j_{\ur }$ and $j'_{\ur }$ are epimorphic. Indeed, 
let $U_{L_s}$, $U_{K_s}$ be principal units in 
$L_s$, resp. $K_s$. Then 
$\Cal M_{L{\ur }}=\mathbin{\underset{s}\to\varprojlim} 
U_{L_s}/U_{L_s}^p$ and 
$\Cal M_{K{\ur }}=\mathbin{\underset{s}\to\varprojlim} 
U_{K_s}/U_{K_s}^p$ contain as dense
subsets the images of the groups of principal 
units $U_{L_{\ur }}$, resp. $U_{K_{\ur }}$, of the fields 
$L_{\ur }$, resp. $K_{\ur }$. By class field theory,  
$j_{\ur }$ is induced by the norm map $N=N_{L_{\ur }/K_{\ur }}$ from 
$L^*_{\ur }$ to $K^*_{\ur }$. By [Iw, Ch.2],  
$N(U_{L_{\ur }})$ is dense in $U_{K_{\ur }}$ and, therefore,  $j_{\ur }$ 
(together with $j'_{\ur }$) is surjective. 

Suppose $L/K$ and $L'/K'$ are Galois extensions. 
Denote their inertia subgroups
by 
$I_{L/K}$ and $I_{L'/K'}$. Then we have identifications 
$I_{L/K}=\Gal (L_{\ur }/K_{\ur })$ and 
$I_{L'/K'}=\Gal (L'_{\ur }/K'_{\ur })$. 

Consider the following condition:

\proclaim{C} There is a group isomorphism 
$\kappa :I_{L/K}\longrightarrow I_{L'/K'}$ 
such that for any $\tau\in I_{L/K}$, 
$\tau ^*_{LL{\ur }}g_{LL'\ur }=g_{LL'\ur}
\kappa (\tau )_{L'L'{\ur }}^*$. 
\endproclaim 
\medskip

\proclaim{Proposition 3.5} Suppose $g_{LL'}$ and $f_{KK'}$ are
comparable and $g_{LL'}$ satisfies the above condition {\bf C}. 
If $g_{LL'}$ is admissible then $f_{KK'}$ is also admissible. 
\endproclaim 

\demo{Proof} Because $g_{LL'}$ is admissible we have the following
commutative diagram  

$$
\CD 
\bar\Cal M_{L{\ur }}  @>{g_{LL'\ur }}>>    \bar\Cal M_{L'{\ur }} \\ 
@V{\pi _{L{\ur }}}VV                   @VV{\pi _{L'{\ur }}}V\\ 
\hat\Omega ^1_{O_{L_{\ur }}}  @>{g_{LL'\infty }}>>    \hat\Omega ^1_{O_{L'_{\ur }}} 
\endCD 
\tag{3.3}$$

If $\tau\in I_{L/K}\subset\Aut ^0(L_{\ur })$ then it follows from the definition
of $\pi _{L_{\ur }}$ that 
$$\tau ^*\pi _{L{\ur }}=\pi _{L{\ur }}\Omega (\tau ). \tag{3.4}$$
This means that $\pi _{L{\ur }}$ transforms the natural action 
of $I_{L/K}$ on $\bar\Cal M_{L{\ur }}$ into 
the natural action of $I_{L/K}$ on $\hat\Omega ^1_{O_{L_{\ur }}}$. 
Because $j_{\ur }$ is induced by the norm map of the field extension 
$L_{\ur }/K_{\ur }$, this gives us the following commutative diagram

$$
\CD 
\bar\Cal M_{L{\ur }}  @>{\pi _{L{\ur }}}>>    \hat\Omega ^1_{O_{L_{\ur }}} \\ 
@V{j_{\ur }}VV                   @VV{\Tr}V \\ 
\bar\Cal M_{K{\ur }}  @>{\pi _{K{\ur }}}>>    \hat\Omega ^1_{O_{K_{\ur }}} 
\endCD 
\tag{3.5}$$
where $\Tr $ is induced by the trace of the extension 
$L_{\ur }/K_{\ur }$. Similarly, we have the commutative diagram  

$$
\CD 
\bar\Cal M_{L'{\ur }}  @>{\pi _{L'{\ur }}}>>    \hat\Omega ^1_{O_{L'_{\ur }}} \\ 
@V{j'_{\ur }}VV                   @VV{\Tr '}V \\ 
\bar\Cal M_{K'{\ur }}  @>{\pi _{K'{\ur }}}>>    \hat\Omega ^1_{O_{K'_{\ur }}} 
\endCD 
\tag{3.6}$$

We have already seen that 
$\pi _{L{\ur }}$, $\pi _{L'{\ur }}$, 
$j_{\ur }$ and $j'_{\ur }$ are surjective. The traces 
$\Tr $ and $\Tr '$ are also surjective. Indeed, 
suppose $t_L$, resp. $t_K$, are uniformising elements 
for $L$, resp. $K$. Then 
$$\hat\Omega ^1_{O_{L_{\ur }}}=
\{f\operatorname{d}(t_L)\ |\ f\in \hat O_{L_{\ur }}\}
=\{g\operatorname{d}(t_K)\ |\ g\in \Cal D(L/K)^{-1}\hat O_{L_{\ur }}\},$$
where $\Cal D(L/K)$ is the different of the extension $L/K$. 
It remains to notice that 
$\Tr (\Cal D(L/K)^{-1}\hat O_{L_{\ur }})=\hat O_{K_{\ur }}$.

Because $g_{LL'}$ and $f_{KK'}$ are comparable, we have the following
commutative diagram  

$$
\CD 
\bar\Cal M_{L{\ur }}  @>{g_{LL'\ur }}>>   \bar\Cal M_{L'{\ur }} \\ 
@V{j_{\ur }}VV                   @VV{j '_{\ur }}V \\ 
\bar\Cal M_{K{\ur }}  @>{f_{KK'\ur }}>>    \bar\Cal M_{K'{\ur }}
\endCD 
\tag{3.7}$$

Suppose $\omega _K\in\hat\Omega ^1_{O_{K_{\ur }}}$. As it has  
been proved there is an $\omega _L\in\hat\Omega ^1_{O_{L_{\ur }}}$
such that 
$$\Tr (\omega _L)=\sum\Sb \tau\in I_{L/K}\endSb 
\Omega (\tau )(\omega _L)=\omega _K.$$
Then 
$$g_{LL'\infty }(\omega _K)=\sum\Sb \tau\in I_{L/K}\endSb 
g_{LL'\infty }(\Omega (\tau )(\omega _L))\tag{3.8}$$
$$=\sum\Sb \tau '\in I_{L'/K'}\endSb 
\Omega (\tau ')(g_{LL'\infty }(\omega _L))=
\Tr '(g_{LL'\infty }(\omega _L))\in\hat\Omega ^1_{O_{K'_{\ur }}}$$
because 
$\Omega (\tau )g_{LL'\infty }=g_{LL'\infty }\Omega (\kappa (\tau ))$, 
for any $\tau\in I_{L/K}$.  
 This equality 
is implied by the following computations (we use the commutative
diagrams (3.3), (3.4) and condition {\bf C})

$$\pi _{L{\ur }}\Omega (\tau )g_{LL'\infty }=\tau ^*\pi _{L{\ur
}}g_{LL'\infty }=\tau ^*g_{LL'\ur }\pi _{L'{\ur }}\qquad\qquad $$
$$=
g_{LL'\ur }\kappa (\tau )^*\pi _{L'{\ur }}=
g_{LL'\ur }\pi _{L'{\ur }}\Omega (\kappa (\tau ))=
\pi _{L{\ur }}g_{LL'\infty }\Omega (\kappa (\tau )),$$ 
because $\pi _{L_{\ur }}$ is surjective. 

Let $f_{KK'\infty }$ be the restriction of 
$g_{LL'\infty }$ on $\hat\Omega ^1_{O_{K_{\ur }}}$. Then 
formula (3.8)  implies that 
$f_{KK'\infty }(\hat\Omega ^1_{O_{K_{\ur }}})
\subset \hat\Omega ^1_{O_{K'_{\ur }}}$ and we have the following
commutative diagram

$$
\CD 
\hat\Omega ^1_{O_{L_{\ur }}}  @>{g _{LL'\infty }}>>    \hat\Omega ^1_{O_{L'_{\ur }}} \\ 
@V{\Tr }VV                   @VV{\Tr '}V \\ 
\hat\Omega ^1_{O_{K_{\ur }}}  @>{f_{KK'\infty }}>>    \hat\Omega ^1_{O_{K'_{\ur }}} 
\endCD
\tag{3.9}$$

We now verify that $f_{KK'\infty }$ satisfies the requirements {\bf A1-A3} from
n.3.3.  

Property {\bf A1} means that we have the following commutative
diagram

$$
\CD 
\bar\Cal M_{K_{\ur }}  @>{f_{KK'\ur }}>>    \bar\Cal M_{K'_{\ur }} \\ 
@V{\pi _{K_{\ur }}}VV                   @VV{\pi _{K'_{\ur }}}V \\ 
\hat\Omega ^1_{O_{K_{\ur }}}  @>{f_{KK'\infty }}>>    \hat\Omega ^1_{O_{K'_{\ur }}} 
\endCD 
$$

Its commutativity 
is implied by the following
computations (we use commutative diagrams (3.2), (3.5), (3.3) and (3.9)) 
$$j_{\ur }f_{KK'\ur }\pi _{K'{\ur }}=
g_{LL'\ur }j'_{\ur }\pi _{K'{\ur }}=
g_{LL'\ur }\pi _{L'{\ur }}\Tr '$$
$$=\pi _{L{\ur }}g_{LL'\infty }\Tr '=
\pi _{L{\ur }}\Tr f_{KK'\infty }=
j_{\ur }\pi _{K{\ur }}f_{KK'\infty }$$
because $j_{\ur }$ is surjective. 

Let $C_K$, $C_{K'}$, $C_L$ and $C_{L'}$ be the Cartier 
operators on, resp., $\hat\Omega ^1_{O_{K_{\ur }}}$, 
 $\hat\Omega ^1_{O_{K'_{\ur }}}$,
 $\hat\Omega ^1_{O_{L_{\ur }}}$ and 
 $\hat\Omega ^1_{O_{L'_{\ur }}}$. 
Clearly, $C_L\Tr =\Tr C_K$ and 
$C_{L'}\Tr '=\Tr 'C_{K'}$. Then it follows 
from the commutative diagram  
(3.9) and property {\bf A2} for $g_{LL'\infty }$ that 
$$\Tr C_Kf_{KK'\infty }=
C_L\Tr f_{KK'\infty }=
C_Lg_{LL'\infty }\Tr $$
$$=g_{LL'\infty }C_{L'}\Tr =
g_{LL'\infty }\Tr C_{K'}=
\Tr f_{KK'\infty }C_{K'}.$$
Property {\bf A2} for $f_{KK'\infty }$ follows 
because $\Tr $ is surjective. 

By condition {\bf C}, the ramification indices $e$ and $e'$ 
of the extensions 
$L_{\ur }/K_{\ur }$ and $L'_{\ur }/K'_{\ur }$ are equal. Then we use 
the condition {\bf A3} for $g_{LL'\infty }$ to deduce that for any 
$n\geqslant 0$,  
$$g_{LL'\infty }(t_K^n\hat\Omega ^1_{O_{L_{\ur }}})=
g_{LL'\infty }(t_L^{en}\hat\Omega ^1_{O_{L_{\ur }}})=
t_L^{\prime e'n}\hat\Omega ^1_{O_{L'_{\ur }}}=
t_{K'}^n\hat\Omega
^1_{O_{L'_{\ur }}}.$$

Therefore, it follows from the commutativity of diagram (3.9) that 

$$t^n_{K'}\hat\Omega ^1_{O_{K'_{\ur }}}=t_{K'}^n
\Tr '(\hat\Omega ^1_{O_{L'_{\ur }}})=\Tr '(g_{LL'\infty
}(t_K^n\hat\Omega ^1_{O_{L_{\ur }}}))$$
$$=f_{KK'\infty }(\Tr (t_K^n\hat\Omega ^1_{O_{L_{\ur }}}))=f_{KK'\infty
}(t_{K'}^n\hat\Omega ^1_{O_{K_{\ur }}}).$$

The proposition is proved. 
\enddemo 

\remark{Remark} Using the embeddings of the Galois groups 
$\Gamma _{L_s}(p)$ and $\Gamma _{K_s}(p)$ into their Magnus's algebras 
from n.1.3, one can prove in addition that if $g_{LL'}$ is special then 
$f_{KK'}$ is also special. In other words, 
under condition {\bf C}, $j_{\ur }(\Cal M_{L\infty })\subset\Cal
M_{K\infty }$. 
\endremark 
\medskip 

Suppose $g_{LL'}$ and $f_{KK'}$ are comparable systems. 
Suppose also that 
$g_{LL'}$ and $f_{KK'}$ are special admissible, 
locally analytic and satisfy  
condition {\bf C}. Then there are 
$\eta _{LL' }\in\Iso ^0(L,L')$ and 
$\eta _{KK' }\in\Iso ^0(K,K')$ such that 
$f_{KK'\infty }|_{\operatorname{d}\hat O_{K_{\ur }}}
=\operatorname{d}(\eta _{KK' })\hat\otimes _kk(p)$ and 
$g_{LL'\infty }|_{\operatorname{d}\hat O_{L_{\ur }}}=
\operatorname{d}(\eta _{LL' })\hat\otimes _{k_L}k_L(p)$. 

\proclaim{Proposition 3.6} With the above notation and assumptions, 
$\eta _{LL' }|_K=\eta _{KK' }$. 
\endproclaim 

\demo{Proof} Clearly, 
for any $\tau\in I_{L/K}$, condition {\bf C} implies that 
$\tau ^*_{LL\infty }g_{LL'\infty }=g_{LL'\infty }\kappa (\tau
)^*_{L'L'\infty }$. 
Restricting this equality to $\operatorname{d}\hat O_{L_{\ur }}$, we obtain 
$$\operatorname{d}(\tau )\operatorname{d}(\eta _{LL' })
=\operatorname{d}(\eta _{LL' })\operatorname{d}(\kappa
(\tau )).$$
Then it follows from proposition 2.7 that  
$\tau\eta _{LL' }=\eta _{LL' }\kappa (\tau )$. 
Therefore, $\eta _{LL' }|_K$ induces 
a ring isomorphism from $\hat O_{K_{\ur }}$ onto $\hat O_{K'_{\ur }}$. 

Suppose $a\in\Tr (\hat O_{L_{\ur }})\subset \hat O_{K_{\ur }}$. If 
$a=\Tr (b)$ with $b\in\hat O_{L_{\ur }}$ then it follows from 
diagram (3.9) and condition {\bf C} that 

$$\operatorname{d}(\eta _{KK' }(a))=\Tr '(\operatorname{d}
(\eta _{LL' }(b)))=
\sum\Sb \tau '\in I_{L'/K'}\endSb 
\operatorname{d}(\tau ')\left (\operatorname{d}
(\eta _{LL' }(b))\right )$$
$$=\sum\Sb \tau\in I_{L/K}\endSb \operatorname{d}(\eta _{LL' 
})(\operatorname{d}(\tau (b)))=
\operatorname{d}\eta _{LL' }(\operatorname{d}a)=\operatorname{d}(\eta _{LL'}(a)).$$
Therefore, for a sufficiently large $M\in\Bbb N$, 
$\operatorname{d}\left (\eta _{LL' }|_K\right )$ and 
$\operatorname{d}\eta _{KK' }$ coincide on 
$t_K^M\hat O_{K_{\ur }}$. Then proposition 2.7 implies that 
$\eta _{LL' }|_K=\eta _{KK' }$.

The proposition is proved.
\enddemo 
\medskip

\subhead 4. Explicit description of the ramification ideals $\Cal
A^{(v)}\operatorname{mod}J^3$ 
\endsubhead 
\medskip 

We return to the notation from n.1. In particular, 
$\Cal A$ is the $\Bbb Z_p$-algebra from n.1.2, 
$\Cal J$ is its augmentation ideal, $\Cal A_k=\Cal A\otimes W(k)$, 
$\Cal J_k=\Cal J\otimes W(k)$, $\Cal A_K=\Cal A\otimes O(K)$,  
etc.  
are the corresponding extensions of scalars, $e\in\Cal A_K$
is the element introduced in n.1.3. We fix an  
$f\in\Cal A_{K(p)}$ such that 
$\sigma f=fe$ and denote the embedding $\psi _f:
\Gamma (p)\longrightarrow (1+\Cal J)^{\times }$ by $\psi $.

\subsubhead {\rm 4.1.} Ramification filtration on $\Cal A$ 
\endsubsubhead 
For any $v\geqslant 0$, consider the ramification 
subgroup $\Gamma (p)^{(v)}$ of $\Gamma (p)$ in the upper numbering. Denote by 
$\Cal A^{(v)}$ the minimal 2-sided closed ideal  in $\Cal A$ 
containing the elements 
$\psi (\tau )-1$, for all $\tau\in\Gamma (p)^{(v)}$.  Then 
$\{\Cal A^{(v)}\ |\ v\geqslant 0\}$ is a decreasing filtration 
by closed ideals of $\Cal A$. In particular, 
if  
$\Cal A_{CM}^{(v)}\operatorname{mod}\Cal J_{CM}^n$ are the projections of 
$\Cal A^{(v)}$ to $\Cal A_{CM}\operatorname{mod}\Cal J_{CM}^n$,  for
$C,M,n\in\Bbb N$, 
then $\Cal A^{(v)}=\mathbin{\underset{C,M,n}\to\varprojlim}\Cal
A^{(v)}_{CM}\operatorname{mod}\Cal J_{CM}^n$. 
Notice also that the ramification filtration $\{\Gamma (p)^{(v)}\}_{v\geqslant 0}$
is left-continuous, i.e.
$\Gamma (p)^{(v_0)}=\mathbin{\underset{v<v_0}\to\bigcap}\Gamma
(p)^{(v)}$, 
 for any $v_0>0$. This
implies a corresponding analogous property  for  
the filtration $\{\Cal A^{(v)}\ |\ v\geqslant 0\}$ on each finite level, 
i.e. for any $C,M,n\in\Bbb N$, we have the following property.  

\proclaim{Proposition 4.1} For any $C,M,n\in\Bbb N$ and $v_0>0$, there is 
a $0<\delta <v_0$ such that  
$\Cal A_{CM}^{(v)}\operatorname{mod}\Cal J_{CM}^n=\Cal
A_{CM}^{(v_0)}\operatorname{mod}\Cal J_{CM}^n$, 
for any $v\in(v_0-\delta ,v_0)$.  
\endproclaim 

\demo{Proof} This follows directly from the definition of the
ramification filtration and the fact that the field of
definition of each projection $f_{CM}\operatorname{mod}\Cal J_{CM}^n$ 
of $f$ to 
\linebreak 
$\Cal A_{CMK(p) }\operatorname{mod}\Cal J_{CMK(p) }^n$ is a
finite 
extension of $K$,
cf. n.1.3. 
\enddemo 

Notice also that the class field theory implies the following property.

\proclaim{Proposition 4.2} If $v\geqslant 0$ and $\Cal A^{(v)}_k:=
\Cal A^{(v)}\otimes W(k)$ then $\Cal A_k^{(v)}\operatorname{mod}\Cal
J_k^2$ is 
topologically generated by 
all elements $p^sD_{an}$, for $n\in\Bbb Z\operatorname{mod}N_0$, 
$a\in\Bbb Z(p)$, $s\geqslant 0$ and $p^sa\geqslant v$.
\endproclaim 

\subsubhead {\rm 4.2.} The filtration $\Cal A(v)$, $v\geqslant 0$ 
\endsubsubhead 

For any $\gamma\geqslant 0$, introduce 
$\Cal F_{\gamma }\in\Cal A_k$ as follows. 

If $\gamma =0$ let $\Cal F_{\gamma }=D_0$. 

If $\gamma >0$ let  
$$\Cal F_{\gamma }=
p^{v_{\gamma }}a_{\gamma }D_{a_{\gamma }v_{\gamma }}-
\sum\Sb a_1,a_2\in\Bbb Z(p) \\ 
n\geqslant 0 \\ p^n(a_1+a_2)=\gamma \endSb 
p^na_1D_{a_1n}D_{a_2n}
-\sum\Sb a_1,a_2\in\Bbb Z(p) \\ n_1\geqslant 0, n_2<n_1\\
p^{n_1}a_1+p^{n_2}a_2=\gamma \endSb
p^{n_1}a_1[D_{a_1n_1},D_{a_2n_2}].$$
Here the first two terms appear only if $\gamma\in\Bbb N$, and the
corresponding 
$v_{\gamma }\in\Bbb Z_{\geqslant 0}$ and $a_{\gamma }\in\Bbb Z(p)$ are
uniquely determined from the equality $\gamma =p^{v_{\gamma }}
a_{\gamma }$. If $\gamma\notin\Bbb Z$ then the above 
formula for $\Cal F_{\gamma }$ contains only the last sum.

For any $v\geqslant 0$, let $\Cal A(v)$ be the minimal closed ideal in 
$\Cal A$ such that $\Cal F_{\gamma }\in \Cal A_k(v):=\Cal
A(v)\otimes W(k)$,  for all 
$\gamma\geqslant v$. 
Equivalently, $\Cal A_k{(v)}$ is the minimal $\sigma $-invariant
closed ideal of $\Cal A_k$, which 
contains all $\Cal F_{\gamma }$ with $\gamma\geqslant v$. 

\remark{Remark}  a) For any $v\geqslant 0$, $\Cal A^{(v)}
\operatorname{mod}\Cal J^2=\Cal A(v)\operatorname{mod}\Cal J^2$.

b) The filtration $\{\Cal A(v)\ |\ v\geqslant 0\}$ is left-continuous.

c) If $C,M\in\Bbb N$ and $\Cal A_{CM}(v)\operatorname{mod}\Cal J_{CM}^n$ 
is the image of $\Cal A(v)$ in 
$\Cal A_{CM}\operatorname{mod}\Cal J_{CM}^n$,  then 
$\Cal A(v)\operatorname{mod}\Cal J^n=
\mathbin{\underset{C,M}\to\varprojlim}\Cal
A_{CM}(v)\operatorname{mod}\Cal J_{CM}^n$. 
\endremark 

If $\gamma\geqslant v_0\geqslant 0$, denote by $\tilde\Cal F_{\gamma }(v_0)$ the elements in 
$\Cal A_k$ given by the same expressions as $\Cal F_{\gamma }$ but 
with the additional restriction $p^{n_1}a_1, p^{n_1}a_2<v_0$ for all 
degree 2 terms 
$p^{n_1}a_1D_{a_1n_1}D_{a_2n_2}$ or 
$p^{n_1}a_1[D_{a_1n_1},D_{a_2,n_2}]$. Clearly, we have the  
following property.

\proclaim{Proposition 4.3} {\rm a}) $\Cal A(v_0)\operatorname{mod}\Cal
J^3$
is 
the minimal ideal of $\Cal A$ such that $\Cal A_k(v_0)$ is 
generated by all elements $\tilde\Cal F_{\gamma }(v_0)$ with
$\gamma\geqslant v_0$. 
\newline 
{\rm b}) If $\gamma\geqslant 2v_0$, then $\tilde\Cal F_{\gamma }(v_0)=
\gamma D_{a_{\gamma }v_{\gamma }}$. 
\endproclaim 

The following theorem is the main technical result about the structure
of the ramification filtration that we need in this paper.

\proclaim{Theorem B} For any $v\geqslant 0$, $\Cal
A^{(v)}\operatorname{mod}
\Cal J^3=\Cal A(v)\operatorname{mod}\Cal J^3$. 
\endproclaim 

This theorem gives an explicit description of the ramification filtration 
$\{\Cal A^{(v)}\}_{v\geqslant 0}$ on the level of $p$-extensions of 
nilpotent class 2. (On the level of abelian $p$-extensions 
such a description is given by the above Remark a).) 
Theorem B can also be stated in the following equivalent
form, where  
we use the index $M+1$ instead of $M$ to simplify the notation in 
its proof below. 

\proclaim{Theorem B'} Suppose $C\in\Bbb N$, $M\in\Bbb Z_{\geqslant 0}$
and $v_0>0$. If, for all 
$v>v_0$, 
$$\Cal A^{(v)}_{C,M+1}\operatorname{mod}\Cal J_{C,M+1}^3=
\Cal A_{C,M+1}(v)\operatorname{mod}\Cal J_{C,M+1}^3,$$
 then 
$$\Cal A^{(v_0)}_{C,M+1}\operatorname{mod}\Cal J_{C,M+1}^3=\Cal A_{C,M+1}
(v_0)\operatorname{mod}\Cal J_{C,M+1}^3.$$ 
\endproclaim 

Clearly, Theorem B' follows from theorem B. 

Conversely, notice first that,  
for a given $C\in\Bbb N$, $M\geqslant 0$  
and $v\gg 0$, 
$$\Cal A_{C,M+1}^{(v)}\operatorname{mod}\Cal J_{C,M+1}^3=
\Cal A_{C,M+1}(v)\operatorname{mod}\Cal J_{C,M+1}^3=0.$$
 Indeed, this is 
obvious for the ideals $\Cal A_{C,M}(v)$, because they are generated by the elements 
obtained from the above elements $\tilde\Cal F_{\gamma }(v)$  by adding the  restrictions 
$a_1,a_2,a_{\gamma }<C$ and $n_1,n_2,v_{\gamma }\leqslant M$. But then, for 
$\gamma\geqslant 2p^MC$, the conditions 
$p^{n_1}a_1+p^{n_2}a_2=\gamma $ (where $n_2\leqslant n_1$) and 
$p^{v_{\gamma }}a_{\gamma }=\gamma $ are never satisfied. For the 
filtration $\{\Cal A^{(v)}\}_{v\geqslant 0}$, we notice, as earlier, that 
the field of definition $K_{C,M+1,3}(f)$ of the image of $f$ in 
$\Cal A_{C,M+1, K(p)}\operatorname{mod}\Cal J_{C,M+1,K(p)}^3$ 
is of finite degree over the basic field $K$. 
Therefore, 
for $v\gg 0$, the ramification subgroup 
$\Gamma (p)^{(v)}$ acts trivially on $K_{C,M+1,3}(f)$ and 
$\Cal A^{(v)}_{C,M+1}\operatorname{mod}\Cal J_{C,M+1}^3=0$. 

Now we can apply descending transfinite induction on $v\geqslant 0$. 
Let 
$$S_{C,M+1}=\{v\geqslant 0\ |\ \Cal
A^{(v)}_{C,M+1}\operatorname{mod}\Cal J_{C,M+1}^3=
\Cal A_{C,M+1}(v)\operatorname{mod}\Cal J_{C,M+1}^3\}.$$
Then $S_{C,M+1}\ne\emptyset $. Let $v_0=\operatorname{inf}S_{C,M+1}$.  

If $v_0>0$ then $\Cal A^{(v_0)}_{C,M+1}\operatorname{mod}\Cal J_{C,M+1}^3=
\Cal A_{CM}(v_0)\operatorname{mod}\Cal J_{C,M+1}^3$ by Theorem B'. 
By the left-continuity property of both filtrations, there is 
a $\delta\in (0,v_0)$ such that 
$\Cal A_{C,M+1}^{(v)}\operatorname{mod}\Cal J_{C,M+1}^3=
\Cal A_{C,M+1}(v)\operatorname{mod}\Cal J_{C,M+1}^3$ whenever
$v\in(v_0-\delta ,v_0)$. 
So, $v_0=\inf S_{C,M+1}\leqslant v_0-\delta $. This is a
contradiction, 
hence $v_0=0$.  
In this case we have   
$\Cal A^{(0)}_{C,M+1}\operatorname{mod}\Cal J_{C,M+1}^3=
\Cal A_{C,M+1}\operatorname{mod}\Cal J_{C,M+1}^3=
\Cal A_{C,M+1}(0)\operatorname{mod}\Cal J_{C,M+1}^3$. 
This implies that $S_{C,M+1}=\Bbb R_{\geqslant 0}$,  and 
Theorem B is deduced from Theorem B'.

The rest of this section is concerned with a proof 
of Theorem B'.

\subsubhead {\rm 4.3.} Auxiliary results  
\endsubsubhead 

\subsubhead {\rm 4.3.1.} The field $K(N^*,r^*)$
\endsubsubhead 

Suppose $N^*\in\Bbb N$, $q=p^{N^*}$ and $r^*=m^*/(q-1)$, where 
$m^*\in\Bbb Z(p)$. Then there is a field $K_1:=K(N^*,r^*)\subset K_{\sep }$ 
such that 

a) $[K_1:K]=q$;

b) the Herbrand function $\varphi _{K_1/K}(x)$ has only one corner point 
$(r^*,r^*)$;

c) $K_1=k((t_{K_1}))$, where $t_{K_1}^qE(-1,t_{K_1}^{m^*})=t_K$ and  
$E$ is the generalised 
Artin-Hasse exponential  introduced in n.1.4. 
\medskip  

The field $K(N^*,r^*)$ 
appears as a subfield of $K(U)$, where $U^q-U=u^{-m^*}$ and
$u^{q-1}=t_K$. 
It is of degree $q$ over $K$. 
Its construction is 
explained in all detail in [Ab2]. 

\subsubhead {\rm 4.3.2.} Relation between liftings of $K$ and $K_1$
modulo $p^{M+1}$, 
$M\geqslant 0$ 
\endsubsubhead 

Recall that we use the uniformiser $t_K$ in $K$ to construct 
the liftings modulo $p^{M+1}$ of $K$, $O_{M+1}(K)=W_{M+1}(k)((t))$ 
and of $K(p)$, $O_{M+1}(K(p))$, where $t=t_{K,M+1}$. 
We use the uniformiser $t_{K_1}$ from above n.4.3.1 c) 
to construct analogous liftings for $K_1$, 
$O'_{M+1}(K_1)=W_{M+1}(k)((t_1))$ 
and for $K_1(p)\supset K(p)$, $O'_{M+1}(K_1(p))$. (Here $t_1
=t_{K_1,M+1}$ is 
the Teichm\" uller representative of $t_{K_1}$ in 
$W_{M+1}(K_1(p))$.) 

Note that, with the above notation  
the field embedding $K\subset K_1$ does not induce an  
embedding $O_{M+1}(K)\subset O'_{M+1}(K_1)$ for $M\geqslant 1$, 
because the Teichm\" uller 
representative 
$t_1=t_{K_1,M+1}=[t_{K_1}]$ cannot be expressed in terms of 
the Teichm\" uller representative $t=t_{K,M+1}=[t_K]$. This difficulty 
can be overcome as follows. Take $t_K^{p^{M}}$ as a 
uniformising element for $\sigma ^{M}K$ and consider the 
corresponding liftings modulo $p^{M+1}$, 
$O_{M+1}(\sigma ^{M}K)=W_{M+1}(k)((t^{p^M}))$ and 
$O_{M+1}(\sigma ^{M}K(p))\subset O_{M+1}(K(p))$. 
 From the definition of liftings it follows that 
$$O_{M+1}(\sigma ^{M}K)\subset W_{M+1}(\sigma ^{M}K)\subset 
W_{M+1}(\sigma ^{M}K_1)\subset O'_{M+1}(K_1)\subset W_{M+1}(K_1),$$
$$O_{M+1}(\sigma ^{M}K(p))\subset W_{M+1}(\sigma ^{M}K(p))
\subset 
W_{M+1}(\sigma ^MK_1(p))
\qquad\qquad\qquad $$
$$\qquad\qquad\qquad\qquad\qquad\qquad 
\subset O'_{M+1}(K_1(p))\subset W_{M+1}(K_1(p)).$$

\proclaim{Lemma 4.4} With respect to the above embedding 
$O_{M+1}(\sigma ^{M}K)\subset O'_{M+1}(K_1)$ we have 
$$t^{p^{M}}=t_1^{qp^{M}}E(-1,t_1^{m^*})^{p^{M}}.$$ 
\endproclaim 

\demo{Proof} If $V$ is the Verschiebung 
morphism on $W_{M+1}(K_1)$ then  property c) 
from n.4.3.1 is equivalent to 
the relation 
$t\equiv t_1^{qp^{M}}E(-1,
t_1^{m^*})\operatorname{mod}VW_{M+1}(K_1)$. 
Then, for any $s\geqslant 0$, we have 
$$t^{p^s}\equiv t_1^{qp^s}E(-1, t_1^{m^*})^{p^s}\operatorname{mod}
V^{s+1}W_{M+1}(K_1).$$
(Using that for any $w_1,w_2\in W_M(K_1)$, $(Vw_1)(Vw_2)=V^2(F(w_1w_2))$
and $pV(w_1)=V^2(Fw_1)$.)
For $s=M$ we obtain the statement of the lemma.
\enddemo 

\subsubhead {\rm 4.3.3.} A criterion 
\endsubsubhead

Consider 
$\sigma ^{M}e=1+\sum\Sb a\in\Bbb Z^0(p) \endSb 
t^{-ap^{M}}D_{a,M}\in\Cal A
\otimes O(\sigma ^{M}K)$, where $O(\sigma ^MK)=
\mathbin{\underset{n}\to\varprojlim }O_n(\sigma ^MK)$. 
Then 
$\sigma ^{M}f\in\Cal A\otimes O(\sigma ^{M}K(p))$ 
satisfies the relation 
$\sigma (\sigma ^{M}f)=(\sigma ^{M}f)(\sigma ^{M}e)$ 
and induces the same morphism 
$\psi :\Gamma (p)\longrightarrow\Cal A$ as $f$.  Indeed, for any 
$\tau\in\Gamma (p)$, 
$$\tau (\sigma ^{M}f)(\sigma ^{M}f)^{-1}=
\sigma ^{M}(\tau (f)f^{-1})=\sigma ^{M}(\psi (\tau ))=\psi (\tau )$$
because $\sigma $ acts trivially on $\Cal A$. 

This means that we can still study the ramification filtration 
$\{\Cal A^{(v)}\operatorname{mod}p^{M+1}\}_{v\geqslant 0}$ by working inside the lifting 
$O'_{M+1}(K_1(p))\supset O_{M+1}(\sigma ^MK(p))$ associated with our auxiliary field 
$K_1$ and its uniformiser $t_{K_1}$.

Set $\Cal B=\Cal A_{C,M+1}\operatorname{mod}\Cal J_{C,M+1}^3$ 
and for any $v\geqslant 0$, $\Cal B^{(v)}=\Cal A^{(v)}_{C,M+1}
\operatorname{mod}\Cal J_{C,M+1}$. 
We shall also use the notation $\Cal B_k=\Cal B\otimes W_{M+1}(k)$, 
$\Cal B_{K_1}=\Cal B\otimes O'_{M+1}(K_1)$, 
and 
$\Cal B_{K_1(p)}=\Cal B\otimes O'_{M+1}(K_1(p))$. 
Denote again by $\Cal J$ the augmentation 
ideal in $\Cal B$. Its extensions of scalars will be denoted 
similarly by 
$\Cal J_k, \Cal J_{K_1}$ and $\Cal J_{K_1(p)}$.

Consider an abstract continuous field isomorphism $\alpha :K\longrightarrow K_1$, which is 
the identity on the residue fields and sends $t_K$ to $t_{K_1}$. 
Consider its extension to the field isomorphism 
$\hat\alpha :K(p)\longrightarrow K_1(p)$. 
Then we have an induced isomorphism of liftings 
$\hat\alpha :O_{M+1}(K(p))\longrightarrow O'_{M+1}(K_1(p))$.  
Use it to define the morphism 
$$\id\otimes\hat\alpha :\Cal A_{C,M+1,K(p)}\longrightarrow \Cal
B_{K_1(p)}$$ and set 
$f_1:=(\id\otimes\hat\alpha )(f)\in\Cal B_{K_1(p)}$. 
Then $\sigma (f_1)=f_1e_1$, where 
$e_1=(\id\otimes\hat\alpha )(e)=1+\sum\Sb a\in\Bbb Z^0(p)\endSb 
t_1^{-a}D_{a0}$. 

If $N^*\equiv 0\operatorname{mod}N_0$, then 
$\sigma ^{M+N^*}(D_{a0})=\sigma ^{M}(D_{a0})=D_{aM}$ and 
we can relate the elements $\sigma ^{M}e=1+\sum\Sb a\in\Bbb Z^0(p) \endSb 
t^{-ap^{M}}D_{a,M}$ and 
$\sigma ^{M+N^*}e_1=1+\sum\Sb a\in\Bbb Z^0(p)\endSb 
t_1^{-ap^{M}q}D_{a,M}$ by the use of the relation between $t$ and
$t_1$ from 
lemma 4.4. So, it will be natural to compare 
the elements $\sigma ^{M}f$ and $\sigma ^{M+N^*}f_1$ 
in $\Cal B_{K_1(p)}$ by introducing 
$X\in\Cal B_{K_1(p)}$ such that 
$(\sigma ^{M}f)(1+X)=\sigma ^{M+N^*}f_1$. 
This element will be used for the  
characterisation of the ideal $\Cal B^{(v_0)}$ in proposition 4.5 
below.

Notice first, that $\Cal B^{(v_0)}$ is the minimal 2-sided ideal in 
$\Cal B$ such that 
the field of definition of 
$f\operatorname{mod}\Cal B_{K_1(p)}^{(v_0)}$ 
is invariant under the action of the group 
$\Gamma (p)^{(v_0)}$. In other words, if $I$ is a 2-sided ideal 
in $\Cal B$ and $K(f,I)$ is the field of definition of 
$f\operatorname {mod}I_{K_1(p)}$, then $I$ contains $\Cal B^{(v_0)}$ if and only if 
the largest upper ramification number $v(K(f,I)/K)$ (= the  2nd coordinate 
of the last vertex of the  graph of the Herbrand function $\varphi
_{K(f,I)/K}$) 
is less than $v_0$. 

With the above notation we have the following criterion.

\proclaim{Proposition 4.5} Suppose $r^*=v(K_1/K)<v_0$. 
Then $\Cal B^{(v_0)}$ is the 
minimal element in the set of all 2-sided ideals 
$I$ such that if $K_1(X,I)$ is the field of definition of 
$X\operatorname{mod}I_{K_1(p)}$ over $K_1$ then 
its largest upper ramification number satisfies 
$v(K_1(X,I)/K_1)<qv_0-r^*(q-1)$. 
\endproclaim 

\demo{Proof} We must prove that for any 2-sided ideal $I$ in $\Cal B$,
$$v:=v(K(f,I)/K)<v_0\ \ \Leftrightarrow\ \ v_1(X):=v(K_1(X,I)/K_1)<qv_0-r^*(q-1).$$
The following proof is similar to the proof of the corresponding statement
from  [Ab1,2]. 

Suppose $v<v_0$. The existence of the field isomorphism $\hat\alpha $  
implies that 
\linebreak 
$v(K_1(f_1,I)/K_1)=v$. Then   
$$v_1:=v(K_1(f_1,I)/K)=\max\{r^*,\varphi _{K_1/K}(v)\}\tag{$4.1$}$$
Indeed, it is sufficient to look at the maximal vertex of the Herbrand function for 
the extension 
$K_1(f_1,I)/K$ and to use the composition property 
for the corresponding Herbrand functions 
$\varphi _{K_1(f_1,I)/K}(x)=\varphi _{K_1/K}(\varphi
_{K_1(f_1,I)/K_1}(x))$.  
This implies that $v_1=r^*$ if $r^*\geqslant v$ and $v_1<v$ if 
$v> r^*$, where we have used that $\varphi _{K_1/K}(v)=r^*+(v-r^*)/q<v$ 
if $v>r^*$. 
Therefore, the largest upper ramification number of the composite 
$K(f,I)$ and $K_1(f_1,I)$ over $K$ is 
$\max \{r^*, v\}<v_0$. Clearly, $K_1(X,I)$ is contained in this
composite 
and, therefore,  $v(X):=v(K_1(X,I)/K)<v_0$. Similarly to formula (4.1) 
we obtain that  
$v(X)=\max\{r^*, \varphi _{K_1/K}(v_1(X))\}$. 
Therefore, $\varphi _{K_1/K}(v_1(X))<v_0$ and $v_1(X)<qv_0-r^*(q-1)$. 
\medskip 

Conversely, assume that $v_1(X)<qv_0-r^*(q-1)$. 
Then 
$$v(X)=\max\{r^*, \varphi _{K_1/K}(v_1(X))\}<v_0.$$ 
Suppose 
$v=v(K(f,I)/K)\geqslant v_0$. As earlier, the existence of 
$\hat\alpha $ implies that  
$v(K_1(f_1,I)/K_1)=v$ and similarly to (4.1) we have  
$$v(K_1(f_1,I)/K)=\max\{r^*,\varphi _{K_1/K}(v)\}=\varphi
_{K_1/K}(v)<v.$$
Therefore, the largest upper ramification number 
of the composite of $K_1(X,I)$ and $K_1(f_1,I)$ over $K$ equals 
$$\max \{v(K_1(X,I)/K), v(K_1(f_1,I)/K)\}=
\max\{v(X), \varphi _{K_1/K}(v)\}.$$ 
Because $K(f,I)$ is contained in 
this composite, we have 
$$v\leqslant\max\{v(X),\varphi _{K_1/K}(v)\}.$$
But $v\geqslant v_0>v(X)$ and $v>\varphi _{K_1/K}(v)$. This 
contradiction proves the proposition. . 
\enddemo

\subsubhead {\rm 4.3.4} Choosing $N^*$ and $r^*$ 
\endsubsubhead 

In order to apply the criterion from Proposition 4.5 
we shall use the special choice of $K_1=K(N^*,r^*)$, where 
$N^*\in\Bbb N$ and $r^*<v_0$ are specified 
as follows.

Introduce 
$\delta _1:=\min\{v_0-p^sa\ |\ p^sa<v_0, a\leqslant C, a\in\Bbb Z^0(p)\}$, and 
$$\delta _2:=\min\{ v_0-(p^{s_1}a_1+p^{s_2}a_2)\ |\ 
p^{s_1}a_1+p^{s_2}a_2<v_0, a_1,a_2\leqslant C, a_1,a_2\in\Bbb Z^0(p), 
s_1,s_2\in\Bbb Z\}.$$

One can see that for sufficiently large $N^*\equiv
0\operatorname{mod}N_0$,  there exists  
\linebreak 
$r^*=m^*/(q-1)<v_0$ with
$q=p^{N^*}$ and  
$ m^*\in\Bbb Z(p)$ such that 
\medskip 
a) $-(v_0-\delta _1)q+r^*(q-1)>Cp^M$; 
\medskip 
 
b) $-(v_0-\delta _2)q+r^*(q-1)>0$;
\medskip 

c) $v_0q<2r^*(q-1)$.  
\medskip 

So, we may assume that $K_1=K(N^*,r^*)$ where $N^*\equiv 0\operatorname{mod}N_0$
and the above inequalities a)-c) hold. 
\medskip 

\subsubhead {\rm 4.4} A recurrence formula for $X$ 
\endsubsubhead 

Set $\Theta ^*=t_1^{r^*(q-1)}$. Then 
$$\omega =\sigma ^{M}e-\sigma ^{M+N^*}e_1=
\sum\Sb a\in\Bbb Z^0(p)\endSb t_1^{-ap^{M}q}(E(a,\Theta
^*)^{p^M}-1)D_{aM}\in\Cal J_{K_1}.$$
The relation $1+X=(\sigma ^Mf)^{-1}(\sigma ^{M+N^*}f_1)$ implies that 
$$1+\sigma X=(\sigma ^Me)^{-1}(1+X)(\sigma ^{M+N^*}e_1)$$ and 
$$X-\sigma X=\omega +(\sigma ^Me-1)\sigma X-X(\sigma ^{M+N^*}e_1-1). \tag {$4.2$}$$ 

If $\bar X:=X\operatorname{mod}\Cal J_{K_1(p)}^2$, then the above relation (4.2) gives 
$\bar X-\sigma \bar X=\omega\operatorname{mod}\Cal J^2_{K_1(p)}$. 
We shall use this relation in n.4.5 below to study $\bar X$. Now 
(4.2) can be rewritten as 
$$X-\sigma X=\omega -\omega (\sigma ^{M+N^*}e_1-1)-[\sigma \bar X,\sigma ^{M+N^*}e_1-1]
+\omega \sigma (\bar X),  \tag{$4.3$}$$
using that $X\equiv \omega +\sigma
X\operatorname{mod}\Cal J^2_{K_1(p)}$. 
We shall use this relation in nn.4.6-4.7 below to study the field of definition of $X$. 
\medskip 

\subsubhead {\rm 4.5} The study of $\bar X$
\endsubsubhead 

For $0\leqslant r\leqslant M$ and $b\in\Bbb Z_p$, introduce $\Cal
E_r(b,T)\in \Bbb Z_p[[T]]$ 
as follows:
\medskip 

$\Cal E_0(b,T)=E(b,T)-1$, where $E(b,T)$ is the 
generalisation of the Artin-Hasse
exponential from n.1.4; 
\medskip 

$\Cal E_1(b,T)=E(b,T)^p-E(b,T^p)=(\exp (pbT)-1)E(b,T^p)$, 
\medskip 

.............................................
\medskip 

$\Cal E_M(b,T)=E(b,T)^{p^M}-E(b,T^p)^{p^{M-1}}
=(\exp (p^MbT)-1)E(b,T^p)^{p^{M-1}}.$
\medskip 

Notice the following simple properties: 
\roster

\item \ $E(b,T)^{p^M}-1=\Cal E_0(b,T^{p^M})+\Cal E_1(b,T^{p^{M-1}})+\dots
+\Cal E_M(b,T)$;
\medskip

\item \ $\Cal E_r(b,T)=p^rT+p^rT^2g_r(T)$, where $0\leqslant r\leqslant M$ and
$g_r\in \Bbb Z_p[[T]]$. 
\endroster 
\medskip 

Consider the decomposition 
$\omega =\sum\Sb r+s=M\endSb \sigma ^r\omega _s$  (cf. n.4.4 for  the
definition of $\omega $), 
where  
$$\omega _s:=\sum\Sb a\in\Bbb Z^0(p) \endSb t_1^{-ap^sq}\Cal
E_s(a,\Theta ^*)D_{as},$$
for $0\leqslant s\leqslant
M$. 
Note that 
 $p^sD_{as}\in\Cal
B_k^{(v_0)}\operatorname{mod}\Cal J_k^2$, whenever 
$p^sa\geqslant v_0$, 
cf. proposition 4.2.
Also,  if $p^sa<v_0$ then $-ap^sq+r^*(q-1)>Cp^M$, cf. n.4.3.4, 
and we have  
$t_1^{-ap^sq}\Cal E_s(a,\Theta ^*)\in t_1^{Cp^M}\m _1$, where $\m
_1:=t_1W_M(k)[[t_1]]$. 
\medskip 

So, for $0\leqslant s\leqslant M$, 
$$\omega _s\in \Cal B_{K_1}^{(v_0)}+t_1^{Cp^M}
\Cal J_{\m _1}+\Cal J_{K_1}^2,\tag {$4.4$}$$
where $\Cal J_{\m _1}=\Cal J\otimes \m _1$. 

For $0\leqslant s\leqslant M$, consider $X_s\in\Cal B_{K_1(p)}$ such
that 
$X_s-\sigma X_s=\omega _s$. Because of (4.4), 
we 
may assume that $X_s\equiv\sum\Sb u\geqslant 0\endSb \sigma ^u\omega _s
\operatorname{mod}(\Cal B^{(v_0)}_{K_1(p)}+\Cal J^2_{K_1(p)})$. 
Notice that 
$$\bar X\equiv\sum\Sb r+s=M\endSb \sigma ^r(X_s)\operatorname{mod}\Cal
J^2_{K_1(p)},$$
and 
after  
replacing 
the infinite sum $\sum_{u\geqslant 0}$ by its first $(N^*-s)$ terms  
in the above congruence for $X_s$, we obtain  
$$\bar X=\sum\Sb u+s\geqslant M\\ u<N^* \endSb 
\sigma ^u\omega _s \operatorname{mod}(\Cal B^{(v_0)}_{K_1(p)}
+\Cal J_{K_1(p)}^2+t_1^{Cp^Mq}\Cal J_{\m _1}).\tag{$4.5$}$$
\medskip 

\subsubhead {\rm 4.6.} The study of $X$
\endsubsubhead 

From the above formulas (4.4) it follows that 
$\bar X$ and $\sigma (\bar X)$ belong to 
$\Cal B_{K_1(p)}^{(v_0)}+t_1^{Cp^M}\Cal J_{\m _1}+\Cal J^2_{K_1(p)}$. 
This implies that 
$$\omega\sigma (\bar X)\in \Cal B^{(v_0)}_{K_1(p)}\Cal J_{K_1(p)}+
\Cal J_{\m _1}.$$
Therefore, when solving  equation (4.3) for $X$, this term will not  have any 
influence on the field of definition of $X\operatorname{mod}
\Cal B^{(v_0)}_{K_1(p)}\Cal J_{K_1(p)}$. 

For a similar reason, we may replace $\bar X$ in (4.3) 
by the right hand side from (4.5) 
without affecting the field of definition of 
$X\operatorname{mod}\Cal B^{(v_0)}_{K_1(p)}\Cal J_{K_1(p)}$. 
The new right hand side will be then equal to 

$$\sum\Sb a\in\Bbb Z^0(p)\\ 0\leqslant s\leqslant M\endSb 
t_1^{-ap^Mq}\Cal E_s(a,\Theta ^{*p^{M-s}})-
\sum\Sb a_1,a_2\in\Bbb Z^0(p) \\ 0\leqslant s\leqslant M\endSb 
t_1^{-(a_1+a_2)p^Mq}\Cal E_s(a_1,\Theta ^{*p^{M-s}})D_{a_1M}D_{a_2M}$$
$$\qquad\qquad\qquad -\sum\Sb 0\leqslant s_1\leqslant M,
a_1,a_2\in\Bbb Z^0(p)\\ N^*>u>M-s_1\endSb
t_1^{-a_1p^{s_1+u}q-a_2p^Mq}\Cal E_{s_1}(a_1,\Theta ^{*p^u})[D_{a_1,s_1+u},D_{a_2,M}].$$
Finally we can apply the Witt-Artin-Schreier equivalence to 
the last formula to deduce that modulo any ideal containing the ideal 
$\Cal B^{(v_0)}_{K_1(p)}\Cal J_{K_1(p)}$, the elements 
$X$ and $X'$, where 
$$X'-\sigma X'=\sum\Sb 0\leqslant s\leqslant M \endSb  
t_1^{-ap^sq}\Cal E_s(a_1, \Theta ^*)D_{as}
-\sum\Sb 0\leqslant s\leqslant M\endSb t_1^{-(a_1+a_2)p^sq}\Cal
E_s(a_1,\Theta ^*)D_{a_1s}D_{a_2s}$$

$$\qquad\qquad\qquad -\sum\Sb 0\leqslant s_1\leqslant M \\ M-N^*<s_2<s_1\endSb 
t_1^{-(a_1p^{s_1}+a_2p^{s_2})q}\Cal E_{s_1}(a,\Theta ^*)[D_{a_1s_1},D_{a_2s_2}]
\tag{4.6} $$
have the same field of definition. 

We can use this relation to find the minimal ideal $I$ in $\Cal B$ 
such that 
$X\operatorname{mod}I_{K_1(p)}$ is defined over an extension of 
$K_1$ with upper ramification number  less than 
$qv_0-r^*(q-1)$. Indeed, we
know that 
$I\operatorname{mod}\Cal J^2=\Cal B^{(v_0)}
\operatorname{mod}\Cal J^2$ and therefore, 
we may always assume
that  
$I\supset \Cal B^{(v_0)}\Cal J$. 
As before, we are also allowed to change the 
right hand side of ($4.6$) by any element of  
$\Cal B\otimes \Cal J_{\m _1}$. We may always assume that 
$I\supset \Cal B(v)$ for any $v>v_0$, because $I$ must 
contain 
all $\Cal B^{(v)}$ with $v>v_0$ and, by the  inductive assumption,  
$\Cal B^{(v)}$ coincides with $\Cal B(v)$. So, we can assume that $I$ 
contains the ideal   
$\Cal B^{(v_0+)}$ generated 
by $\Cal B^{(v_0)}\Cal J$ and all $\Cal B^{(v)}$ 
with $v>v_0$. 
\medskip 

\subsubhead {\rm 4.7.} Final simplification of {\rm (4.6)}
\endsubsubhead 

For $0\leqslant s\leqslant M$, consider the identity 
$\Cal E_s(a,\Theta ^*)=p^sat_1^{r^*(q-1)}+p^st_1^{2r^*(q-1)}g_r(t_1)$ 
from n.4.5.  

\proclaim{Lemma 4.6}
$p^st_1^{-(a_1+a_2)p^sq+2r^*(q-1)}D_{a_1s}D_{a_2s}\in\Cal
B_{K_1}^{(v_0)}\Cal J_{K_1}
+\Cal J_{\m _1}$. 
\endproclaim 

\demo{Proof} Indeed, if $p^{s}a_1\geqslant v_0$ 
(resp. if  
$p^sa_2\geqslant v_0$) then 
$p^sD_{a_1s}$ (resp. $p^sD_{a_2s}$) belongs to 
$\Cal B_k^{(v_0)}\operatorname{mod}\Cal J_k^2$.

If both $p^sa_1, p^sa_2$ are less than $v_0$ then  we use 
the fact that 
$$-(a_1+a_2)p^sq+2r^*(q-1)>Cp^M+Cp^M>0,$$ 
cf. n 4.3.4, to conclude that the corresponding term belongs to $\Cal J_{\m
_1}$. 

The lemma is proved
\enddemo  

The following lemma deals with the terms coming from the third sum 
and can be proved similarly. 

\proclaim{Lemma 4.7} 
 $p^{s_1}t_1^{-(p^{s_1}a_1+p^{s_2}a_2)q+2r^*(q-1)}[D_{a_1s_1},D_{a_2s_2}]
\in\Cal B_{K_1}^{(v_0)}\Cal J_{K_1}+\Cal J_{\m _1}$. 
\endproclaim 

The next lemma deals with the terms coming from the first sum.

\proclaim{Lemma 4.8} 
$p^st_1^{-ap^sq+2r^*(q-1)}D_{as}\in\Cal B^{(v_0+)}_{K_1}
+\Cal J_{\m _1}$. 
\endproclaim 

\demo{Proof} There is nothing to prove if $-ap^sq+2r^*(q-1)>0$. 

Assume now that $ap^sq\geqslant 2r^*(q-1)$. Consider the expression for 
$\Cal F_{ap^s}$, cf. n.4.2. Notice that $ap^s>v_0$ 
(use estimate c) from n.4.3.4) and, 
therefore,  $\Cal F_{ap^s}\in\Cal
B_k{(ap^s)}=\Cal B_k^{(ap^s)}$.

It will be sufficient to show that any term of degree 2 in the 
expression of $\Cal F_{ap^s}$ belongs to 
$\Cal B_k^{(v_0)}\Cal J_k$. Indeed, it then follows that 
the linear term 
$p^saD_{as}$ of $\Cal F_{ap^s}$ belongs to 
$\Cal B_k^{(ap^s)}+\Cal B_k^{(v_0)}\Cal J_k\subset\Cal B_k^{(v_0+)}$ 
and the statement of our lemma is proved. 

In order to prove this property of degree 2 terms notice that 
all of them  contain as
a 
factor either a product $p^{s_1}D_{a_1s_1}D_{a_2s_2}$ 
or a product $p^{s_1}D_{a_2s_2}D_{a_1s_1}$, where $s_1\geqslant s_2$ and 
$p^{s_1}a_1+p^{s_2}a_2=p^sa$. Then we have the following two cases: 
\roster

\item \ if either $p^{s_1}a_1\geqslant v_0$ or $p^{s_1}a_2\geqslant  v_0$ then 
this product belongs 
to $\Cal B^{(v_0)}_{k}\Cal J_{k}$; 
\medskip  

\item \ if both $p^{s_1}a_1$ and $p^{s_1}a_2$ are less than $v_0$, then 
$p^{s_1}a_1<v_0-\delta _1$ and 
\linebreak 
$p^{s_2}a_2\leqslant p^{s_1}a_2 <v_0-\delta _1$. 
Therefore, 
$$2r^*(q-1)\leqslant p^saq=(p^{s_1}a_1+p^{s_2}a_2)q<2q(v_0-\delta _1).$$
This contradicts the assumption a)  from n.4.3.4.  
\endroster  
\medskip 

The lemma is completely proved.
\enddemo 

By the above three lemmas, we can everywhere replace the factors 
$\Cal E_s(a,\Theta ^*)$ by $p^sat_1^{r^*(q-1)}$ and, therefore, the
right hand side 
of (4.6) is congruent modulo 
$\Cal B_{K_1}^{(v_0+)}+\Cal J_{\m _1}$ to the sum 
$\sum\Sb \gamma \geqslant 0\endSb t_1^{-q\gamma +r^*(q-1)}\Cal
F'_{\gamma }$, 
where $\Cal F'_{\gamma }$ is given by the same formula as $\Cal
F_{\gamma }$, cf. n.4.2, 
but with the additional restriction $n_2>M-N^*$ in the last sum. 

\proclaim{Lemma 4.9} If $\gamma\geqslant v_0$ then 
$\Cal F'_{\gamma }\equiv\Cal F_{\gamma }\operatorname{mod}\Cal
B_k^{(v_0)}\Cal J_k$. 
\endproclaim 

\demo{Proof} Suppose the term 
$p^{n_1}a_1[D_{a_1n_1},D_{a_2n_2}]$ enters into the formula for $\Cal
F_{\gamma }$ but does not enter into the formula for $\Cal F'_{\gamma
}$. 

Then $a_1,a_2\leqslant C$, $p^{n_1}a_1+p^{n_2}a_2=\gamma\geqslant v_0$ 
and $n_2\leqslant M-N^*$. Then 
$$p^{n_1}a_1=\gamma -p^{n_2}a_2\geqslant v_0-p^Mq^{-1}C
>r^*(1-q^{-1})-p^Mq^{-1}C>v_0-\delta _1$$
(use 4.3.2 a)). Therefore, $p^{n_1}a_1\geqslant v_0$, 
$p^{n_1}D_{a_1n_1}\in\Cal B_{k}^{(v_0)}\Cal J_k^2$ and 
$p^{n_1}a_1[D_{a_1n_1},D_{a_2n_2}]\in\Cal B_k^{(v_0)}\Cal J_k$.

The lemma is proved. 
\enddemo 

Now notice that: 
\medskip 

$\bullet $\ \ if $\gamma >v_0$, then the term 
$t_1^{-q\gamma +r^*(q-1)}\Cal F_{\gamma }$ belongs to 
$\Cal B_{K_1}(\gamma )=\Cal B_{K_1}^{(\gamma )}$; 
\medskip 

$\bullet $\ \ if $\gamma <v_0$, then the term 
$t_1^{-q\gamma +r^*(q-1)}\Cal F'_{\gamma }$ belongs to $\Cal J_{\m _1}$.

So, the ideal $\Cal B^{(v_0)}$ 
appears as the minimal ideal $I$ of $\Cal B$ such that 
$I$ contains  
the ideal 
$\Cal B^{(v_0+)}$ 
and such that 
the largest upper 
ramification number of the field of definition over $K_1$ 
of the solution 
$X''\in\Cal B_{K_1(p)}\operatorname{mod}I_{K_1(p)}$ 
of the equation 
$$X''-\sigma X''=\Cal F_{v_0}t_1^{-qv_0+r^*(q-1)}
\operatorname{mod}I_{K_1(p)}$$
is less than $qv_0-r^*(q-1)$.
\medskip  

It only remains to notice that $p\Cal F_{v_0}\in\Cal B^{(v_0+)}_k$,  and 
if $\Cal F_{v_0}\notin I_k$ then the upper ramification number of 
the field of definition 
$K_1(X'', I)$ over $K_1$ is equal to 
\newline 
$qv_0-r^*(q-1)$. 

The theorem is proved. 
\medskip 
\medskip

\subhead 5. Compatibility with ramification filtration 
\endsubhead 

In this section with the notation from n.1, 
$A=\Cal A\operatorname{mod}\Cal J^3$, $A_k=A\otimes
W(k)$. 
For any $v\geqslant 0$, 
$A^{(v)}=\Cal A^{(v)}\operatorname{mod}\Cal J^3$, 
$A^{(v)}_k:=A^{(v)}\otimes W(k)$. We also set 
$J=\Cal J\operatorname{mod}\Cal J^3$ 
with the corresponding extension of scalars $J_k=J\otimes W(k)$. 
Suppose $f$ 
is a continuous automorphism of  
the $\Bbb Z_p$-algebra A such that,  
for any $v\geqslant 0$, $f(A^{(v)})=A^{(v)}$. 
Consider the identification 
$\Cal J\operatorname{mod}\Cal J^2=\Gamma (p)^{\ab }$ 
from part b) of proposition 1.2 
and denote again by $f$ the continuous automorphism of 
$\Cal M=I (p)^{\ab }\operatorname{mod}p$ induced by $f$. 
Consider the standard topological generators 
$D_{an}$, $a\in\Bbb Z(p), n\in\Bbb Z\operatorname{mod}N_0$, 
for $\Cal M$ and set, for any $a\in\Bbb Z(p)$, 
$$f(D_{a0})=\sum\Sb b,m \endSb \alpha _{abm}(f)D_{bm},$$
where the coefficients $\alpha _{abm}(f)\in k$. With the above
notation, the principal results of this section are:
\medskip 

if $\alpha _{110}(f)\ne 0$ and $N_0\geqslant 3$ then  
\medskip 

$\bullet $ \ \ there is an $\eta\in\Aut ^0K$ such that 
for any $a,b\in\Bbb Z(p)$ and $a\leqslant b<
p^{N_0-3}$, it holds 
$\alpha _{ab0}(f)=\alpha _{ab0}(\eta ^*)$;
\medskip 

$\bullet $\ \ if $a\leqslant b<p^{N_0-3}$ and $m\in\Bbb N$ is such
that $b/p^m<a$ then $\alpha _{a,b,-m\operatorname{mod}N_0}(f)=0$. 
\medskip 

\subsubhead {\rm 5.1.} The elements $\Cal F_{\gamma }{(v)}$ 
\endsubsubhead 

 By Theorem B, cf. n.4.2,  for any $v\geqslant 0$, the ideal  $A_k^{(v)}$ is 
the minimal closed $\sigma $-invariant ideal in $A_k$ containing 
the explicitly given elements $\Cal F_\gamma $, 
for all $\gamma\geqslant v$. 
For any $a\in\Bbb Z(p)$ and 
$n\in\Bbb Z\operatorname{mod}N_0$, set $\Delta _{a0}=(1/a)\Cal F_a$ 
and $\Delta _{an}=\sigma ^n\Delta _{a0}$. Then 
$\Delta _{an}\equiv D_{an}\operatorname{mod}\Cal J_k^2$ and 
$\{\Delta _{an}\ |\ a\in\Bbb Z(p), n\in\Bbb Z\operatorname{mod}N_0\}
\cup\{D_0\}$ 
 is a new system of topological generators for $A_k$. 
The elements of this  new set of generators 
together with their pairwise products form a topological basis 
of the $W(k)$-module $A_k$.  
\medskip 

For any $\gamma \geqslant v\geqslant 0$, consider the following 
elements $\Cal F_{\gamma }{(v)}$ (these elements have already been 
mentioned in n.4.2):
\medskip  

If $\gamma =ap^m$ with $a\in\Bbb Z(p)$ and $m\in\Bbb Z_{\geqslant
 0}$ set 
$$\Cal F_{\gamma }{(v)}=p^ma\Delta _{am}-
\sum\Sb n\geqslant 0, a_1,a_2\in\Bbb Z(p) \\ 
p^{n}(a_1+a_2)=\gamma \\
p^na_1,p^na_2<v \endSb 
p^{n}a_1\Delta _{a_1n}\Delta _{a_2n};$$

\medskip 

If $\gamma \notin\Bbb Z$ set 
$$\Cal F_{\gamma }{(v)}=-\sum\Sb n_1\geqslant 0, a_1,a_2\in\Bbb Z(p) \\ 
p^{n_1}a_1+p^{n_2}a_2=\gamma \\ 
p^{n_1}a_1,p^{n_1}a_2<v \endSb 
p^{n_1}a_1[\Delta _{a_1n_1},\Delta _{a_2n_2}].$$

Similarly to n.4.2, we have  the following property.  

\proclaim{Proposition 5.1} For any  $v\geqslant 0$, 
$A_k^{(v)}$ is the minimal $\sigma $-invariant 
closed ideal of $A_k$ containing the elements $\Cal F_{\gamma }{(v)}$ for all 
$\gamma\geqslant v$. 
\endproclaim

\subsubhead {\rm 5.2.} The submodules $A_{\tr }^{(v)}$ and $A_{\adm }^{(v)}$ 
\endsubsubhead 

Suppose $v\geqslant 0$. 

Let $A_{\tr}^{(v)}$ be the $W(k)$-submodule in $A_k$ generated by 
the following elements:
\medskip 

$\tr _1$) $p^s\Delta _{an}$ with $s\geqslant 0$ and $p^sa\geqslant
2v$; 
\medskip 

$\tr _2$) $p^s\Delta _{a_1n_1}\Delta _{a_2n_2}$ with $a_1,a_2\in\Bbb
Z(p)$, $s\geqslant 0$ and $n_1,n_2\in\Bbb Z\operatorname{mod}N_0$ 
such that 
\linebreak 
$\max\{p^sa_1,p^sa_2\}\geqslant v$. 
\medskip 

Let $A^{(v)}_{\adm }$ be the minimal closed 
$W(k)$-submodule in $A_k$ containing $A_{\tr }^{(v)}$ and the
following elements: 
\medskip 

$\adm _1$) $p^s\Delta _{an}$, with $s\geqslant 0$, $a\in\Bbb Z(p)$ and
$p^sa\geqslant v$;
\medskip 

$\adm _2$) $p^s\Delta _{a_1n_1}\Delta _{a_2n_2}$, where
$a_1,a_2\in\Bbb Z(p)$, $n_1,n_2\in\Bbb Z\operatorname{mod}N_0$ and 
$s=s(a_1,a_2)\in\Bbb Z_{\geqslant 0}$ are such that:   
\roster 
\item \ \ $v/p\leqslant\max\{p^sa_1,p^sa_2\}<v$;
\item 
\ \  
$\dsize \max\left\{ p^s\left (a_1+\frac{a_2}{p^{n_{12}}}\right ), 
p^s\left (\frac{a_1}{p^{n_{21}}}+a_2\right )\right \}\geqslant v,$
where $0\leqslant n_{12},n_{21}<N_0$, $n_{12}\equiv
n_1-n_2\operatorname{mod}N_0$ and 
$n_{21}\equiv n_2-n_1\operatorname{mod}N_0$; 
\item 
\ \ if $n_1=n_2$ then $a_1+a_2\equiv 0\operatorname{mod}p$.
\endroster 

\proclaim{Proposition 5.2} For any $v\geqslant 0$, 
\newline 
{\rm 1)} $f(A^{(v)}_{\tr })=A_{\tr }^{(v)}$; 
\newline 
{\rm 2)} 
$A^{(v)}_{\adm }\supset A_k^{(v)}\supset A_{\tr }^{(v)}\supset
pA_{\adm }^{(v)}$; 
\newline 
{\rm 3)} the elements from $\adm _1)$ and $\adm _2)$ form a $k$-basis
of  
$A_{\adm }^{(v)}\operatorname{mod}A_{\tr }^{(v)}$. 
\endproclaim 

\demo{Proof} 1) It is sufficient to notice that $A^{(v)}_{\tr }$ is
the minimal $\sigma $-invariant $W(k)$-submodule in $A$ containing 
$\sum\Sb \gamma\geqslant 2v \endSb \Cal F_{\gamma }(v)W(k)
+\sum\Sb \gamma\geqslant v \endSb \Cal F_{\gamma }(v)J_k$. 

2) From the above n.1) it follows that $A_k^{(v)}\supset A_{\tr
   }^{(v)}$. 
The embedding $A_k^{(v)}\subset A_{\adm }^{(v)}$ follows from the
   definition 
of $A^{(v)}_{\adm }$: as a matter of fact, 
$A_{\tr }^{(v)}$ is spanned by all summands of elements 
$\sigma ^s\Cal F_{\gamma }$ with $s\in\Bbb Z\operatorname{mod}N_0$ and
   $\gamma\geqslant v$. The embedding $pA_{\adm }^{(v)}\subset A_{\tr
   }^{(v)}$ follows from the fact that each element listed in $\adm
   _1)$ and $\adm _2)$ belongs to $A_{\tr }^{(v)}$ after
   multiplication  by $p$. 

3) It is easy to see that any $k$-linear combination of the elements from 
$\adm _1)$ and $\adm _2)$ does not belong to $A^{(v)}_{\tr
}\operatorname{mod}pA_{\adm }^{(v)}$. 
\enddemo

\proclaim{Proposition 5.3} Suppose $v\geqslant 0$ and 
$p^s\Delta _{a_1n_1}\Delta _{a_2n_2}$ is one of elements 
listed in $\adm _2)$. Let $n=\min\{n_{12},n_{21}\}$. 
If 
$$v/p^{N_0-n}\leqslant d(v):=\min\{v-a\ |\ a\in\Bbb Z,a<v\}$$
then there are unique $m\in \Bbb Z\operatorname{mod}N_0$ and 
$\gamma\geqslant v$ such that $p^sa_1\Delta _{a_1n_1}\Delta _{a_2n_2}$ 
appears (with non-zero coefficient) 
in the expression of $\sigma ^m\Cal F_{\gamma }(v)$.
\endproclaim 

\remark{Remark} We are going to apply this proposition in the
following situations:
\roster 

\item \ $v\in\Bbb N$, $v<p^{N_0}$, $n_1=n_2=0$;
\medskip 

\item \ $v=c+1/p$, $n_1=0$, $n_2=-1$, where 
$c\in\Bbb N$ and $c<p^{N_0-2}$.
\endroster 
\endremark
\medskip 

\demo{Proof} By symmetry we may assume that $n=n_{12}$. 

If $n_{12}\ne 0$ we have $\dsize p^s\left (a_1+\frac{a_2}{p^n}\right
)=\gamma\geqslant v$, because of property $\adm _2)(2)$, and 
$$p^s\left (\frac{a_1}{p^{N_0-n}}+a_2\right
)<\frac{v}{p^{N_0-n}}+p^sa_2
\leqslant d(v)+(v-d(v))=v\leqslant\gamma .$$
Therefore, the term $p^s\Delta _{a_1n_1}\Delta _{a_2n_2}$ appears in 
the expression of 
$\sigma ^{n_1-s}\Cal F_{\gamma }(v)$. This term will appear 
in the expression of another 
$\sigma ^{n'}\Cal F_{\gamma '}(v)$, where $\gamma '\geqslant v$, 
if and only if 
$\dsize p^s\left (a_1+\frac{a_2}{p^{n+mN_0}}\right )\geqslant v$ or 
$\dsize p^s\left (\frac{a_1}{p^{mN_0-n}}+a_2\right )\geqslant v$, 
where $m\in\Bbb N$. But the condition 
$v/p^{N_0-n}<d(v)$ implies that all such numbers are 
less than $v$.

If $n_{12}=0$ then $\gamma =p^s(a_1+a_2)\geqslant v$ and 
$p^{s}\Delta _{a_1n_1}\Delta _{a_2n_2}$ appears in the expression of 
$\sigma ^{n_1-s}\Cal F_{\gamma }(v)$. This element can appear in 
the expression of another 
$\sigma ^{n'}\Cal F_{\gamma '}(v)$, where $\gamma '\geqslant v$, if and only if 
$\dsize\gamma '=p^s\left (a_1+\frac{a_2}{p^{mN_0}}\right )\geqslant v$ 
or 
$\dsize\gamma '=p^s\left (\frac{a_1}{p^{mN_0}}+a_2\right )\geqslant v$, 
where $m\in\Bbb N$. As earlier, 
$\gamma '<v$ in both cases. 

The proposition is proved. 
\enddemo  

\remark{Remark} 
If $v/p^{N_0/2}<d(v)$, then 
elements of the set 
$$\{ \sigma ^s\Cal F_{\gamma }^{(v)}
\operatorname{mod}A_{\adm}^{(v)}\ |\ 0\leqslant s<N_0, 
\gamma\geqslant v\}$$
are linear combinations of disjoint groups of elements listed 
in $\adm _1)$ and $\adm _2)$. 
\endremark 
\medskip

5.3. Denote by the same 
symbol $f$ the  morphism 
of $W(k)$-modules 
$$A^{(v)}\operatorname{mod}A^{(v)}_{\tr }\longrightarrow
A^{(v)}\operatorname{mod}A^{(v)}_{\tr },$$
which is induced by $f:A\longrightarrow A$. As earlier, 
denote again by $f$   
the $k$-linear extension of 
the automorphism of $\Cal M$, which is induced by $f$. Because 
the images of $D_{an}$ and $\Delta _{an}$ coincide in $\Cal M_k$, 
we have, for any $a\in\Bbb Z(p)$, 
$$f(\Delta _{a0})=\sum\Sb b\in\Bbb Z(p) \\ m\in\Bbb Z\operatorname{mod}N_0 
\endSb \alpha _{abm}(f)\Delta _{bm}.$$
It will be convenient sometimes to set 
$\alpha _{ab0}(f)=0$ if $a$ or $b$ are divisible by $p$.

\proclaim{Proposition 5.4} 
Suppose $\alpha _{110}(f)=\alpha \in k^*$. 
Then  $\alpha _{aa0}(f)=\alpha ^{a}$, 
for any $a\in\Bbb Z(p)$ such that 
$a<p^{N_0-1}$ if $p\ne 2$ and $N_0\geqslant 2$, 
and such that $a<2^{N_0}$ if $p=2$ and $N_0\geqslant 3$. 
\endproclaim 

\demo{Proof} By proposition 5.3, for any $v\leqslant p^{N_0}$ such that $v\equiv 0\operatorname{mod}p$, 
$f(\Cal F_v(v))\operatorname{mod}A_{\tr }^{(v)}$ must contain all terms 
$a_1\Delta _{a_10}\Delta _{a_20}$, for which $a_1+a_2=v$, and the term 
$p^sa\Delta _{as}$, where $p^sa=v$ and $a\in\Bbb Z(p)$, with the same coefficient. In 
other words, for such indices $a_1,a_2,a\in\Bbb Z(p)$, 
$$\alpha _{a_1a_10}(f)\alpha _{a_2a_20}(f)=\sigma ^s\alpha
_{aa0}(f).\tag {5.1}$$
For $a\in\Bbb Z(p)$, $a<p^{N_0}$, 
set $\gamma (a)=\alpha _{aa0}(f)\alpha
_{110}(f)^{-1}$. Then 
$\gamma (1)=1$ and $\gamma (a_1)\gamma (a_2)=\gamma (a)^{p^s}$ 
if $a_1+a_2=p^sa$. 

Suppose $p\ne 2$. 

First, we prove that for $n\in\Bbb Z(p)$ satisfying 
$1\leqslant n<p^{N_0-1}$, we
have 
$$\gamma (n)=\gamma (2)^{n-1}.\tag{5.2}$$

This is obviously true for $n=1$ and $n=2$. 

Assume that $n\geqslant 2$ and that  
$\gamma (m)=\gamma (2)^{m-1}$ holds for all $m\in\Bbb Z(p)$ such that 
$m\leqslant n$. 
Consider a special case of relation (5.1) with $n\in\Bbb Z(p)$ 
$$\gamma (1)\gamma (np-1)=\gamma (n)^p\tag {5.3}$$
If $n\not\equiv -1\operatorname{mod}p$ then use the relation 
$\gamma (p-1)\gamma (p+1)=\gamma (2)^p$, which is again 
a special case of (5.1),   
to deduce from (5.3) that 
$$\gamma (n+1)=\gamma (1)\gamma (n+1)=\gamma (n)\gamma (2)=\gamma
(2)^n.$$ 

If $n\equiv -1\operatorname{mod}p$ and $p\ne 3$ then $n\geqslant 4$ 
and by the inductive assumption $\gamma (3)=\gamma (2)^2$.  
Apply the relation 
$\gamma (p-1)\gamma (2p+1)=\gamma (3)^p=\gamma (2)^{2p}$ 
to deduce from (5.3) that 
$$\gamma (n+1)=\gamma (1)\gamma (n+2)=\gamma (n)\gamma (2)^2=\gamma
(2)^{n+1}.$$ 

If $p=3$ then $\gamma (p-1)\gamma (2p+1)=\gamma (1)^{p^2}$ and we
obtain from (5.3) that 
$$\gamma (n+1)=\gamma (1)\gamma (n+2)=\gamma (n)=\gamma
(2)^{n-1}=\gamma (2)^{n+1},$$
because $\gamma (2)=1$ (using that $\gamma (1)\gamma (2)=\gamma
(1)^3$). 

So, relation (5.2) is proved. 

Still assuming that $p\ne 2$ prove that $\gamma
(2)=1$. 
The relation $\gamma (1)\gamma (p-1)=\gamma (1)^p$ implies that 
$\gamma (2)^{p-2}=\gamma (p-1)=1$. 
The equality 
\linebreak 
$\gamma (1)\gamma (p^2-1)=\gamma (1)^{p^2}$ implies that 
$\gamma (2)^{p^2-2}=\gamma (p^2-1)=1$. Then $\gamma (2)=1$ because 
$p^2-2$ and $p-2$ are coprime. This completes the case 
$p\ne 2$.  

Consider now the case $p=2$. 

Notice that for any $n\in\Bbb Z(2)$ such that $1<n<2^{N_0}$, we have 
$n+1=2^sa$, where $a\in\Bbb Z(2)$, 
$s\in\Bbb N$ and $a<n$. Therefore, 
$\gamma (1)\gamma (n)=\gamma (a)^{2^s}$ 
and the equality $\gamma (n)=1$ follows by induction on $n\geqslant
1$ for all $n<2^{N_0}$. 
\enddemo

\proclaim{Corollary 5.5} If $\alpha _{110}(f)=1$ then 
$\alpha _{aa0}(f)=1$ whenever 
$a<p^{N_0-1}$, $p\ne 2$ or $a<2^{N_0}$, $p=2$. 
\endproclaim

\proclaim{Proposition 5.6} Suppose $N_0\geqslant 3$, 
$\alpha _{110}(f)\in k^*$, 
$a,b\in\Bbb Z(p)$, $a,b<p^{N_0-2}$. 
If 
$0\leqslant m<N_0$ and $b/p^m<a$ then   
$\alpha _{a,b,-m\operatorname{mod}N_0}(f)=0$.
\endproclaim 

\demo{Proof} 
For a given $b\in\Bbb Z(p)$, $b<p^{N_0-2}$ and $1\leqslant m<N_0$, 
let $a\in\Bbb Z(p)$ be the minimal integer such that $\alpha _{a',b,-m}(f)=0$ 
if $a'>a$. If such an $a$ does not exist then $\alpha _{a,b,-m}(f)=0$
for all $a$ and there is nothing to prove. 

If $p\ne 2$ put $v=p^{N_0-1}$ and consider $f(\Cal
F_v(v))\operatorname{mod}(A_{\tr }^{(v)}+pA_{\adm }^{(v)})$. 

We prove that the term $\Delta _{v-a,0}\Delta _{b,-m}$ enters in 
$f(\Cal F_v(v))$ with the coefficient 
$$(v-a)\alpha _{v-a,v-a,0}(f)\alpha _{a,b,-m}(f)=
-a\alpha _{v-a,v-a,0}(f)\alpha _{a,b,-m}(f).\tag {5.4}$$
Indeed, $\Cal F_v(v)\operatorname{mod}(A_{\tr }^{(v)}+pA^{(v)}_{\adm })$ 
is a sum of the terms of the form 
$a_1\Delta _{a_10}\Delta _{a_20}$ with 
$a_1,a_2\in\Bbb Z(p)$ such that 
$a_1+a_2=v$. 
Therefore, $f(a_1\Delta _{a_10}\Delta _{a_20})$ 
contains $\Delta _{v-a,0}\Delta _{b,-m}$ with coefficient 
$$a_1\alpha _{a_1,v-a,0}(f)\alpha _{a_2,b,-m}(f).$$
Now notice that $\alpha _{a_2,b,-m}(f)=0$ if $a_2>a$, 
and $\alpha _{a_1,v-a,0}(f)=0$ if $a_1>v-a$ or, equivalently, if 
$a_2<a$. So, $a_1=v-a$ and the coefficient is given by 
formula (5.4).  

By the choice of $a$, the coefficient (5.4) is not zero. 
Therefore, $\Delta _{v-a,0}\Delta _{b,-m}\in A_{\adm }^{(v)}$. 
Notice that  
$$\max\left\{v-a+\frac{b}{p^m}, \frac{v-a}{p^{N_0-m}}+b\right\}
=v-a+\frac{b}{p^m}$$
and $b/p^m\geqslant a$. Indeed, we can use that 
$$\frac{v-a}{p^{N_0-m}}+b<\frac{p^{N_0-1}}{p}+p^{N_0-2}<
2p^{N_0-2}<p^{N_0-1}-p^{N_0-2}
<v-a+\frac{b}{p^m}.$$
Therefore, $v-a+b/p^m\geqslant v$, i.e. $b/p^m\geqslant a$ and the
proposition is proved in the case $p\ne 2$. 

If $p=2$ we can take $v=2^{N_0}$ and repeat the above arguments 
by using in the last step the inequality 
$$\frac{v-a}{2^{N_0-m}}+b<\frac{2^{N_0}}{2}+2^{N_0-2}<2^{N_0}-a\left
(1-\frac{1}{2^m}\right )\leqslant 
v-a+\frac{b}{2^m}.$$
The proposition is completely proved. 
\enddemo 
\medskip 

5.4. Suppose $r\in\Bbb N$ is such that 
$\alpha _{aa'0}(f)=0$ for any 
$a,a'\in\Bbb Z(p)$ such that  
$a<a'<a+r<p^{N_0-2}$.

Let $\delta (p)$ be $p$ if $p\ne 2$ and $\delta (p)=4$ if $p=2$.

\proclaim{Proposition 5.7} Assume that 
$\alpha _{110}(f)=1$. If $b,b_1\in\Bbb Z(p)$, $b_1=b+r$ and 
$b_1+\delta (p)<p^{N_0-2}$ 
then  
$\alpha _{bb_10}(f)=\alpha _{b-\delta (p),b_1-\delta (p),0}(f)$.   
\endproclaim 

\demo{Proof}  

Let $a_0=p^{N_0-2}-1$, $\dsize v_0=a_0+1/p$, $\dsize
v=a_0+\frac{b}{p}$. We need the following lemma. 

\proclaim{Lemma } If $a',b',c\leqslant a_0$ and 
$a'+b'/p=v$ then $\alpha _{a',c,-1}(f)=0$.
\endproclaim 

\demo{Proof of lemma } It follows from the inequalities 
$$\frac{c}{p}\leqslant \frac{a_0}{p}\leqslant a_0-\frac{a_0}{p}
<v-\frac{b'}{p}=a'$$
and proposition 5.6. 
\enddemo 

We continue the proof of proposition 5.7. Consider 
$$\Cal F_v(v_0)=-\sum\Sb a'+b'/p=v \\ a',b'\leqslant a_0 \endSb 
a'[\Delta _{a'0},\Delta _{b',-1}]
\operatorname{mod}pA_{\adm }^{(v)}.$$

Using  that $v_0/p^{N_0-1}< d(v_0)=1/p$, cf. proposition 5.3,  
we can find now the coefficient for $[\Delta _{a_00},\Delta _{b_1,-1}]$ 
in $f(\Cal F_v(v_0))$. By the above lemma 
\linebreak 
$\alpha _{a',b,-1}(f)=0$, therefore the 
image of the term 
$a'[\Delta _{a'0},\Delta _{b',-1}]$ gives a coefficient 
$$a'\alpha _{a'a_00}(f)\sigma ^{-1}(\alpha _{b'b_10}(f)).$$
If $a'<a_0$ and $\alpha _{a'a_00}(f)\ne 0$ then 
$a'\leqslant a_0-r$, $b'\geqslant b+rp>b_1$ 
and $\alpha _{b'b_10}(f)=0$. So, the coefficient is non-zero 
only for $a'=a_0$. Then by Corollary 5.5 
$\alpha _{a'a_00}(f)=1$ and the coefficient will be equal to 
$a_0\sigma ^{-1}(\alpha _{bb_10}(f))$.

If $p\ne 2$ we can proceed similarly to 
find the coefficient for $[\Delta _{a_0-1,0},\Delta _{b_1+p,-1}]$ 
in $f(\Cal F_v(v_0))$. It equals $(a_0-1)\sigma ^{-1}(\alpha
_{b+p,b_1+p,0}(f))$. Therefore, by proposition 5.3 
$$\alpha _{bb_10}(f)=\alpha _{b+p,b_1+p,0}(f)$$
and the case $p\ne 2$ is completely considered. 

If $p=2$, we similarly find similarly the coefficient for 
$[\Delta _{a_0-2,0},\Delta _{b_1+4,-1}]$ 
in $f(\Cal F_v(v_0))$. It equals 
$(a_0-2)\sigma ^{-1}(\alpha _{b+4,b_1+4,0}(f))$ and we obtain 
$$\alpha _{bb_10}(f)=\alpha _{b+4,b_1+4,0}(f).$$
The proposition is proved.
\enddemo  
\medskip 

5.5. Now we come to the central point of this section.

\proclaim{Proposition 5.8} Suppose $\alpha _{110}(f)\ne 0$ and
$N_0\geqslant 3$. Then there is an 
$\eta\in\Aut ^0K$ such that 
$\alpha _{ab0}(f\eta ^*)=\delta _{ab}$, 
for any $a,b\in\Bbb Z(p)$ 
with $a\leqslant b<p^{N_0-3}$, 
 where $\delta _{ab}$ is the Kronecker symbol. 
\endproclaim 

\demo{Proof} Proposition 5.4 together with part 2) of 
proposition 2.1 imply that after replacing $f$ by
$f\eta ^*$ for some $\eta\in\Aut ^0K$ such that 
$\eta (t)=\alpha _{110}(f)t$, we can assume that 
$\alpha _{aa0}(f)=1$ if $a<p^{N_0-1}$. 

Let $r=r(f)\in\Bbb N$ be the maximal subject to the condition that 
$\alpha _{ab0}(f)=0$, 
for any $a,b\in\Bbb Z(p)$ with
$a,b<p^{N_0-2}$ and $a<b<a+r$.  
. 

If $r\geqslant p^{N_0-3}-1$ then there is nothing to prove. 
Therefore, we can
assume that $r\leqslant p^{N_0-3}-2$. For $1\leqslant a<p^{N_0-2}$, set 
$\alpha _a(r)=\alpha _{a,a+r,0}(f)$ if $a\in\Bbb Z(p)$ and 
$\alpha _a(r)=0$, otherwise. 

By proposition 5.7 $\alpha _a(r)$ depends only on the residue 
$a\operatorname{mod}\delta (p)$ and by the choice of $r$ the function 
$a\mapsto\alpha _a(r)$ is not identically zero. The proposition 
will be proved if we show the existence of 
$\eta\in\Aut ^0K$ such that $r(f\eta ^*)>r(f)$. 

In the case $p\ne 2$ apply proposition 2.5 with $w_0=1+r$. 
Let $\eta $ will be the corrsponding character. If 
$r(f\eta ^*)>r(f)$, then the proposition is proved. 
So, assume that $r(f\eta ^*)=r(f)$. Therefore, by replacing $f$ by $f\eta ^*$ we
can assume the following normalisation conditions:
\medskip 
a) \ $\alpha _1(r)=0$ if $r\not\equiv -1\operatorname{mod}p$;
\medskip 
b)\  $\alpha _2(r)=0$ if $r\equiv -1\operatorname{mod}p$.
\medskip 

In the case $p=2$, apply proposition 2.6 with either $w_0=r+2$ if $r\equiv
2\operatorname{mod}4$ or $w_0=r$ if $r\equiv 0\operatorname{mod}4$. 
In the first case we have the normalisation condition 
\medskip 
c)\ $\alpha _1(r)=\alpha _3(r)=0$;
\medskip 
in the second case we obtain only that 
\medskip 
d)\  $\alpha _1(r)=0$.
\medskip 

{\it The case $p\ne 2$.}

If $r=p^{N_0-3}-2$ then 
$\alpha _1(r)=\alpha _{ab0}(f)=0$ if 
$a=1,b=p^{N_0-3}-1$. For all other couples 
$a,b\in\Bbb Z(p)$ such that 
$a<b<p^{N_0-3}$, we have 
$\alpha _{ab0}(f)=0$ because $b-a<r$. Therefore, 
we can assume that $r\leqslant p^{N_0-3}-3$.

Let $c_j=p(r+1)+j$ for $j=1,2,\dots ,p-1$. 
Then  $c_j\leqslant p(p^{N_0-3}-2)+p-1<p^{N_0-2}$, for all $j$.  
Set $v_j=c_j+1/p$ and consider the coefficient 
for $\Cal F_{v_j+r}(v_j)$ in the image $f(\Cal F_{v_j}(v_j))\in  
A_{\adm }^{(v_j)}\operatorname{mod}A_{\tr }^{(v_j)}+pA_{\adm }^{(v_j)}$.

Similarly to the proof of proposition 5.7, we see that the term 
$[\Delta _{c_j0},\Delta _{1+rp,-1}]$ from the expression of  
$\Cal F_{v_j+r}(v_j)$ can appear with non-zero coefficient 
only as image of one of the following two terms 
from $\Cal F_{v_j}(v_j)$: 
$(c_j-r)[\Delta _{c_j-r,0},\Delta _{1+rp,-1}]$ 
and $c_j[\Delta _{c_j0},\Delta _{1,-1}]$. This coefficient is 
equal to 
$$(c_j-r)\alpha _{c_j-r}(r)+c_j\alpha _{1,1+rp,0}(f).$$  

Similarly, the term 
$[\Delta _{c_j-1,0},\Delta _{1+(r+1)p,-1}]$ from the expression of  
$\Cal F_{v_j+r}(v_j)$ can appear with non-zero coefficient 
only as image of   
$(c_j-1-r)[\Delta _{c_j-1-r,0},\Delta _{1+(r+1)p,-1}]$ 
and $(c_j-1)[\Delta _{c_j-1,0},\Delta _{1+p,-1}]$. This coefficient is 
$$(c_j-1-r)\alpha _{c_j-1-r}(r)+(c_j-1)\sigma ^{-1}\alpha
_{1+p,1+(r+1)p,0}(f).$$

Therefore, we have the following relation 
$$\frac{c_j-r}{c_j}\alpha _{c_j-r}(r)=
\frac{c_j-1-r}{c_j-1}\alpha _{c_j-1-r}(r)+X, \tag{5.5}$$
where $X=\sigma ^{-1}(\alpha _{1+p,1+(r+1)p,0}(f))-
\sigma ^{-1}(\alpha _{1,1+rp,0}(f))$.  

For $j=1,\dots ,p-1$, set $\beta _j=\dsize\frac{c_j-r}{c_j}\alpha
_{j-r}(r)$. Then the above relation (5.5) implies that 
$\beta _2=\beta _1+X, \beta _3=\beta _2+X, \dots ,\beta _{p-1}=\beta
_{p-2}+X$. 
\medskip 

{\it The case  
$r\not\equiv 0\operatorname{mod}p$, $p\ne 2$.} 

In this case the normalisation conditions imply that 
\medskip 

---\ if $r\not\equiv -1\operatorname{mod}p$ then 
$\beta _{r+1}=0$;
\medskip 

---\ if $r\equiv -1\operatorname{mod}p$ then $\beta
_{r+2}=0$.

In both cases  $\beta _r=0$. This implies that 
$\beta _1=\dots =\beta _{p-1}=0$. Therefore, 
$\alpha _a(r)=0$, for all $a$. This is a contradiction. 

So, in the case $r\not\equiv 0\operatorname{mod}p$, $p\ne 2$ the  
proposition is proved. 
\medskip

{\it The case $r\equiv 0\operatorname{mod}p$, $p\ne 2$}

In this case  we only have the normalisation condition  
$\beta _1=0$. Therefore, for $i=1,\dots ,p-1$, we have 
$\beta _i=(i-1)X$ and 
$\alpha _a(r)=(a-1)X$ for any $a\in\Bbb Z(p)$, $a<p^{N_0-3}$.

Let $v=(p-1)r+p$ and consider the coefficient for $\Cal
F_{v+r}(v)$ in the image $f(\Cal F_v(v))$. Following the images 
of terms of
degree 2 we see that this coefficient equals $-2X$. Now notice that the linear
terms in $\Cal F_v(v)$ (resp. $\Cal F_{v+r}(v)$) have coefficients 
with $p$-adic valuation 
$v_p((p-1)r+p)$ (resp. $v_p(pr+p)$). Clearly, 
if $1=v_p(pr+p)$ and if $1<v_p((p-1)r+p)$ then the linear term of 
$\Cal F_{v+r}(v)$ cannot appear in the image 
$f(\Cal F_v(v))$. Therefore, 
$1=v_p(pr+p)=v_p((p-1)r+p)$ and the linear terms 
in $\Cal F_v(v)$ (resp. $\Cal F_{v+r}(v)$) are multiples of 
$\Delta _{r+1-r/p,1}$ (resp.  
$\Delta _{r+1,1}$). But then 
$(r+1)-(r+1-r/p)=r/p<r$ and by the definition of $r$, 
$\Delta _{r+1,1}$ will not appear in the image 
$F(\Delta _{r+1-r/p,1})$. This contradiction proves the  
proposition in the case $r\equiv 0\operatorname{mod}p$, 
$p\ne 2$. 
\medskip 

{\it The case $p=2$.} 

Here $r\equiv 0\operatorname{mod}2$. 
If $r\equiv 2\operatorname{mod}4$ then the normalisation 
conditions imply that $\alpha _a(r)=0$ for all $a$ and 
the proposition is proved. 

If $r\equiv 0\operatorname{mod}4$ then we only have 
one normalisation condition 
$\alpha _a(r)=0$ if $a\equiv 1\operatorname{mod}4$. 
Let $\alpha _a(r)=\alpha $ where  
$a\equiv 3\operatorname{mod}4$. Consider 
$$\Cal F_{r+4}(r+4)=(r+4)\Delta _{\frac{r+4}{2^s},s}
+\sum\Sb a+b=r+4 \\ a,b<r+4\endSb 
\Delta _{a0}\Delta _{b0}
\in A^{(r+4)}_{\adm }\operatorname{mod}A_{\tr }^{(r+4)},$$
where $s=v_2(r+4)\geqslant 2$. 
Then $f(\Cal F_{r+4}(r+4))$ contains 
$\Delta _{r+1,0}\Delta _{r+3,0}$ with coefficient 
$$\alpha _{1,r+1,0}(f)+\alpha _{3,3+r,0}(f)=\alpha ,$$ 
and therefore it contains 
$\Cal F_{2r+4}(r+4)$ with coefficient $\alpha $. 
Similarly to the case $p\ne 2$, 
we obtain the equality 
$v_2(r+4)=v_2(2r+4)=2$ and consequently the fact 
that $f(\Delta _{r/2+1,2})$ cannot contain 
$\Delta _{r/4+1,2}$ with non-zero coefficient because 
$(r/2+1)-(r/4+1)=r/4<r$. The proposition 
is completely proved. 
\enddemo 
\medskip 

\subhead 6. Proof of the main theorem --- the characteristic $p$ case 
\endsubhead

Suppose $\char E=p$. 

Then $\char E'=p$ because the topological groups 
$\Gamma _E(p)^{\ab }$ and $\Gamma _{E'}(p)^{\ab }$ are isomorphic. 
Looking at the ramification filtrations of these groups we deduce that 
the residue fields of $E$ and $E'$ are isomorphic. Therefore, 
$E$ and $E'$ are isomorphic 
complete discrete valuation fields and  
we can identify the maximal $p$-extensions $E(p)$
of $E$ and $E'(p)$ of $E'$. 

Let $K$ be a finite Galois extension of $E$ in $E(p)$. Then $E(p)$ is
a maximal $p$-extension of $K$ and 
$\Gamma _K(p)=\Gal (E(p)/K)$. Let $K'$ be the extension 
of $E'$ in $E(p)$ such that 
$g(\Gamma _K(p))=\Gamma _{K'}(p)$ (recall that $g$ is a group isomorphism). 
If $s\geqslant 0$ and $K_s$ is the unramified extension of $K$ in $E(p)$
such that $[K_s:K]=p^s$ then 
$g(\Gamma _{K_s}(p))=\Gamma _{K'_s}(p)$, where 
$K'_s$ is the unramified extension of $K'$ in $E(p)$ of degree $p^s$. 
Therefore, with the notation from n.3 we have a compatible system 
$g_{KK'}=\{g_{KK's}\}_{s\geqslant 0}$ of $\Bbb F_p$-linear continuous 
automorphisms 
$g_{KK's}:\bar\Cal M_{Ks}\longrightarrow\bar\Cal M_{K's}$. 

Now choose 
uniformising elements $t_K$ and $t_{K'}$ in $K$ and, resp., $K'$. 
Consider the corresponding standard 
generators $D^{(s)}_{an}$ (resp. $D^{\prime
(s)}_{an}$), where $a\in\Bbb Z(p)$ and $n\in\Bbb
Z\operatorname{mod}N_s$, 
 of $\bar\Cal M_{Ks}=\Cal M_{Ks}\hat\otimes _kk(p)$ 
(resp., 
$\bar\Cal M_{K's}=\Cal M_{Ks}\hat\otimes _kk(p)$). 
Here, as usual,  $k\simeq\Bbb F_{q_0}$ is the residue field of 
$K$, $q_0=p^{N_0}$, $N_s=N_0p^s$.  
Then 
$$g_{KK's}(D^{(s)}_{a0})=\sum\Sb b\in\Bbb Z(p),m\in
\Bbb Z\operatorname{mod}N_s \endSb 
\alpha _{abm}(g_{KK's})D^{\prime (s)}_{bm}$$
with $\alpha _{abm}(g_{KK's})\in k_s\subset k(p)$. 

For each $s\geqslant 0$, choose $n_s\in\Bbb Z\operatorname{mod}N_s$ 
such that $\alpha _{11n_s}(g_{KK's})\ne 0$: 
$n_s$ exists, because $g_{KK's}$ induces a    
$k(p)$-linear isomorphism of  
$\bar\Cal M_{Ks}\operatorname{mod}\bar\Cal M^{(2)}_{Ks}$ 
and $\bar\Cal M_{K's}\operatorname{mod}\bar\Cal
M^{(2)}_{K's}$. 

Let $\Fr (t_{K'})\in\Aut K'_{\ur }$ be such that 
$\Fr (t_{K'}):t_{K'}\mapsto t_{K'}$ and 
$\Fr (t_{K'})|_{k(p)}=\sigma $. Let 
$\xi \in\Iso ^0(K'_{\ur }, K_{\ur })$ be such that 
$\xi (t_{K'})=t_K$. 

For any $s\geqslant 0$, $\Fr (t_{K'})$ (resp.  
$\xi $) induces a continuous 
field isomorphism $K'_s\longrightarrow K'_s$ 
(resp.    
$K'_s\longrightarrow K_s$). It will be denoted by 
$\Fr (t_{K'})_s$ (resp. $\xi _s$).  
With notation from n.3, we introduce 
continuous group isomorphisms 
$$g^0_{KK's}=g_{KK's}\Fr (t_{K'})_s^{n_s*}:
\bar\Cal M_{Ks}\longrightarrow\bar\Cal M_{K's}.$$

Clearly, $h_s:=g^0_{KK's}\xi _s^*$ is induced 
by an automorphism of $\Gamma _{K_s}(p)$ which is compatible with 
the ramification filtration.
Notice also that, by proposition 2.1, if $a\in\Bbb Z(p)$, $n\in\Bbb
Z\operatorname{mod}N_s$ and  
$$h_{s}(D^{(s)}_{a0})=\sum\Sb b,m\endSb 
\alpha _{abm}(h_s)D^{(s)}_{bm},$$
then $\alpha _{a,b,m-n_s}(h_s)=\sigma ^{n_s}\alpha
_{abm}(g_{KK's})$. In particular,  
$\alpha _{110}(h_s)\ne 0$. 
Therefore, applying proposition 5.6,   
we obtain that for all $s\geqslant 0$, 
$$h_s\in\Aut _{\adm }\Cal M_{Ks}
\operatorname{mod}\Cal M_{K_s}^{(p^{N_s-2})},$$ 
 the residues $n_s\in\Bbb Z\operatorname{mod}N_s$ are unique, and
$n_{s+1}\operatorname{mod}N_s=n_s$. Here we use that  
$D^{(s+1)}_{an}\mapsto D^{(s)}_{an}$ under 
the natural morphism from $\bar\Cal M_{K,{s+1}}$ to $\bar\Cal
M_{Ks}$. 
Then 
$h_{KK}:=\{h_s\}_{s\geqslant 0}$ and 
$g^0_{KK'}:=\{g^0_{KK's}\}_{s\geqslant 0}$ 
are compatible systems and, by propositions 3.3 and 5.8,  
they are special admisible locally analytic 
systems. By proposition 3.4 there is an 
$\eta _{KK'}\in\Iso ^0(K,K')$ such that 
$g^0_{KK'\an }=\operatorname{d}(\eta _{KK'})\hat\otimes _kk(p)$. 
Notice also that if 
$\bar n_{KK'}:=\mathbin{\underset{s}\to\varprojlim }n_s\in
\mathbin{\underset{s}\to\varprojlim}\Bbb Z/N_s\Bbb Z$ 
then $g_{KK'}=g^0_{KK'}\Fr (t_{K'})^{-\bar n_{KK'}*}$, 
where $\Fr (t_{K'})^*=\{\Fr (t_{K'})_s\}_{s\geqslant 0}$ is 
the compatible system from n.3.5. 

Suppose $L$ is a finite Galois extension of $E$ in $E(p)$ 
containing $K$. Proceed similarly to obtain $L'\subset E(p)$ such that 
$g$ induces an isomorphism of $\Gamma _L(p)$ and $\Gamma _{L'}(p)$, the
corresponding 
compatible system $g_{LL'}=\{g_{LL's}\}_{s\geqslant 0}$ and the special
admissible locally analytic system 
$g^0_{LL'}=\{g^0_{LL's}\}_{s\geqslant 0}$, where 
$g_{LL'}=g^0_{LL'}\Fr (t_{L'})^{-\bar n_{LL'}*}$, together with the 
corresponding $\eta _{LL'}\in\Iso ^0(L,L')$ such that 
$g^0_{LL'\an }=\operatorname{d}(\eta _{LL'})\hat\otimes _{k_L}k_L(p)$.  
Here $k_L$ is the residue field of $L$, 
$k_L\simeq\Bbb F_{p^{M_0}}$ and $\bar n_{LL'}\in\varprojlim\Bbb
Z/p^{M_0p^s}\Bbb Z$. Notice that all these maps depend on some 
choice of uniformising elements 
$t_L$ and $t_{L'}$ in, respectively, $L$ and $L'$. 

The systems  $g_{LL'}$ and $g_{KK'}$ are comparable because 
both come from the group isomorphisms 
$\Gamma _L(p)\longrightarrow\Gamma _{L'}(p)$ and 
$\Gamma _K(p)\longrightarrow\Gamma _{K'}(p)$ which are induced by
$g$. If $I_{L/K}$ is the inertia subgroup of  
$\Gal (L/K)$ then there is a natural group embedding 
$I_{L/K}\subset\Aut ^0(L)\subset\Aut ^0(L_{\ur })$. 
Similarly, we have a group embedding for the inertia subgroup 
$I_{L'/K'}$ of $\Gal (L'/K')$ into $\Aut ^0(L')$. 

Let $\kappa :I_{L/K}\longrightarrow I_{L'/K'}$ be the group isomorphism 
induced by $g$. Then 
$\tau ^*g_{LL's}=g_{LL's}\kappa (\tau )^*$, 
 for any 
$\tau\in I_{L/K}$ and any $s\geqslant 0$.  
This implies that 
$$\tau ^*g_{LL'\ur }=g_{LL'\ur }\kappa (\tau )^*,$$ 
i.e. condition C from n.3.7 holds in this case.

Let $\mu _{KK'}=\eta _{KK'}\Fr (t_{K'})^{-\bar n_{KK'}}\in\Iso (K,K')$ and 
$\mu _{LL'}=\eta _{LL'}\Fr (t_{K'})^{-\bar n_{LL'}}\in\Iso (L,L')$. 

\proclaim{Proposition 6.1} With the above notation: 
\newline 
{\rm a)} $\mu _{LL'}|_K=\mu _{KK'}$; 
\newline 
{\rm b)} for any $\tau\in I_{L/K}$, $\tau\mu _{LL'}=\mu _{LL'}\kappa
(\tau )$. 
\endproclaim 

\demo{Proof} Let 
$\alpha =\Fr (t_{L'})^{\bar n_{LL'}}$. 
Consider $K'_{\ur }$ as a subfield in $L'_{\ur }$ and 
set $K''_{\ur }=\alpha (K'_{\ur })\subset L'_{\ur
}$. Then $K''_{\ur }$ is the maximal unramified $p$-extension 
of the complete discrete valuation field 
$K'':=\alpha (K')\subset E(p)$ in $E(p)$.  

Let 
$\beta =\alpha |_{K'_{\ur }}$. Consider the following 
commutative diagramm

$$
\CD 
\bar\Cal M_{L\ur }  @>{g_{LL'\ur }}>> \bar\Cal M_{L'\ur } 
@>{\alpha ^*_{L'L'\ur }}>>   \bar\Cal M_{L'{\ur }} \\ 
@VVV       @VVV            @VVV \\ 
\bar\Cal M_{K\ur }  @>{g_{KK'\ur }}>>   \bar\Cal M_{K'\ur } 
@>{\beta ^{*}_{K'K''\ur }}>> \bar\Cal M_{K''\ur }  
\endCD 
$$
where the vertical arrows come from natural embeddings of the
corresponding Galois groups. 

The systems $g^0_{LL'}=g_{LL'}\alpha ^*_{L'L'}$ 
and $f_{KK''}:=g_{KK'}\beta ^{*}_{K'K''}$ 
are comparable, because they come from the compatible 
group isomorphisms $\Gamma _L(p)\longrightarrow\Gamma _{L'}(p)$ 
and 
\linebreak 
$\Gamma _K(p)\mathbin{\overset{f}\to\longrightarrow}\Gamma _{K''}(p)$. 
In this situation, condiditon {\bf C} is automatically 
satisfied and, by 
proposition 3.5, the admissibility of 
$g^0_{LL'}$ implies the admissibility of 
$f_{KK''}$. Because the group homomorphism $f$ is 
compatible with ramification filtrations, we can apply 
the results of section 5 to deduce that 
$f_{KK''}$ is special admissible locally analytic 
and that there is an $\eta ^1_{KK''}\in\Iso ^0(K,K'')$ 
such that $f_{KK''\an }=
\operatorname{d}(\eta ^1_{KK''})\hat\otimes _kk(p)$ and 
$\eta _{LL'}|_K=\eta ^1_{KK''}$. 

Consider $\psi :=\eta ^{-1}_{KK'}\eta _{LL'}|_K\in\Iso ^0(K',K'')$. 
Then 

$$\psi _{\an }=\eta ^{-1}_{KK'\an }\eta ^1_{KK''\an }=
(g^0_{KK'\an })^{-1}(g_{KK'}\beta ^{*}_{K'K''})_{KK''\an }$$
$$=\left ({g^0_{KK'}}^{-1}g_{KK'}\beta ^{*}_{K'K''}\right
)_{K'K''\an }=
\left (\Fr (t_{K'})^{-\bar n_{KK'}}\beta\right )_{\an }.$$
Therefore by proposition 2.7, 
$$\eta ^{-1}_{KK'}\eta _{LL'}|_{K}=
\Fr (t_{K'})^{-\bar n_{KK'}}\Fr (t_{L'})^{\bar n_{LL'}}|_K$$
or $\mu _{LL'}|_K=\mu _{KK'}$. 

Part a) of our proposition is proved.

Consider the inertia subgroups 
$I_{L/K}\subset\Gal (L_{\ur }/K_{\ur })$, 
$I_{L'/K'}\subset\Gal (L'_{\ur }/K'_{\ur })$ and 
$I_{L'/K''}\subset\Gal (L'_{\ur }/K''_{\ur })$. As it was noticed 
earlier, the correspondence 
$$\tau ^*\mapsto \tau ^{\prime *}=
g^{-1}_{LL'\ur }\tau ^*g_{LL'\ur }$$
induces a group isomorphism 
$\kappa :I_{L/K}\longrightarrow  I_{L'/K'}$ such that 
$\kappa (\tau )=\tau '$. 

We use the correspondence 
$$\alpha ^*:\tau '\mapsto\tau ''=\alpha ^{-1}\tau '\alpha $$
to define the  group isomorphism $\kappa _{\alpha }:I_{L'/K'}\longrightarrow
I_{L'/K''}$ such that $\kappa _{\alpha }(\tau ')=\tau ''$. With this
notation 
we have the following equality of compatible systems 
$$\tau ^*_{LL}g^0_{LL'}=g^0_{LL'}{\tau ''_{L'L'}}^*,$$
where as earlier, $g^0_{LL'}=g_{LL'}\alpha ^*_{L'L'}$. 

Therefore, the equality 
$(\tau \eta _{LL'})_{\an }=(\tau ^*_{LL}g^0_{LL'})_{\an }
=(g^0_{LL'}{\tau ''_{L'L'}}^*)_{\an }=(\eta _{LL'}\tau '')_{\an }$ 
together with proposition 2.7 and the definition of $\tau ''$ imply  that 
$\tau \eta _{LL'}=\eta _{LL'}\tau ''=\eta _{LL'}\alpha ^{-1}\tau
'\alpha $, 
i.e. $\tau\mu _{LL'}=\mu _{LL'}\tau '$. 

The proposition is proved. 
\enddemo 

Let $\mu :=\mathbin{\underset{\rightarrow}\to\lim} \mu _{KK'}:E(p)\longrightarrow E(p)$. 
Clearly, it is a continuous field isomorphism and $\mu (E)=E'$. 

\proclaim{Proposition 6.2} $\mu ^*=g$.
\endproclaim 

\demo{Proof} As earlier, let $K$ and $K'$ be Galois extensions of 
$E$ and $E'$, respectively, such that 
$g(\Gamma _K(p))=\Gamma _{K'}(p)$.

By part b) of the above proposition 6.1, 
the correspondences 
$\mu ^*:\tau\mapsto\mu ^{-1}\tau\mu $ and 
$g:\tau\mapsto g(\tau )$ induce the same isomorphism of the inertia
subgroups 
$I_K(p)\longrightarrow I_{K'}(p)$. Consider the induced 
isomorphism $I_K(p)^{\ab }\longrightarrow I_{K'}(p)^{\ab }$. 
With respect to 
the identifications of class field theory 
$I_K(p)^{\ab }=U_K$ and 
$I_{K'}(p)^{\ab }=U_{K'}$, where $U_K$ and $U_{K'}$ are groups 
of principal units in $K$ and $K'$, respectively, this homomorphism 
is induced by the restriction of the field isomorphism 
$\mu _{KK'}$ on $U_K$. 
In addition, $\mu _{KK'}$ transforms the natural action 
of any $\tau \in\Gamma _E(p)$ on $U_K$ into the natural action 
of $g(\tau )\in\Gamma _{E'}(p)$ on $U_{K'}$. Therefore, the 
two field automorphisms 
$\mu ^{-1}\tau\mu |_{K'}$ and $g(\tau )|_{K'}$ of $K'$ 
become equal after restricting on 
$U_{K'}$. This implies that they 
coincide on the whole field $K'$, i.e.  
$\mu ^{-1}\tau\mu \equiv g(\tau )\operatorname{mod}\Gamma _{K'}(p)$, 
for any $\tau\in\Gamma _E(p)$.   
Because  
$K$ is an arbitrary Galois extension of $E$ in $E(p)$ this implies that 
$g=\mu ^*$. 

So, proposition 6.2 together with the characteristic $p$ case of the Main
Theorem are completely proved. 
\enddemo

\subhead 7. Proof of the main theorem --- the mixed characteristic case 
\endsubhead 
\medskip 

In this section $\char E=0$. Clearly, this implies that $\char E'=0$. 
\medskip 

7.1. Following the paper [Wtb] introduce the categories  
$\Psi $, $\widetilde{\Psi }$ and the functor $\Phi
:\Psi\longrightarrow\widetilde{\Psi }$. 

The objects of $\Psi $ are the 
field extensions $L/K$, where $[K:\Bbb Q_p]<\infty $, $L$ is 
an infinite Galois extension of $K$ in a fixed maximal $p$-extension
$K(p)$ of $K$ and $\Gamma _{L/K}=\Gal (L/K)$ is a $p$-adic Lie group. 
A morphism from $L/K$ to an object $L'/K'$ in $\Psi $ is 
a continuous field embedding $f:L\longrightarrow L'$ such
that $[L':f(L)]<\infty $ and $f|_K$ is a field isomorphism of $K$ and 
$K'$. 

The objects of $\widetilde{\Psi }$ are couples $(\Cal K,G)$ where 
$\Cal K$ is a complete discrete valuation field 
of characteristic $p$ with finite residue field and $G$ is a closed
subgroup of the group of all continuous automorphisms of $\Cal K$. 
In addition, with respect to the induced topology $G$, is a compact 
finite dimensional $p$-adic Lie group. A morphism from 
$(\Cal K,G)$ to an object $(\Cal K',G')$ in $\widetilde{\Psi }$
is a closed field embedding $f:\Cal K\longrightarrow\Cal K'$ such that 
$\Cal K'$ is a finite separable extension of $f(\Cal K)$. In addition, 
$f(\Cal K)$ is $G'$-invariant and the corrspondence 
$\tau\mapsto \tau|_{f(\Cal K)}$ induces a group epimorphism 
from $G'$ to $G$. 

Let $X$ be the Fontaine-Wintenberger field-of-norm functor,
cf. [Wi2]. 
Then the correspondence $L/K\mapsto (X(L),G_{L/K})$, 
where $G_{L/K}=\{X(\tau )\ |\ \tau\in\Gamma _{L/K}\}$, 
induces the functor $\Phi :\Psi\longrightarrow\widetilde{\Psi }$. 

One of main results in [Wi1] states that the functor $\Phi $ is 
fully faithful. 
\medskip 

7.2. Let $\{E_{\alpha }/E, i_{\alpha\beta }\}_{\Cal I}$ be an
     inductive system of objects in the category $\Psi $. 
From now on $\Cal I$ is a set of indices $\alpha $ 
with a suitable partial ordering.  The connecting morphisms 
$i_{\alpha\beta }\in\Hom _{\Psi }
(E_\alpha , E_{\beta })$ are the natural field embeddings 
defined for suitable couples $\alpha ,\beta\in\Cal I$. 
We can choose this inductive system to be large enough 
to satisfy the requirement $\mathbin{\underset\rightarrow\to\lim}E_{\alpha }
=E(p)$. 

By applying the functor $\Phi $, we obtain the inductive system 
$\{(\Cal E_{\alpha }, G_{\alpha }), \tilde i_{\alpha\beta }\}_{\Cal I}$ 
in the category $\widetilde{\Psi }$, where  
$(\Cal E_{\alpha },
G_{\alpha })=\Phi (E_{\alpha }/E)$ and $\tilde\imath _{\alpha \beta }=\Phi
(i_{\alpha\beta })$, for all $\alpha\in\Cal I$.  Then 
$\mathbin{\underset\rightarrow\to\lim}\Cal E_{\alpha }=\Cal E(p)$ is a
maximal $p$-extension for each field $\Cal E_{\alpha }$,
$\alpha\in\Cal I$.          

Notice that the field embeddings 
$\tilde\imath _{\alpha\beta }$ induce group epimorphisms 
$\tilde\jmath _{\alpha\beta }:G_{\beta }\longrightarrow G_{\alpha }$ 
with corresponding projective system  
$\{G_{\alpha },\tilde\jmath _{\alpha\beta }\}_{\Cal I}$ such that 
$\varprojlim G_{\alpha }$ is identified via the functor $X$ with 
$\Gamma _E(p)$. For any $\alpha\in\Cal I$, 
we then have the identifications 
$\Gamma _{E_{\alpha }}(p)=\Gamma _{\Cal E_{\alpha }}(p)$. These
identifications are compatible with the ramification filtrations. 
This means that one can define the Herbrand function 
$\varphi _{\alpha }$ for the infinite extension $E_{\alpha }/E$ 
as the limit of Herbrand functions of all finite
subextensions in $E_{\alpha }$ over $E$ and 
$$\Gamma _E(p)^{(v)}\cap\Gamma _{E_{\alpha }}(p)=
\Gamma _{\Cal E_{\alpha }}(p)^{(\varphi _{\alpha }(v))},$$
for all $v\geqslant 0$. 
\medskip 

7.3. Consider the group isomorphism 
$g:\Gamma _E(p)\longrightarrow\Gamma _{E'}(p)$ from the statement of
the Theorem. For $\alpha\in\Cal I$, 
let $E'_{\alpha }\subset E'(p)$ be such that 
$g(\Gamma _{E_{\alpha }}(p))=\Gamma _{E'_{\alpha }}(p)$. Then we have 
the corresponding injective system 
$\{E'_{\alpha },i'_{\alpha\beta }\}_{\Cal I}$ and 
$\mathbin{\underset\rightarrow\to\lim}E'_{\alpha }=E'(p)$. 

Clearly, for any $\alpha\in\Cal I$, 
\medskip 

$\bullet $\ \ $E'_{\alpha }/E'$ is an object of $\Psi $;
\medskip 

$\bullet $\ \ $\bar g_{\alpha }:=g_{\alpha }\operatorname{mod}
\Gamma _{E_{\alpha }}(p):\Gamma _{E_{\alpha }/E}\longrightarrow\Gamma
_{E'_{\alpha }/E'}$ is a group 
isomorphism which is compatible with the ramification filtrations; 
in particular, this implies that the Herbrand functions 
for the infinite extensions $E_{\alpha }/E$ and 
$E'_{\alpha }/E'$ are equal;
\medskip 

$\bullet $\ \ for any $v\geqslant 0$, 
$g_{\alpha }:=g|_{\Gamma _{E_{\alpha }}(p)}$ induces a continuous 
group isomorphism 
of $\Gamma _E(p)^{(v)}\cap\Gamma _{E_{\alpha }}(p)$ and 
$\Gamma _{E'}(p)^{(v)}\cap\Gamma _{E'_{\alpha }}(p)$.
\medskip 

For $\alpha\in\Cal I$, set $\Phi (E'_{\alpha }/E')=
(\Cal E'_{\alpha },G'_{\alpha })$ and 
$\Phi (i'_{\alpha\beta })=\tilde \imath '_{\alpha\beta }$. 
Then $\{(\Cal E'_{\alpha },G'_{\alpha }), \tilde\imath '_{\alpha\beta
}\}_{\Cal I}$ is an inductive system, 
$\mathbin{\underset\rightarrow\to\lim}\Cal E'_{\alpha }:=\Cal E'(p)$ 
is a maximal $p$-extension for each $\Cal E'_{\alpha }$. As earlier, 
we obtain the projective system 
$\{G'_{\alpha }, \tilde\jmath '_{\alpha\beta }\}_{\Cal I}$ and 
the field-of-norms functor allows us to identify 
the topological groups 
$\varprojlim G'_{\alpha }$ and $\Gamma _{E'}(p)$. 
Therefore, for any $\alpha\in\Cal I$, we have an identification of
the 
groups $\Gamma _{E'_{\alpha }}(p)$ and $\Gamma _{\Cal E'_{\alpha
}}(p)$. 

This implies that for all $\alpha\in\Cal I$, we have the following 
isomorphisms of topological groups:
\medskip 

$\bullet $ \ \ $\tilde g_{\alpha }:=X(g_{\alpha }):\Gamma _{\Cal E_{\alpha
}}(p)\longrightarrow\Gamma _{\Cal E'_{\alpha }}(p)$ 
such that, for any $v\geqslant 0$, 
$\tilde g_{\alpha }(\Gamma _{\Cal E_{\alpha }}(p)^{(v)})=
\Gamma _{\Cal E'_{\alpha }}(p)^{(v)}$;
\medskip 

$\bullet $\ \ $X(\bar g_{\alpha }):G_{\alpha }\longrightarrow
G'_{\alpha }$ which maps the projective system 
$\{G_{\alpha },\tilde\jmath _{\alpha\beta }\}_{\Cal I}$  
to the projective system 
$\{G'_{\alpha },\tilde\jmath '_{\alpha\beta }\}_{\Cal I}$. 
\medskip 

7.4. By the characteristic $p$ case of the Main Theorem 
for all $\alpha\in\Cal I$, 
there are continuous field isomorphisms 
$\tilde\mu _{\alpha }:\Cal E_{\alpha }\longrightarrow \Cal E'_{\alpha
}$ 
such that 
\medskip 

$\bullet $\ \ $\{\tilde\mu _{\alpha }\}_{\alpha\in\Cal I}$ maps the
inductive system $\{\Cal E_{\alpha }, \tilde\imath _{\alpha\beta }\}_{\Cal
I}$ to the inductive system 
$\{\Cal E_{\alpha }',\tilde\imath '_{\alpha\beta }\}_{\Cal I}$;
\medskip 

$\bullet $\ \ $X(\bar g_{\alpha })$ is induced by $\tilde\mu _{\alpha
}$, i.e. if $\tau\in G_{\alpha }$ and $\tau '=X(\bar g_{\alpha })\in
G'_{\alpha }$ then $\tau\tilde\mu _{\alpha }=\tilde\mu _{\alpha }\tau
'$.
\medskip 

Because $\Phi $ is fully faithful for all $\alpha\in\Cal I$, there is
a 
$\mu _{\alpha }\in\Hom _{\Psi }(E_{\alpha }/E, E'_{\alpha }/E')$ 
such that 
\medskip 

$\bullet $\ \ $\{\mu _{\alpha }\}_{\alpha\in\Cal I}$ transforms the
inductive system 
$\{E_{\alpha }/E, i_{\alpha\beta }\}_{\Cal I}$ into the inductive
system $\{E'_{\alpha }/E', i'_{\alpha\beta }\}_{\Cal I}$;
\medskip 

$\bullet $\ \ if $\tau\in\Gamma _{E_{\alpha }/E}$ and $\tau '=\bar
g_{\alpha }(\tau )\in\Gamma _{E'_{\alpha }/E'}$ then 
$\tau\mu _{\alpha }=\mu _{\alpha }\tau '$.
\medskip 

Therefore, $\mu :=\mathbin{\underset\rightarrow\to\lim }\mu _{\alpha
}$ 
is a continuous field isomorphism from $E(p)$ to $E'(p)$ such that 
$\tau\mu =\mu g(\tau )$, i.e. $g(\tau )=\mu ^{-1}\tau\mu $, 
for $\tau\in\varprojlim\Gamma _{E_{\alpha }/E}=\Gamma _E(p)$ and 
$g(\tau )\in\varprojlim\Gamma _{E'_{\alpha }/E'}=\Gamma _{E'}(p)$.

The Main Theorem is completely proved.

\enddocument

\enddocument